\documentclass[12pt,reqno]{amsart}

\usepackage{amsthm,amsfonts,amssymb,euscript}
\numberwithin{equation}{section}
\begin{document}

\title[Global well-posedness of the Benjamin--Ono equation]
{Global well-posedness of the Benjamin--Ono equation in low-regularity spaces}
\author{Alexandru D. Ionescu}
\address{University of Wisconsin--Madison}
\email{ionescu@math.wisc.edu}
\author{Carlos E. Kenig}
\address{University of Chicago}
\email{cek@math.uchicago.edu}
\thanks{The first author was
supported in part by an NSF grant, an Alfred P. Sloan
research fellowship, and a David and Lucile Packard fellowship.
The second author was supported in part by an NSF grant.}
\begin{abstract}
We prove that the Benjamin--Ono initial value problem is glo\-ba\-lly well-posed in the Banach spaces $H^\sigma_r(\mathbb{R})$, $\sigma\geq 0$, of real-valued Sobolev functions.
\end{abstract}
\maketitle

\section{Introduction}\label{section1}

In this paper we consider the Benjamin--Ono initial value problem
\begin{equation}\label{eq-1}
\begin{cases}
\partial_tu+\mathcal{H}\partial_x^2u+\partial_x(u^2/2)=0;\\
u(0)=\phi,
\end{cases}
\end{equation}
where $\mathcal{H}$ is the Hilbert transform operator defined (on the spaces $C(\mathbb{R}:H^\sigma)$, $\sigma\in\mathbb{R}$) by
the  Fourier multiplier $-i\,\mathrm{sgn}(\xi)$. The Benjamin--Ono equation is a model for one-dimensional long waves in deep stratified fluids (\cite{Be}  and \cite{On}), and is completely integrable. The  initial value problem for this equation has been studied extensively for data in the Sobolev spaces $H^\sigma_r(\mathbb{R})$, $\sigma\geq 0$.\footnote{In this paper $H^\sigma_r=H^\sigma_r(\mathbb{R})$ denotes the space of {\it{real-valued}} functions $\phi$ with the usual norm $||\phi||_{H^\sigma_r}=||\phi||_{H^\sigma}=||(1+|\xi|^2)^{\sigma/2}\widehat{\phi}(\xi)||_{L^2_\xi}$. All the other Banach spaces of functions, such as $L^2$, $H^\sigma$, $\widetilde{H}^\sigma$, $F^\sigma$, $N^\sigma$ etc, are defined as spaces of {\it{complex-valued}} functions.} It is known that the Benjamin--Ono initial value problem has weak solutions in $H^0_r(\mathbb{R})$, $H^{1/2}_r(\mathbb{R})$, and $H^1_r(\mathbb{R})$ (see \cite{GiVe}, \cite{To}, and \cite{Sa}), and is globally well-posed in $H^\sigma_r(\mathbb{R})$, $\sigma\geq 1$ (see \cite{Tao1}, as well as \cite{Io}, \cite{Po}, \cite{KoTz2}, and \cite{KeKo} for earlier local and global well-posedness results in higher regularity spaces). In this paper we prove that the Benjamin--Ono initial value problem is globally well-posed in $H^\sigma_r(\mathbb{R})$, $\sigma\geq 0$.

Let $H^\infty_r(\mathbb{R})=\cap_{\sigma=0}^\infty H^\sigma_r(\mathbb{R})$ with the induced metric. Let $S^\infty:H^\infty_r(\mathbb{R})\to C(\mathbb{R}:H^\infty_r(\mathbb{R}))$ denote the (nonlinear) mapping that associates to any data $\phi\in H^\infty_r$ the corresponding classical solution $u\in C(\mathbb{R}:H^\infty_r)$ of the initial value problem \eqref{eq-1}. We will use the $L^2$ conservation law: if $\phi\in H^\infty_r$ and $u=S^\infty(\phi)$ then
\begin{equation}\label{conserve}
\int_{\mathbb{R}}u(x,t)^2\,dx=\int_{\mathbb{R}}\phi(x)^2\,dx\,\text{ for any }t\in\mathbb{R}.
\end{equation}
For $T>0$ let $S^\infty_T:H^\infty_r(\mathbb{R})\to C([-T,T]:H^\infty_r(\mathbb{R}))$ denote the restriction of the  mapping $S^\infty$ to the time interval $[-T,T]$.

\newtheorem{Main1}{Theorem}[section]
\begin{Main1}\label{Main1}
(a) Assume $T>0$. Then the mapping $S_T^\infty:H^\infty_r\to C([-T,T]:H^\infty_r)$ extends uniquely to a continuous mapping $S^0_T:H^0_r\to C([-T,T]:H^0_r)$ and
\begin{equation*}
||S^0_T(\phi)(t)||_{H^0_r}=||\phi||_{H^0_r}\,\text{ for any }\,t\in[-T,T],\,\phi\in H^0_r.
\end{equation*}
The function $S^0_T(\phi)$ solves the initial value problem \eqref{eq-1} in $\mathcal{S}'(\mathbb{R}\times(-T,T))$ for any $\phi\in H^0_r$. 

(b) In addition, for any $\sigma\geq 0$, $S^0_T(H^\sigma_r)\subseteq C([-T,T]:H^\sigma_r)$,
\begin{equation*}
||S^0_T(\phi)||_{C([-T,T]:H^\sigma_r)}\leq C(T,\sigma,||\phi||_{H^\sigma_r}),
\end{equation*}
and the mapping $S^\sigma_T={S^0_T}|_{H^\sigma_r}:H^\sigma_r\to C([-T,T]:H^\sigma_r)$ is continuous.
\end{Main1}

We mention that the flow map $\phi\to S^\sigma_T(\phi)$ fails to
be uniformly continuous on bounded sets in $H^\sigma_r$ for any
$\sigma>0$, see \cite{KoTz}. In a forthcoming we will consider well-posedness theorems for complex-valued data, as well as (local) ill-posedness theorems for data in $H^\sigma_r$, $\sigma<0$. 

We discuss now some of the ingredients in the proof of Theorem \ref{Main1}. The main obstruction to simply using a fixed-point argument in some
$X^{s,b}$ space (in a way similar to the case of the KdV equation,
see J. Bourgain \cite{Bo}) is the lack of control of the
interaction between very high and very low frequencies of
solutions (cf. \cite{MoSaTz} and \cite{KoTz}). Following \cite{Tao1}, we first construct
a gauge transformation that weakens this interaction, in the sense
that we will be able to assume that low frequency functions
have some additional structure (see the space $Z_0$ defined in
section \ref{section2}). Even with this low-frequency assumption,
the use of standard $X^{s,b}$ spaces for high-frequency
functions (i.e. spaces defined by suitably weighted norms in the
frequency space) seems to lead inevitably to logarithmic
divergences in the modulation variable (see \cite{KeCoSt}). To avoid these logarithmic
divergences we work with high-frequency spaces that have two
components: an $X^{s,b}$-type component in the frequency space and
a normalized $L^1_xL^2_t$ component in the physical space. This
type of spaces have been used in the context of wave maps by D.
Tataru \cite{Tat1} (and refined versions in \cite{Tat2},
\cite{Tao2}, and \cite{Tao3}); we remark that for the physical
space component we use a suitable normalization of the local
smoothing space $L^1_xL^2_t$ instead of the energy space
$L^1_tL^2_x$. Then we prove suitable linear and bilinear estimates
in these spaces, and conclude the proof of Theorem \ref{Main1}
using a fixed-point argument.

The rest of the paper is organized as follows: in section
\ref{gauge} we construct our gauge transformation and reduce
solving the initial-value problem \eqref{eq-1} to solving three
easier initial-value problems. The point of this reduction is
that the initial datum of the resulting three initial-value problems
have some special structure at very low frequencies (see the spaces $\widetilde{H}^\sigma$ defined in \eqref{def4}). In sections
\ref{section2} and \ref{prop} we construct our main Banach spaces and prove
some of their elementary properties. In section \ref{linear} we prove
several linear estimates using these Banach spaces. In sections
\ref{bilin}, \ref{bilinear} and \ref{bilinear2} we prove our main
bilinear estimates. In section \ref{mult} we prove several bounds
for operators defined by multiplication with certain smooth
bounded functions (such estimates are delicate in the context of $X^{s,b}$ spaces). Finally, in section \ref{proof} we combine all
these estimates and a fixed-point argument to complete the proof
of Theorem \ref{Main1}.

We would like to thank S. Herr and H. Koch for useful discussions in the early stages of this work. 

\section{The gauge transformation}\label{gauge}

The first step is to construct a gauge transformation to weaken significantly the contribution coming from the low frequencies of the data. Assume $\phi\in H^\infty_r$ and $u=S^\infty(\phi)\in C(\mathbb{R}:H^\infty_r)$. On $L^2(\mathbb{R})$ we define the operators
\begin{equation*}
\begin{split}
&P_{\mathrm{low}}\,\text{ defined by the Fourier multiplier }\,\xi\to\mathbf{1}_{[-2^{10},2^{10}]}(\xi);\\
&P_{\pm\mathrm{high}}\,\text{ defined by the Fourier multiplier }\,\xi\to\mathbf{1}_{[2^{10},\infty)}(\pm\xi);\\
&P_{\pm}\,\text{ defined by the Fourier multiplier }\,\xi\to\mathbf{1}_{[0,\infty)}(\pm\xi).
\end{split}
\end{equation*}
Let $\phi_0=P_{\mathrm{low}}\phi\in H^\infty_r$, $u_0=S^\infty(\phi_0)$, $\widetilde{u}=u-u_0$. Since $||\phi_0||_{H^\sigma_r}\leq C_\sigma||\phi||_{L^2}$ for any $\sigma\geq 0$,
\begin{equation}\label{fg0}
\sup_{t\in[-2,2]}||\partial_t^{\sigma_1}\partial_x^{\sigma_2}
u_0(.,t)||_{L^2_x}\leq
C_{\sigma_1,\sigma_2}||\phi||_{L^2},\,\,\sigma_1,\sigma_2\in[0,\infty)\cap\mathbb{Z}.
\end{equation}

Using the equation \eqref{eq-1},
\begin{equation}\label{fg1}
\begin{cases}
&\partial_t\widetilde{u}+\mathcal{H}\partial_x^2\widetilde{u}+\partial_x(u_0\cdot \widetilde{u})+\partial_x(\widetilde{u}^2/2)=0;\\
&\widetilde{u}(0)=P_{+\mathrm{high}}\phi+P_{-\mathrm{high}}\phi.
\end{cases}
\end{equation}
We apply $P_{+\mathrm{high}}$, $P_{-\mathrm{high}}$, and $P_{\mathrm{low}}$ to \eqref{fg1} to obtain
\begin{equation}\label{fg2}
\begin{cases}
&\partial_t(P_{\pm\mathrm{high}}\widetilde{u})\mp i\cdot\partial_x^2(P_{\pm\mathrm{high}}\widetilde{u})+P_{\pm\mathrm{high}}\partial_x(u_0\cdot \widetilde{u})+P_{\pm\mathrm{high}}\partial_x(\widetilde{u}^2/2)=0;\\
&(P_{\pm\mathrm{high}}\widetilde{u})(0)=P_{\pm\mathrm{high}}\phi,
\end{cases}
\end{equation}
and
\begin{equation}\label{fg3}
\begin{cases}
&\partial_t(P_{\mathrm{low}}\widetilde{u})+\mathcal{H}\partial_x^2(P_{\mathrm{low}}\widetilde{u})+P_{\mathrm{low}}\partial_x(u_0\cdot \widetilde{u})+P_{\mathrm{low}}\partial_x(\widetilde{u}^2/2)=0;\\
&(P_{\mathrm{low}}\widetilde{u})(0)=0.
\end{cases}
\end{equation}

We now let
\begin{equation}\label{fg4}
\begin{cases}
&P_{+\mathrm{high}}\widetilde{u}=e^{-iU_0}w_+;\\
&P_{-\mathrm{high}}\widetilde{u}=e^{iU_0}w_-;\\
&P_{\mathrm{low}}\widetilde{u}=w_0,
\end{cases}
\end{equation}
where $U_0$ is a suitable gauge that depends only on $u_0$. As in \cite{Tao1}, we define first $U(0,t)$ on the time axis $x=0$ by the formula
\begin{equation}\label{fg5}
\partial_tU_0(0,t)+\frac{1}{2}\mathcal{H}\partial_xu_0(0,t)+\frac{1}{4}u_0^2(0,t)=0,\,\,\,U_0(0,0)=0,
\end{equation}
and then we construct $U_0(x,t)$ using the formula
\begin{equation}\label{fg6}
\partial_xU_0(x,t)=\frac{1}{2}u_0(x,t).
\end{equation}
It is important to notice that $U_0$ is real-valued, since $\phi_0$ and $u_0$ are both real-valued. Using the equation \eqref{eq-1} for $u_0=S^\infty(\phi_0)$ and \eqref{fg6}, we have
\begin{equation*}
\partial_x[\partial_tU_0+\mathcal{H}\partial_x^2U_0+(\partial_xU_0)^2]=0\,\,\text{ on }\,\,\mathbb{R}\times\mathbb{R}.
\end{equation*}
Using \eqref{fg5} and \eqref{fg6} it follows that
\begin{equation}\label{fg8}
\partial_tU_0=-\frac{1}{2}\mathcal{H}\partial_xu_0-\frac{1}{4}u_0^2\,\,\text{ on }\,\,\mathbb{R}\times\mathbb{R}.
\end{equation}
In particular, in view of \eqref{fg6} and \eqref{fg8}, $U_0\in C^\infty(\mathbb{R}\times\mathbb{R})$. In fact, it follows from \eqref{fg0}, \eqref{fg6}, and \eqref{fg8} that for any integers $\sigma_1,\sigma_2\geq 0$, $(\sigma_1,\sigma_2)\neq (0,0)$,
\begin{equation}\label{fg9}
\sup_{t\in[-2,2]}||\partial_t^{\sigma_1}\partial_x^{\sigma_2}
U_0(.,t)||_{L^2_x}\leq C_{\sigma_1,\sigma_2}||\phi||_{L^2}.
\end{equation}

We substitute now the formulas $P_{+\mathrm{high}}\widetilde{u}=e^{-iU_0}w_+$ and $\widetilde{u}=e^{-iU_0}w_++e^{iU_0}w_-+w_0$ in the equation \eqref{fg2} for $P_{+\mathrm{high}}\widetilde{u}$; the term $P_{+\mathrm{high}}(u_0e^{-iU_0}\partial_xw_+)$ cancels (using \eqref{fg6}), and the result is
\begin{equation}\label{fg20}
\begin{cases}
&(\partial_t+\mathcal{H}\partial_x^2)w_+=E_+(w_+,w_-,w_0);\\
&w_+(0)=e^{iU_0(.,0)}P_{+\mathrm{high}}\phi,
\end{cases}
\end{equation}
where
\begin{equation*}
\begin{split}
E_+(w_+,&w_-,w_0)=-e^{iU_0}P_{+\mathrm{high}}[\partial_x(e^{-iU_0}w_++e^{iU_0}w_-+w_0)^2/2]\\
&-e^{iU_0}P_{+\mathrm{high}}[\partial_x[u_0(e^{iU_0}w_-+w_0)]]\\
&+e^{iU_0}(P_{-\mathrm{high}}+P_{\mathrm{low}})(u_0e^{-iU_0}\partial_xw_+)+2iP_-(\partial_x^2w_+)\\
&-e^{iU_0}P_{+\mathrm{high}}[\partial_x(u_0e^{-iU_0})\cdot w_+]+i(\partial_tU_0-i\partial_x^2U_0-(\partial_xU_0)^2)\cdot w_+.
\end{split}
\end{equation*}
Since $w_+=e^{iU_0}P_{+\mathrm{high}}(e^{-iU_0}w_+)$, $w_-=e^{-iU_0}P_{-\mathrm{high}}(e^{iU_0}w_-)$, and $w_0=P_{\mathrm{low}}(w_0)$ (see \eqref{fg4}), we use \eqref{fg6} and \eqref{fg8} to rewrite $E_+(w_+,w_-,w_0)$ in the form
\begin{equation}\label{fg10}
\begin{split}
E_+(w_+,&w_-,w_0)=-e^{iU_0}P_{+\mathrm{high}}[\partial_x(e^{-iU_0}w_++e^{iU_0}w_-+w_0)^2/2]\\
&-e^{iU_0}P_{+\mathrm{high}}[\partial_x[u_0\cdot P_{-\mathrm{high}}(e^{iU_0}w_-)+u_0\cdot P_{\mathrm{low}}(w_0)]]\\
&+e^{iU_0}(P_{-\mathrm{high}}+P_{\mathrm{low}})[\partial_x(u_0\cdot P_{+\mathrm{high}}(e^{-iU_0}w_+))]\\
&+2iP_-[\partial_x^2(e^{iU_0}P_{+\mathrm{high}}(e^{-iU_0}w_+))]\\
&-P_+\partial_xu_0\cdot w_+.
\end{split}
\end{equation}

A similar computation using the equation \eqref{fg2} for $P_{-\mathrm{high}}\widetilde{u}$ gives
\begin{equation}\label{fg21}
\begin{cases}
&(\partial_t+\mathcal{H}\partial_x^2)w_-=E_-(w_+,w_-,w_0);\\
&w_-(0)=e^{-iU_0(.,0)}P_{-\mathrm{high}}\phi,
\end{cases}
\end{equation}
where
\begin{equation}\label{fg11}
\begin{split}
E_-(w_+,&w_-,w_0)=-e^{-iU_0}P_{-\mathrm{high}}[\partial_x(e^{-iU_0}w_++e^{iU_0}w_-+w_0)^2/2]\\
&-e^{-iU_0}P_{-\mathrm{high}}[\partial_x[u_0\cdot P_{+\mathrm{high}}(e^{-iU_0}w_+)+u_0\cdot P_{\mathrm{low}}(w_0)]]\\
&+e^{-iU_0}(P_{+\mathrm{high}}+P_{\mathrm{low}})[\partial_x(u_0\cdot P_{-\mathrm{high}}(e^{iU_0}w_-))]\\
&-2iP_+[\partial_x^2(e^{-iU_0}P_{-\mathrm{high}}(e^{iU_0}w_-))]\\
&-P_-\partial_xu_0\cdot w_-.
\end{split}
\end{equation}
Finally, using \eqref{fg3},
\begin{equation}\label{fg22}
\begin{cases}
&(\partial_t+\mathcal{H}\partial_x^2)w_0=E_0(w_+,w_-,w_0);\\
&w_0(0)=0,
\end{cases}
\end{equation}
where
\begin{equation}\label{fg12}
E_0(w_+,w_-,w_0)=-\frac{1}{2}P_{\mathrm{low}}[\partial_x[(e^{-iU_0}w_++e^{iU_0}w_-+w_0+u_0)^2-u_0^2]].
\end{equation}
We summarize our construction in the following lemma:

\newtheorem{Lemmaw1}{Lemma}[section]
\begin{Lemmaw1}\label{Lemmaw1}
Assume $\phi\in H^\infty_r$ and $u=S^\infty(\phi)\in C(\mathbb{R}:H^\infty_r)$. Then
\begin{equation*}
u=e^{-iU_0}w_++e^{iU_0}w_-+w_0+u_0,
\end{equation*}
where $u_0=S^\infty(P_{\mathrm{low}}(\phi))$ satisfies
\eqref{fg0}, $U_0$ satisfies \eqref{fg9}, and $w_+$, $w_-$, and
$w_0$ satisfy the equations \eqref{fg20}, \eqref{fg21}, and
\eqref{fg22}, where $E_+$, $E_-$, and $E_0$ are as in
\eqref{fg10}, \eqref{fg11}, and \eqref{fg12}.
\end{Lemmaw1}

{\bf{Remark:}} The expressions $E_+$ and $E_+$ in \eqref{fg10} and \eqref{fg11} appear complicated due to the various terms. We observe however that only the nonlinear terms in the first lines are difficult to handle: the terms in the second, third, and fourth lines are essentially of the form 
\begin{equation*}
P_{\pm}[\mathrm{smooth}\cdot P_{\mp}(\mathrm{rough})].  
\end{equation*}
Such expressions have a strong smoothing effect on the rough function. Also, in the term in the fifth line, the derivative acts on the smooth function.

\section{The Banach spaces}\label{section2}

Let $\eta_0:\mathbb{R}\to[0,1]$ denote an even smooth function
supported in $[-8/5,8/5]$ and equal to $1$ in $[-5/4,5/4]$. For
$l\in\mathbb{Z}$ let
$\chi_l(\xi)=\eta_0(\xi/2^l)-\eta_0(\xi/2^{l-1})$, $\chi_l$
supported in $\{\xi:|\xi|\in[(5/8)\cdot2^{l},(8/5)\cdot2^{l}]\}$, and
\begin{equation*}
\chi_{[l_1,l_2]}=\sum_{l=l_1}^{l_2}\chi_l\text{ for any }l_1\leq l_2\in\mathbb{Z}.
\end{equation*}
For simplicity of notation, let $\eta_l=\chi_l$ if $l\geq 1$ and
$\eta_l\equiv 0$ if $l\leq -1$. Also, for $l_1\leq
l_2\in\mathbb{Z}$ let
\begin{equation*}
\eta_{[l_1,l_2]}=\sum_{l=l_1}^{l_2}\eta_l\text{ and }\eta_{\leq l_2}=\sum_{l=-\infty}^{l_2}\eta_l.
\end{equation*}
For any integer $k\geq 0$ and $\phi\in L^2(\mathbb{R})$ we define
the operator $P_k$ by the formula
\begin{equation*}
\widehat{P_k\phi}(\xi)=\eta_k(\xi)\widehat{\phi}(\xi).
\end{equation*}
By a slight abuse of notation we also define the operators $P_k$ on $L^2(\mathbb{R}\times\mathbb{R})$ by the formula $\mathcal{F}(P_ku)(\xi,\tau)=\eta_k(\xi)\mathcal{F}(u)(\xi,\tau)$.

For $l\in\mathbb{Z}$ let
$I_l=\{\xi\in\mathbb{R}:|\xi|\in[2^{l-1},2^{l+1}]\}$. For
$l\in[0,\infty)\cap\mathbb{Z}$ let $\widetilde{I}_l=[-2,2]$ if
$l=0$ and $\widetilde{I}_l=I_l$ if $l\geq 1$. For $k\in\mathbb{Z}$
and $j\geq 0$ let
\begin{equation*}
\begin{cases}
&D_{k,j}=\{(\xi,\tau)\in\mathbb{R}\times\mathbb{R}:\xi\in
I_k,\,\tau-\omega(\xi)\in\widetilde{I}_j\}\text{ if }k\geq 1;\\
&D_{k,j}=\{(\xi,\tau)\in\mathbb{R}\times\mathbb{R}:\xi\in
I_k,\,\tau\in\widetilde{I}_j\}\text{ if }k\leq 0.
\end{cases}
\end{equation*}

Let $\mathbb{Z}_+=\mathbb{Z}\cap[0,\infty)$. For $\xi\in\mathbb{R}$ let
\begin{equation}\label{omega}
\omega(\xi)=-\xi|\xi|.
\end{equation}
We define
first the Banach spaces $X_k=X_k(\mathbb{R}\times\mathbb{R})$,
$k\in\mathbb{Z}_+$: for $k\geq 1$ we define
\begin{equation}\label{def1}
\begin{split}
X_k=&\{f\in L^2:\, f \text{ supported in }I_k\times\mathbb{R}\text{ and }\\
&||f||_{X_k}:=\sum_{j=0}^\infty 2^{j/2}\beta_{k,j}||\eta_j(\tau-\omega(\xi))f(\xi,\tau)\,||_{L^2_{\xi,\tau}}<\infty\},
\end{split}
\end{equation}
where
\begin{equation}\label{def1'}
\beta_{k,j}=1+2^{(j-2k)/2}.
\end{equation}
The precise choice of the coefficients $\beta_{k,j}$ is important in order for all the bilinear estimates \eqref{bj1}, \eqref{bj2}, \eqref{bk1}, and \eqref{bk2} to hold. Notice that $2^{j/2}\beta_{k,j}\approx 2^j$ when $k$ is small. For $k=0$ we define
\begin{equation}\label{def1''}
\begin{split}
X_0=&\{f\in L^2:\, f \text{ supported in }\widetilde{I}_0\times\mathbb{R}\text{ and }\\
&||f||_{X_0}:=\sum_{j=0}^\infty\sum_{k'=-\infty}^12^{j-k'}||\eta_j(\tau)\chi_{k'}(\xi)f(\xi,\tau)\,||_{L^2_{\xi,\tau}}<\infty\}.
\end{split}
\end{equation}
The spaces $X_k$ are not sufficient for our purpose, due to
various logarithmic divergences involving the modulation variable. For $k\geq 100$ and $k=0$ we also define the
Banach spaces $Y_k=Y_k(\mathbb{R}\times\mathbb{R})$. Let $\mathcal{F}$ and $\mathcal{F}_1$ denote the Fourier transform operators on $\mathcal{S}'(\mathbb{R}\times\mathbb{R})$ and $\mathcal{S}'(\mathbb{R})$ respectively. For $k\geq 100$ we define
\begin{equation}\label{def2}
\begin{split}
Y_k=\{f\in L^2:\,&f\text{ supported in } \bigcup_{j=0}^{k-1}D_{k,j}\text{ and }\\
&||f||_{Y_k}:=2^{-k/2}||\mathcal{F}^{-1}[(\tau-\omega(\xi)+i)
f(\xi,\tau)]||_{L^1_xL^2_t}<\infty\}.
\end{split}
\end{equation}
For $k=0$ we define
\begin{equation}\label{def2''}
\begin{split}
Y_0=\{f\in L^2:\,&f\text{ supported in } \widetilde{I}_0\times\mathbb{R}\text{ and }\\
&||f||_{Y_0}:=\sum_{j=0}^\infty2^j||\mathcal{F}^{-1}[\eta_j(\tau)
f(\xi,\tau)]||_{L^1_xL^2_t}<\infty\}.
\end{split}
\end{equation}
Then we define
\begin{equation}\label{def3}
Z_k:=X_k\text{ if }1\leq k\leq 99\text{ and }Z_k:=X_k+Y_k\text{ if }k\geq 100\text{ or }k=0.
\end{equation}
The spaces $Z_k$ are our basic Banach spaces. The spaces $X_k$ are
$X^{s,b}$-type spaces; the spaces $Y_k$ are relevant due to the
local smoothing inequality
\begin{equation*}
||\partial_xu||_{L^\infty_xL^2_t}\leq
C||(\partial_t+\mathcal{H}\partial_x^2)u||_{L^1_xL^2_t}\text{ for
any }u\in\mathcal{S}(\mathbb{R}\times\mathbb{R}).
\end{equation*}

{\bf{Remark:}} For $k\in[1,99]\cap\mathbb{Z}$ we could define the
spaces $Y_k$ as in \eqref{def2} and let $Z_k:=X_k+Y_k$. This is
not necessary, however, in view of Lemma \ref{Lemmaa1} (b) below.

In some estimates we will also need the space $\overline{Z}_0$, $Z_0\subseteq\overline{Z}_0$,
\begin{equation}\label{def1'''}
\begin{split}
\overline{Z}_0=&\{f\in L^2(\mathbb{R}\times\mathbb{R}):\, f \text{ supported in }\widetilde{I}_0\times\mathbb{R}\text{ and }\\
&||f||_{\overline{Z}_0}:=\sum_{j=0}^\infty2^{j}||\eta_j(\tau)f(\xi,\tau)\,||_{L^2_{\xi,\tau}}<\infty\}.
\end{split}
\end{equation}
We also define the space $B_0(\mathbb{R})$ by
\begin{equation}\label{def4'}
\begin{split}
B_0=&\{f\in L^2(\mathbb{R}):\, f \text{ supported in }\widetilde{I}_0\text{ and }\\
&||f||_{B_0}:=\inf_{f=g+h}||\mathcal{F}_1^{-1}(g)||_{L^1_x}+\sum_{k'=-\infty}^{1}2^{-k'}||\chi_{k'}\cdot h||_{L^2_{\xi}}<\infty\}.
\end{split}
\end{equation}

For $k\in\mathbb{Z}_+$ let
\begin{equation*}
\begin{cases}
&A_k(\xi,\tau)=\tau-\omega(\xi)+i\text{ if }k\geq 1;\\
&A_k(\xi,\tau)=\tau+i\text{ if }k=0.\\
\end{cases}
\end{equation*}
For $\sigma\geq 0$ we define the Banach spaces
$\widetilde{H}^{\sigma}=\widetilde{H}^{\sigma}(\mathbb{R})$,
$F^{\sigma}=F^{\sigma}(\mathbb{R}\times\mathbb{R})$, and
$N^{\sigma}=N^{\sigma}(\mathbb{R}\times\mathbb{R})$:
\begin{equation}\label{def4}
\begin{split}
&\widetilde{H}^{\sigma}=\Big\{\phi\in L^2:||\phi||_{\widetilde{H}^{\sigma}}^2:=||\eta_0\cdot\mathcal{F}_1(\phi)||_{B_0}^2+\sum_{k=1}^\infty 2^{2\sigma k}||\eta_k\cdot\mathcal{F}_1(\phi)||_{L^2}^2<\infty\Big\},
\end{split}
\end{equation}
\begin{equation}\label{def5}
\begin{split}
F^{\sigma}&=\Big\{u\in\mathcal{S}'(\mathbb{R}\times\mathbb{R}):||u||_{F^{\sigma}}^2:=\sum_{k=0}^\infty 2^{2\sigma k}
||\eta_k(\xi)(I-\partial_\tau^2)\mathcal{F}(u)||_{Z_k}^2
<\infty\Big\},
\end{split}
\end{equation}
and
\begin{equation}\label{def6}
\begin{split}
N^{\sigma}&=\Big\{u\in\mathcal{S}'(\mathbb{R}\times\mathbb{R}):||u||_{N^{\sigma}}^2:=\sum_{k=0}^\infty 2^{2\sigma k}||\eta_k(\xi)
A_k(\xi,\tau)^{-1}\mathcal{F}(u)||_{Z_k}^2<\infty\Big\}.
\end{split}
\end{equation}

\section{Properties of the spaces $Z_k$}\label{prop}

We start with some basic properties of the spaces
$Z_k$. Using the definitions, if $k\geq 1$ and $f_k\in Z_k$ then $f_k$ can be written in the form
\begin{equation}\label{repr1}
\begin{cases}
&f_k=\sum\limits_{j=0}^\infty f_{k,j}+g_k;\\
&\sum\limits_{j=0}^\infty 2^{j/2}\beta_{k,j}||f_{k,j}||_{L^2}+||g_k||_{Y_k}\leq 2||f_k||_{Z_k},
\end{cases}
\end{equation}
such that $f_{k,j}$ is supported in $D_{k,j}$ and $g_k$ is supported in $\bigcup_{j=0}^{k-1} D_{k,j}$ (if $k\leq 99$ then $g_k\equiv 0$).
If $f_0\in Z_0$ then $f_0$ can be written in the form
\begin{equation}\label{repr2}
\begin{cases}
&f_0=\sum\limits_{j=0}^\infty \sum\limits_{k'=-\infty}^1f^{k'}_{0,j}+\sum\limits_{j=0}^\infty g_{0,j};\\
&\sum\limits_{j=0}^\infty \sum\limits_{k'=-\infty}^12^{j-k'}||f^{k'}_{0,j}||_{L^2}+\sum\limits_{j=0}^\infty 2^j||\mathcal{F}^{-1}(g_{0,j})||_{L^1_xL^2_t}\leq 2||f_0||_{Z_0},
\end{cases}
\end{equation}
such that $f^{k'}_{0,j}$ is supported in $D_{k',j}$ and $g_{0,j}$ is supported in $\widetilde{I}_0\times\widetilde{I}_j$.

\newtheorem{Lemmaa1}{Lemma}[section]
\begin{Lemmaa1}\label{Lemmaa1}
(a) If $m,m':\mathbb{R}\to\mathbb{C}$, $k\geq 0$, and $f_k\in Z_k$ then
\begin{equation}\label{lb1}
\begin{cases}
&||m(\xi)f_k(\xi,\tau)||_{Z_k}\leq C||\mathcal{F}_1^{-1}(m)||_{L^1(\mathbb{R})}||f_k||_{Z_k};\\
&||m'(\tau)f_k(\xi,\tau)||_{Z_k}\leq C||m'||_{L^\infty(\mathbb{R})}||f_k||_{Z_k}.
\end{cases}
\end{equation}

(b) If $k\geq 1$, $j\geq 0$, and $f_k\in Z_k$ then
\begin{equation}\label{lb2}
||\eta_j(\tau-\omega(\xi))f_k(\xi,\tau)||_{X_k}\leq C||f_k||_{Z_k}.
\end{equation}

(c) If $k\geq 1$, $j\in[0,k]$, and $f_k$ is supported in $I_k\times\mathbb{R}$ then
\begin{equation}\label{lb7}
||\mathcal{F}^{-1}[\eta_{\leq j}(\tau-\omega(\xi))f_k(\xi,\tau)]||_{L^1_xL^2_t}\leq C||\mathcal{F}^{-1}(f_k)||_{L^1_xL^2_t}.
\end{equation}
\end{Lemmaa1}

\begin{proof}[Proof of Lemma \ref{Lemmaa1}]
Part (a) follows directly from Plancherel theorem and the
definitions.

For part (b), we may assume $k\geq 100$, $f_k=g_k\in Y_k$, and
$j\leq k$. We notice that if $g_k\in Y_k$ then $g_k$ can be
written in the form
\begin{equation}\label{pr10}
\begin{cases}
&g_k(\xi,\tau)=2^{k/2}\chi_{[k-1,k+1]}(\xi)(\tau-\omega(\xi)+i)^{-1}
\eta_{\leq k}(\tau-\omega(\xi))\int_\mathbb{R}e^{-ix\xi}h(x,\tau)\,dx;\\
&||g_k||_{Y_k}= C||h||_{L^1_xL^2_\tau}.
\end{cases}
\end{equation}
The inequality in part (b) follows easily since $|\{\xi\in
I_k:|\tau-\omega(\xi)|\leq 2^{j+1}\}|\leq C2^{j-k}$, see
\eqref{omega}.

For part (c), using Plancherel theorem, it suffices to prove that
\begin{equation}\label{pr42}
\Big|\Big|\int_{\mathbb{R}}e^{ix\xi}\chi_{[k-1,k+1]}(\xi)\eta_{\leq j}(\tau-\omega(\xi))\,d\xi\Big|\Big|_{L^1_xL^\infty_\tau}\leq C.
\end{equation}
In proving \eqref{pr42} we may assume $k\geq 100$. Then the function in the left-hand side of \eqref{pr42} is not zero only if $\tau\approx 2^{2k}$. Simple estimates using the change of variable $\tau-\omega(\xi)=\alpha$ and integration by parts show that
\begin{equation*}
\Big|\int_{\mathbb{R}}e^{ix\xi}\chi_{[k-1,k+1]}(\xi)\eta_{\leq j}(\tau-\omega(\xi))\,d\xi\Big|\leq C\frac{2^{j-k}}{1+(2^{j-k}x)^2}
\end{equation*}
if $\tau\approx 2^{2k}$, which suffices to prove \eqref{pr42}.
\end{proof}

Using \eqref{repr1} and Lemma \ref{Lemmaa1} (b), (c), it follows easily (see the proof of Lemma \ref{Lemmab3} for a similar argument), that if $k\geq 1$ and $(I-\partial_\tau^2)f_k\in Z_k$ then $f_k$ can be written in the form
\begin{equation}\label{repr3}
\begin{cases}
&f_k=\sum\limits_{j=0}^\infty f_{k,j}+g_k;\\
&\sum\limits_{j=0}^\infty 2^{j/2}\beta_{k,j}||(I-\partial_\tau^2)f_{k,j}||_{L^2}+||(I-\partial_\tau^2)g_k||_{Y_k}\leq C||(I-\partial_\tau^2)f_k||_{Z_k},
\end{cases}
\end{equation}
such that $f_{k,j}$ is supported in $D_{k,j}$ and $g_k$ is supported in $\bigcup_{j=0}^{k-20} D_{k,j}$ (if $k\leq 99$ then $g_k\equiv 0$). We prove now several estimates using the spaces $Z_k$.

\newtheorem{Lemmaa2}[Lemmaa1]{Lemma}
\begin{Lemmaa2}\label{Lemmaa2}

(a) If $k\geq 0$, $t\in\mathbb{R}$, and $f_k\in Z_k$ then
\begin{equation}\label{lb4}
\begin{cases}
&\big|\big|\int_\mathbb{R}f_k(\xi,\tau)e^{it\tau}\,d\tau\big|\big|_{L^2_\xi}
\leq C||f_k||_{Z_k}\text{ if }k\geq 1;\\
&\big|\big|\int_\mathbb{R}f_0(\xi,\tau)e^{it\tau}\,d\tau\big|\big|_{B_0}
\leq C||f_0||_{Z_0}\text{ if }k=0.
\end{cases}
\end{equation}
As a consequence,
\begin{equation}\label{hh80}
F^\sigma\subseteq C(\mathbb{R}:\widetilde{H}^\sigma)\text{ for any }\sigma\geq 0.
\end{equation}

(b) If $k\geq 1$ and $(I-\partial_\tau^2)f_k\in Z_k$ then
\begin{equation}\label{pr40}
||\mathcal{F}^{-1}(f_k)||_{L^2_xL^\infty_t}\leq
C2^{k/2}||(I-\partial_\tau^2)f_k||_{Z_k}.
\end{equation}

(c) If $k\geq 1$ and $f_k\in Z_k$ then
\begin{equation}\label{lb6}
||\mathcal{F}^{-1}(f_k)||_{L^\infty_xL^2_t}\leq C2^{-k/2}||f_k||_{Z_k}.
\end{equation}
\end{Lemmaa2}

\begin{proof}[Proof of Lemma \ref{Lemmaa2}] For part (a), $k\geq 1$, we use the representation \eqref{repr1}. Assume first that $f_k=f_{k,j}$. Then
\begin{equation*}
\begin{split}
\Big|\Big|\int_\mathbb{R}f_{k,j}(\xi,\tau)e^{it\tau}\,d\tau\Big|\Big|_{L^2_\xi}\leq C||f_{k,j}(\xi,\tau)||_{L^2_\xi L^1_\tau}\leq C2^{j/2}||f_{k,j}||_{L^2_{\xi,\tau}},
\end{split}
\end{equation*}
which proves \eqref{lb4} in this case.

Assume now that $k\geq 100$, $f_k=g_k\in Y_k$, and
write $g_k$ as in \eqref{pr10}. We define the modified Hilbert
transform operator
\begin{equation}\label{pr15}
\mathcal{L}_k(g)(\mu)=\int_\mathbb{R}g(\tau)(\tau-\mu+i)^{-1}
\eta_{[0,k]}(\tau-\mu)\,d\tau,\,\,g\in L^2(\mathbb{R}).
\end{equation}
Clearly, $||\mathcal{L}_k||_{L^2\to L^2}\leq C$, uniformly in $k$.
We examine the formula \eqref{pr10} and let $h^\ast(x,\mu)=\mathcal{L}_k[e^{it\tau}h(x,\tau)](\mu)$,
$||h^\ast||_{L^1_xL^2_\mu}\leq C||h||_{L^1_xL^2_\tau}$.
Then, using \eqref{pr10}
\begin{equation}\label{pr60}
\begin{split}
\Big|\Big|\int_\mathbb{R}g_k(\xi,\tau)e^{it\tau}\,d\tau\Big|\Big|_{L^2_\xi}
&\leq C2^{k/2}\Big|\Big|\chi_{[k-1,k+1]}(\xi)\int_\mathbb{R}e^{-ix\xi} h^\ast(x,\omega(\xi))\,dx\Big|\Big|_{L^2_\xi}\\
&\leq C2^{k/2}\cdot 2^{-k/2}||h^\ast||_{L^1_xL^2_\mu}\\
&\leq C||g_k||_{Y_k},
\end{split}
\end{equation}
which completes the proof of \eqref{lb4} in the case $k\geq 1$.

Assume now $k=0$. We use the representation \eqref{repr2}. Assume first that $f_0=f_{0,j}^{k'}$ is supported in $D_{k',j}$, $||f_0||_{Z_0}\approx 2^{j-k'}||f_{0,j}^{k'}||_{L^2}$. Then
\begin{equation*}
\big|\big|\int_\mathbb{R}f_{0,j}^{k'}(\xi,\tau)e^{it\tau}\,d\tau\big|\big|_{B_0}\leq C2^{-k'}\big|\big|\int_\mathbb{R}|f_{0,j}^{k'}(\xi,\tau)|\,d\tau\big|\big|_{L^2_\xi}\leq C2^{-k'}2^{j/2}||f_{0,j}^{k'}||_{L^2},
\end{equation*}
which suffices. 

Assume now that $f_0=g_{0,j}$ is supported in $\widetilde{I}_0\times\widetilde{I}_j$, $||f_0||_{Z_0}\approx 2^j||\mathcal{F}^{-1}(g_{0,j})||_{L^1_xL^2_t}$. Then
\begin{equation*}
\big|\big|\int_\mathbb{R}g_{0,j}(\xi,\tau)e^{it\tau}\,d\tau\big|\big|_{B_0}\leq C||\mathcal{F}^{-1}(g_{0,j})||_{L^1_xL^\infty_t}\leq C2^{j/2}||\mathcal{F}^{-1}(g_{0,j})||_{L^1_xL^2_t},
\end{equation*}
which completes the proof of part (a). 

For part (b) we use the representation \eqref{repr3}. Assume first that $f_k=f_{k,j}$ and let $f^\#_{k,j}(\xi,\mu)=f_{k,j}(\xi,\mu+\omega(\xi))$. By integration by parts, the left-hand side of \eqref{pr40} is dominated by
\begin{equation*}
\begin{split}
\sum_{n\in\mathbb{Z}}\frac{C}{n^2+1}\int_{\widetilde{I}_j}\Big|\Big|\int_{\mathbb{R}}(I-\partial_\mu^2)f^\#_{k,j}(\xi,\mu)e^{ix\xi}e^{it\omega(\xi)}\,d\xi\Big|\Big|_{L^2_xL^\infty_{t\in[n-1/2,n+1/2]}}\,d\mu.
\end{split}
\end{equation*}
The bound \eqref{pr40} now follows from the standard maximal function estimate
\begin{equation}\label{ar400}
\Big|\Big|\int_{\mathbb{R}}g(\xi)e^{ix\xi}e^{it\omega(\xi)}\,d\xi\Big|\Big|_{L^2_xL^\infty_{t\in[-1/2,1/2]}}\leq C2^{k/2}||g||_{L^2_\xi},
\end{equation}
for any function $g$ supported in $I_k$, see \cite[Theorem 2.7]{KePoVe1}. In fact, the  argument above and \eqref{ar400} show that if $f_k\in X_k$ then
\begin{equation}\label{ar399}
\Big|\Big|\int_{\mathbb{R}^2}f_k(\xi,\tau)e^{ix\xi}e^{it\tau}\,d\xi d\tau\Big|\Big|_{L^2_xL^\infty_{t\in[-1/2,1/2]}}\leq  C2^{k/2}||f_k||_{X_k}.
\end{equation}

{\bf{Remark:}} The inequality \eqref{pr40} is relevant only when $j\leq k$. For $j\geq k$ the Sobolev imbedding theorem easily gives a stronger estimate.

Assume now that $k\geq 100$, $f_k=g_k$, $(I-\partial_\tau^2)g_k\in Y_k$. By integration by parts, the left-hand side of \eqref{pr40} is dominated by
\begin{equation*}
\begin{split}
\sum_{n\in\mathbb{Z}}\frac{C}{n^2+1}\Big|\Big|\int_{\mathbb{R}^2}(I-\partial_\tau^2)g_k(\xi,\tau)e^{ix\xi}e^{it\tau}\,d\xi d\tau\Big|\Big|_{L^2_xL^\infty_{t\in[n-1/2,n+1/2]}}.
\end{split}
\end{equation*}
We write now $(I-\partial_\tau^2)g_k$ as in \eqref{pr10}. In view of \eqref{pr10}, it suffices to prove that if
\begin{equation}\label{ar402}
f(\xi,\tau)=2^{k/2}\chi_{[k-1,k+1]}(\xi)(\tau-\omega(\xi)+i)^{-1}
\eta_{\leq k}(\tau-\omega(\xi))\cdot h(\tau)
\end{equation}
then
\begin{equation}\label{ar401}
\Big|\Big|\int_{\mathbb{R}^2}f(\xi,\tau)e^{ix\xi}e^{it\tau}\,d\xi d\tau\Big|\Big|_{L^2_xL^\infty_{t\in[-1/2,1/2]}}\leq  C2^{k/2}||h||_{L^2}.
\end{equation}
Since $k\geq 100$ and $|\xi|\in[2^{k-2},2^{k+2}]$, we may assume that the function $h$ in \eqref{ar402} is supported in the set $\{t:|\tau|\in[2^{2k-10},2^{2k+10}]\}$. Let $h_+=h\cdot\mathbf{1}_{[0,\infty)}$, $h_-=h\cdot\mathbf{1}_{(-\infty,0]}$, and define the corresponding functions $f_+$ and $f_-$ as in \eqref{ar402}. By symmetry, it suffices to prove the bound \eqref{ar401} for the function $f_+$, which is supported in the set $\{(\xi,\tau):\xi\in[-2^{k+2},-2^{k-2}],\,\tau\in[2^{2k-10},2^{2k+10}]\}$. In view of \eqref{omega}, $\tau-\omega(\xi)=\tau-\xi^2$ on the support of $f_+$, and $f_+(\xi,\tau)=0$ unless $|\sqrt\tau+\xi|\leq C$. Let
\begin{equation}\label{ar403}
\begin{split}
f'_+(\xi,\tau)=2^{k/2}\chi_{[k-1,k+1]}(-\sqrt\tau)(\tau-&\xi^2+(\sqrt\tau+\xi)^2+i\sqrt\tau2^{-k})^{-1}\\
&\eta_{0}(\sqrt\tau+\xi)\cdot h_+(\tau).
\end{split}
\end{equation}
Using Lemma \ref{Lemmaa1} (b), it is easy to see that
\begin{equation*}
||f_+-f'_+||_{X_k}\leq C||h_+||_{L^2}.
\end{equation*}
Thus, using \eqref{ar399}, $||\mathcal{F}^{-1}(f_+-f'_+)||_{L^2_xL^\infty_{t\in[-1/2,1/2]}}\leq  C2^{k/2}||h_+||_{L^2}$. To  estimate $||\mathcal{F}^{-1}(f'_+)||_{L^2_xL^\infty_{t\in[-1/2,1/2]}}$ we make the change of variables $\xi=-\sqrt\tau+\mu$. Then
\begin{equation}\label{ar420}
\begin{split}
\mathcal{F}^{-1}(f'_+)(x,t)=2^{k/2}\int_\mathbb{R}h_+(\tau)&(2\sqrt\tau)^{-1}\chi_{[k-1,k+1]}(-\sqrt\tau)e^{it\tau}e^{-ix\sqrt\tau}\,d\tau\\
&\times\int_\mathbb{R}\eta_0(\mu)(\mu+i/2^{k+1})^{-1}e^{ix\mu}\,d\mu.
\end{split}
\end{equation}
The absolute value of the integral in $\mu$ in \eqref{ar420} is bounded by $C$. We make the change of variables $\tau=\theta^2$ in the first integral, and  use the bound \eqref{ar400}. It follows that $||\mathcal{F}^{-1}(f'_+)||_{L^2_xL^\infty_{t\in[-1/2,1/2]}}\leq  C2^{k/2}||h_+||_{L^2}$, which completes the proof of \eqref{ar401}.

For part (c) we use the representation \eqref{repr1}. Assume first that $f_k=f_{k,j}$ and let $f^\#_{k,j}(\xi,\mu)=f_{k,j}(\xi,\mu+\omega(\xi))$. It suffices to prove the stronger bound
\begin{equation*}
\Big|\Big|\int_{D_{k,j}}f_{k,j}(\xi,\tau)e^{ix_0\xi}e^{it\tau}\,d\xi d\tau\Big|\Big|_{L^2_t}\leq C2^{-k/2}2^{j/2}||f_{k,j}||_{L^2},
\end{equation*}
for any $x_0\in\mathbb{R}$. Using Plancherel theorem, duality, and the Cauchy--Schwartz inequality, the left-hand side of the inequality above is dominated by
\begin{equation*}
\begin{split}
&C\sup_{||h||_{L^2(\mathbb{R})=1}}\int_{I_k\times\widetilde{I}_j}|f^\#_{k,j}(\xi,\mu)|\cdot |h(\mu+\omega(\xi))|\,d\xi d\mu\\
&\leq C\sup_{||h||_{L^2(\mathbb{R})=1}}\int_{\widetilde{I}_j}\Big(\int_{I_k}|f^\#_{k,j}(\xi,\mu)|^2\,d\xi\Big)^{1/2}\Big(\int_{I_k}|h(\mu+\omega(\xi))|^2\,d\xi\Big)^{1/2}\,d\mu\\
&\leq C2^{-k/2}2^{j/2}\Big(\int_{I_k\times\widetilde{I}_j}|f^\#_{k,j}(\xi,\mu)|^2\,d\xi d\mu\Big)^{1/2},
\end{split}
\end{equation*}
as desired.

Assume now that $k\geq 100$, $f_k=g_k\in Y_k$ and write $g_k$ as in \eqref{pr10}. Using Plancherel theorem, it suffices to prove that
\begin{equation}\label{pr41}
\int_{\mathbb{R}}e^{ix_0\xi}\chi_{[k-1,k+1]}(\xi)(\tau-\omega(\xi)+i)^{-1}\eta_{\leq k}(\tau-\omega(\xi))\,d\xi\leq C2^{-k},
\end{equation}
uniformly in $x_0$ and $\tau$ (assuming $k\geq 100$). We  may assume $|\tau|\in[2^{2k-10},2^{2k+10}]$ and, by symmetry, $\tau \geq 0$. Then the variable $\xi$ in the integral in \eqref{pr41} is in the interval $[-\sqrt\tau-C,-\sqrt\tau+C]$ and $\tau-\omega(\xi)=\tau-\xi^2$. As in part (b), see \eqref{ar403}, we replace the integrand  $\mathbf{1}_{(-\infty,,0]}(\xi)\chi_{[k-1,k+1]}(\xi)(\tau-\xi^2+i)^{-1}\eta_{\leq k}(\tau-\xi^2)$ with $\chi_{[k-1,k+1]}(-\sqrt\tau)(\tau-\xi^2+(\sqrt\tau+\xi)^2+i\sqrt\tau2^{-k})^{-1}\eta_{0}(\sqrt\tau+\xi)$ at the expense  of an error dominated by
\begin{equation*}
C[2^{-k}+(2^{2k}|\sqrt\tau+\xi|^2+1)^{-1}]\mathbf{1}_{[0,C]}(|\sqrt\tau+\xi|).
\end{equation*}
The $L^1_\xi$ norm of this error is $\leq C2^{-k}$. Then we make the change of  variables $\xi=-\sqrt\tau+\mu$ and  use the uniform boundedness of the integral in $\mu$ in \eqref{ar420}. The bound \eqref{pr41} follows.
\end{proof}

\section{Linear estimates}\label{linear}

For any $u\in C(\mathbb{R}:L^2)$ let $\widetilde{u}(.,t)\in C(\mathbb{R}:L^2)$
denote its partial Fourier transform with respect to the variable
$x$. For $\phi\in L^2(\mathbb{R})$ let $W(t)\phi\in
C(\mathbb{R}:L^2)$ denote the solution of the free Benjamin--Ono evolution
given by
\begin{equation}\label{ni1}
[W(t)\phi]\,\,\widetilde{}\,\,(\xi,t)=e^{it\omega(\xi)}\widehat{\phi}(\xi),
\end{equation}
where $\omega(\xi)$ is defined in \eqref{omega}. Assume
$\psi:\mathbb{R}\to[0,1]$ is an even smooth function supported in
the interval $[-8/5,8/5]$ and equal to $1$ in the interval
$[-5/4,5/4]$ and let
$\varphi=\widehat{\psi}-\widehat{\psi}''\in\mathcal{S}(\mathbb{R})$.

\newtheorem{Lemmab1}{Lemma}[section]
\begin{Lemmab1}\label{Lemmab1}
If $\sigma\geq 0$ and $\phi\in \widetilde{H}^{\sigma}$ then
\begin{equation*}
||\psi(t)\cdot (W(t)\phi)||_{F^{\sigma}}\leq C||\phi||_{\widetilde{H}^{\sigma}}.
\end{equation*}
\end{Lemmab1}

\begin{proof}[Proof of Lemma \ref{Lemmab1}] A straightforward computation shows that
\begin{equation}\label{ni3}
\mathcal{F}[\psi(t)\cdot (W(t)\phi)](\xi,\tau)=
\widehat{\phi}(\xi)\widehat{\psi}(\tau-\omega(\xi)).
\end{equation}
Then, directly from the definitions,
\begin{equation}\label{ar200}
\begin{split}
&||\psi(t)\cdot (W(t)\phi)||_{F^{\sigma}}^2=\sum_{k\in\mathbb{Z}_+}2^{2\sigma k}||\eta_k(\xi)\widehat{\phi}(\xi)\varphi(\tau-\omega(\xi))||_{Z_k}^2\\
&\leq \sum_{k=1}^\infty2^{2\sigma k}||\eta_k(\xi)\widehat{\phi}(\xi)\varphi(\tau-\omega(\xi))||_{X_k}^2+||\eta_0(\xi)\widehat{\phi}(\xi)\varphi(\tau-\omega(\xi))||_{Z_0}^2.
\end{split}
\end{equation}
Since $\varphi\in\mathcal{S}(\mathbb{R})$, for any $k\geq 1$
\begin{equation*}
||\eta_k(\xi)\widehat{\phi}(\xi)\varphi(\tau-\omega(\xi))||_{X_k}\leq C||\eta_k\cdot\widehat{\phi}||_{L^2}.
\end{equation*}
For $k=0$, write $\eta_0\cdot\widehat{\phi}=g+\sum_{k'\leq 1}h_{k'}$, $h_{k'}$ supported in $I_{k'}$ and
\begin{equation}\label{hh20}
||\mathcal{F}_1^{-1}(g)||_{L^1_x}+\sum_{k'\leq 1}2^{-k'}||h_{k'}||_{L^2}\leq 2||\eta_0\cdot\widehat{\phi}||_{B_0}.
\end{equation}
Then
\begin{equation*}
\begin{split}
&||g(\xi)\varphi(\tau-\omega(\xi))||_{Z_0}\leq ||g(\xi)\varphi(\tau)||_{Y_0}+||g(\xi)[\varphi(\tau-\omega(\xi))-\varphi(\tau)]||_{X_0}\\
&\leq C||\mathcal{F}_1^{-1}(g)||_{L^1_x}+C||g(\xi)\xi^2(1+|\tau|)^{-4}||_{X_0}\leq C||\mathcal{F}_1^{-1}(g)||_{L^1_x}.
\end{split}
\end{equation*}
Also,
\begin{equation*}
||h_{k'}(\xi)\varphi(\tau-\omega(\xi))||_{Z_0}\leq ||h_{k'}(\xi)\varphi(\tau-\omega(\xi))||_{X_0}\leq C2^{-k'}||h_{k'}||_{L^2}.
\end{equation*}
Lemma \ref{Lemmab1} follows from \eqref{hh20}.
\end{proof}

\newtheorem{Lemmab3}[Lemmab1]{Lemma}
\begin{Lemmab3}\label{Lemmab3}
If $\sigma\geq 0$ and $u\in N^{\sigma}\cap C(\mathbb{R}:H^{-2})$
then
\begin{equation*}
\Big|\Big|\psi(t)\cdot \int_0^tW(t-s)(u(s))\,ds\Big|\Big|_{F^{\sigma}}\leq C||u||_{N^{\sigma}}.
\end{equation*}
\end{Lemmab3}

\begin{proof}[Proof of Lemma \ref{Lemmab3}] A straightforward computation shows that
\begin{equation}\label{ni2}
\mathcal{F}\Big[\psi(t)\cdot \int_0^tW(t-s)(u(s))ds\Big](\xi,\tau)=
c\int_\mathbb{R}\mathcal{F}(u)(\xi,\tau')\frac{\widehat{\psi}(\tau-\tau')-\widehat{\psi}(\tau-\omega(\xi))}{\tau'-\omega(\xi)}d\tau'.
\end{equation}
For $k\in\mathbb{Z}_+$ let
$f_k(\xi,\tau')=\mathcal{F}(u)(\xi,\tau')\eta_k(\xi)A_k(\xi,\tau')^{-1}$.
For $f_k\in Z_k$ let
\begin{equation}\label{ar202}
T(f_k)(\xi,\tau)=\int_\mathbb{R}f_k(\xi,\tau')\frac{\varphi(\tau-\tau')-
\varphi(\tau-\omega(\xi))}{\tau'-\omega(\xi)}A_k(\xi,\tau')\,d\tau'.
\end{equation}
In view of the definitions, it suffices to prove that
\begin{equation}\label{ni5}
||T||_{Z_k\to Z_k}\leq C\text{ uniformly in }k\in\mathbb{Z}_+.
\end{equation}

We consider first the case $k\geq 1$. To prove \eqref{ni5} we use the representation \eqref{repr1}. Assume first that $f_k=f_{k,j}$ is a function supported in $D_{k,j}$. Let $f_{k,j}^\#(\xi,\mu')=f_{k,j}(\xi,\mu'+\omega(\xi))$ and $T(f_{k,j})^\#(\xi,\mu)=T(f_{k,j})(\xi,\mu+\omega(\xi))$. Then,
\begin{equation}\label{ni6}
T(f_{k,j})^\#(\xi,\mu)=\int_\mathbb{R}f_{k,j}^\#(\xi,\mu')\frac{\varphi(\mu-\mu')-
\varphi(\mu)}{\mu'}(\mu'+i)\,d\mu'.
\end{equation}
We use the elementary bound
\begin{equation*}
\Big|\frac{\varphi(\mu-\mu')-\varphi(\mu)}{\mu'}(\mu'+i)\Big|\leq C[(1+|\mu|)^{-4}+(1+|\mu-\mu'|)^{-4}].
\end{equation*}
Then, using \eqref{ni6},
\begin{equation*}
\begin{split}
|T(f_{k,j})^\#(\xi,\mu)|&\leq C(1+|\mu|)^{-4}2^{j/2}\Big[\int_{\widetilde{I}_j}|f_{k,j}^\#(\xi,\mu')|^2\,d\mu'\Big]^{1/2}\\
&+C\eta_{[j-2,j+2]}(\mu)\int_{\widetilde{I}_j}|f_{k,j}^\#(\xi,\mu')|(1+|\mu-\mu'|)^{-4}\,d\mu'.
\end{split}
\end{equation*}
It follows from the definition of the spaces $X_k$ that
\begin{equation}\label{ni7}
||T||_{X_k\to X_k}\leq C\text{ uniformly in }k\geq 1,
\end{equation}
as desired.

Assume now that $f_k=g_k\in Y_k$, so $k\geq 100$. In view of Lemma
\ref{Lemmaa1} (b), (c), and \eqref{ni7}, we  may assume that $g_k$
is supported in the set $\{(\xi,\tau'):|\tau'-\omega(\xi)|\leq
2^{k-20}\}$. We write
\begin{equation*}
g_k(\xi,\tau')=\frac{\tau'-\omega(\xi)}{\tau'-\omega(\xi)+i}g_k(\xi,\tau')+\frac{i}{\tau'-\omega(\xi)+i}g_k(\xi,\tau').
\end{equation*}
Using Lemma \ref{Lemmaa1} (b), $||i(\tau'-\omega(\xi)+i)^{-1}g_k(\xi,\tau')||_{X_k}\leq C||g_k||_{Y_k}$. In view of \eqref{ni7}, it suffices to prove that
\begin{equation}\label{ni8}
\Big|\Big|\int_\mathbb{R}g_k(\xi,\tau')\varphi(\tau-\tau')\,d\tau'\Big|\Big|_{Z_k}+
\Big|\Big|\varphi(\tau-\omega(\xi))\int_\mathbb{R}g_k(\xi,\tau')\,d\tau'\Big|\Big|_{X_k}\leq C||g_k||_{Y_k}.
\end{equation}
The bound for the second term in the left-hand side of \eqref{ni8} follows from \eqref{pr60} with $t=0$. To bound the first term we write
\begin{equation*}
g_k(\xi,\tau')=g_k(\xi,\tau')\Big[\frac{\tau'-\omega(\xi)+i}{\tau-\omega(\xi)+i}+\frac{\tau-\tau'}{\tau-\omega(\xi)+i}\Big].
\end{equation*}
The first term in the left-hand side of \eqref{ni8} is dominated by
\begin{equation}\label{ni9}
\begin{split}
&C\Big|\Big|\eta_{[0,k-1]}(\tau-\omega(\xi))(\tau-\omega(\xi)+i)^{-1}\int_\mathbb{R}g_k(\xi,\tau')(\tau'-\omega(\xi)+i)\varphi(\tau-\tau')\,d\tau'\Big|\Big|_{Y_k}\\
&+C\sum_{j\leq k}2^{j/2}\Big|\Big|\eta_{j}(\tau-\omega(\xi))(\tau-\omega(\xi)+i)^{-1}\int_\mathbb{R}g_k(\xi,\tau')\varphi(\tau-\tau')(\tau-\tau')\,d\tau'\Big|\Big|_{L^2}\\
&+C\sum_{j\geq k-1}2^{j/2}\beta_{k,j}\Big|\Big|\eta_{j}(\tau-\omega(\xi))\int_\mathbb{R}g_k(\xi,\tau')\varphi(\tau-\tau')\,d\tau'\Big|\Big|_{L^2}.
\end{split}
\end{equation}
For the first term in \eqref{ni9}, we use Lemma \ref{Lemmaa1} (c) to bound it by
\begin{equation*}
C2^{-k/2}||\mathcal{F}_1^{-1}(\varphi)\cdot\mathcal{F}^{-1}[(\tau'-\omega(\xi)+i)g_k(\xi,\tau')]||_{L^1_xL^2_t}\leq C||g_k||_{Y_k},
\end{equation*}
as desired. Let $g_k^\#(\xi,\mu')=g_k(\xi,\mu'+\omega(\xi))$ and for $j'\in[0,k-20]$ let $g_{k,j'}^\#(\xi,\mu')=g_k^\#(\xi,\mu')\eta_{j'}(\mu')$. In view of Lemma \ref{Lemmaa1} (b), $2^{j'/2}||g_{k,j'}||_{L^2}\leq  C||g_k||_{Y_k}$, so the second term in \eqref{ni9} is dominated by
\begin{equation*}
C\sum_{j=0}^k\sum_{j'=0}^{k-20}2^{-j/2}2^{-j'/2}[2^{j'/2}||g_{k,j'}||_{L^2}]\leq C||g_k||_{Y_k}.
\end{equation*}
The third term in \eqref{ni9} is dominated by
\begin{equation*}
C\sum_{j=k-1}^\infty\sum_{j'=0}^{k-20}2^{-3j}2^{j'/2}||g_{k,j'}||_{L^2}\leq C||g_k||_{Y_k}.
\end{equation*}
This completes the  proof of \eqref{ni8}.

We consider now the case $k=0$. To prove \eqref{ni5} we use the
representation \eqref{repr2}. Assume first that $f_0=f_{0,j}^{k'}$
is a function supported in $D_{k',j}$, $||f_0||_{Z_0}\approx
2^{j-k'}||f_{0,j}^{k'}||_{L^2}$. For $|\xi|\leq 2$ we have the
elementary bound
\begin{equation*}
\Big|\frac{\varphi(\tau-\tau')-\varphi(\tau-\omega(\xi))}{\tau'-\omega(\xi)}(\tau'+i)\Big|\leq C[(1+|\tau|)^{-4}+(1+|\tau-\tau'|)^{-4}].
\end{equation*}
Then, using the formula \eqref{ar202},
\begin{equation*}
\begin{split}
|T(f_{0,j}^{k'})(\xi,\tau)|&\leq C(1+|\tau|)^{-4}2^{j/2}\Big[\int_{\widetilde{I}_j}|f_{0,j}^{k'}(\xi,\tau')|^2\,d\tau'\Big]^{1/2}\\
&+C\eta_{[j-4,j+4]}(\tau)\int_{\widetilde{I}_j}|f_{0,j}^{k'}(\xi,\tau')|(1+|\tau-\tau'|)^{-4}\,d\tau'.
\end{split}
\end{equation*}
It follows from the definition of the spaces $X_0$ that $||T||_{X_0\to X_0}\leq C$, as desired.

Assume now that $f_0=g_{0,j}$ is supported in $\widetilde{I}_0\times\widetilde{I}_j$. We can write
\begin{equation}\label{ar209}
\begin{cases}
&g_{0,j}(\xi,\tau')=2^{-j}\eta_{[0,1]}(\xi)\eta_{[j-1,j+1]}(\tau')\int_{\mathbb{R}}e^{-ix\xi}h(x,\tau')\,dx;\\
&2^j||\mathcal{F}^{-1}(g_{0,j})||_{L^1_xL^2_t}=C||h||_{L^1_xL^2_{\tau'}}.
\end{cases}
\end{equation}

We have two cases: $j\leq 5$ and $j\geq 6$. If $j\leq 5$ we write
\begin{equation*}
\frac{\varphi(\tau-\tau')-\varphi(\tau-\omega(\xi))}{\tau'-\omega(\xi)}=c\int_0^1\varphi'(\tau-\alpha\tau'-(1-\alpha)\omega(\xi))\,d\alpha.
\end{equation*}
For \eqref{ni5}, it suffices to prove that
\begin{equation}\label{ar210}
\Big|\Big|\int_\mathbb{R}g_{0,j}(\xi,\tau')\varphi'(\tau-\alpha\tau'-(1-\alpha)\omega(\xi))(\tau'+i)\,d\tau'\Big|\Big|_{Z_0}\leq C||\mathcal{F}^{-1}(g_{0,j})||_{L^1_xL^2_t}\end{equation}
for any $\alpha\in[0,1]$. For $|\xi|\leq 2$ and $|\tau'|\leq C$ we write
\begin{equation*}
\varphi'(\tau-\alpha\tau'-(1-\alpha)\omega(\xi))(\tau'+i)=\varphi'(\tau-\alpha\tau')(\tau'+i)+R(\xi,\tau,\tau'),
\end{equation*}
where
\begin{equation*}
|R(\xi,\tau,\tau')|\leq C\xi^2(1+|\tau|)^{-4}.
\end{equation*}
The left-hand side of \eqref{ar210} is dominated by
\begin{equation*}
\Big|\Big|\int_\mathbb{R}g_{0,j}(\xi,\tau')\varphi'(\tau-\alpha\tau')(\tau'+i)\,d\tau'\Big|\Big|_{Y_0}+C\Big|\Big|\xi^2(1+|\tau|)^{-4}\int_{\widetilde{I}_j}|g_{0,j}(\xi,\tau')|\,d\tau'\Big|\Big|_{X_0},
\end{equation*}
which is easily seen to be dominated by $||\mathcal{F}^{-1}(g_{0,j})||_{L^1_xL^2_t}$ (using the representation \eqref{ar209}). This completes the proof of \eqref{ni5} in the case $j\leq 5$.

Assume now that $j\geq 6$. Since $|\tau'|\geq C$ and $|\xi|\leq 2$ we can write
\begin{equation*}
\frac{\varphi(\tau-\tau')-\varphi(\tau-\omega(\xi))}{\tau'-\omega(\xi)}(\tau'+i)=\frac{\varphi(\tau-\tau')-\varphi(\tau)}{\tau'}(\tau'+i)+R'(\xi,\tau,\tau'),
\end{equation*}
where
\begin{equation*}
|R'(\xi,\tau,\tau')|\leq C\xi^2[(1+|\tau|)^{-4}+(1+|\tau-\tau'|)^{-4}].
\end{equation*}
Using the representation \eqref{ar209} and the definitions, it follows as before that
\begin{equation*}
\Big|\Big|\int_\mathbb{R}g_{0,j}(\xi,\tau')\frac{\varphi(\tau-\tau')-\varphi(\tau)}{\tau'}(\tau'+i)\,d\tau'\Big|\Big|_{Y_0}+\Big|\Big|\int_{\mathbb{R}}|g_{0,j}(\xi,\tau')|\cdot|R'(\xi,\tau,\tau')|\,d\tau'\Big|\Big|_{X_0}
\end{equation*}
is dominated by $C2^j||\mathcal{F}^{-1}(g_{0,j})||_{L^1_xL^2_t}$, which completes the proof of \eqref{ni5}.
\end{proof}

\section{A symmetric estimate}\label{bilin}

We start with a symmetric estimate for nonnegative functions. For
$\xi_1,\xi_2\in\mathbb{R}$ and $\omega:\mathbb{R}\to\mathbb{R}$ as
in \eqref{omega} let
\begin{equation}\label{om1}
\Omega(\xi_1,\xi_2)=-\omega(\xi_1+\xi_2)+\omega(\xi_1)+\omega(\xi_2).
\end{equation}
For compactly supported functions $f,g,h\in L^2(\mathbb{R}\times\mathbb{R})$ let
\begin{equation}\label{om2}
J(f,g,h)=\int_{\mathbb{R}^4}f(\xi_1,\mu_1)g(\xi_2,\mu_2)h(\xi_1+\xi_2,\mu_1+\mu_2+\Omega(\xi_1,\xi_2))\,d\xi_1d\xi_2d\mu_1d\mu_2.
\end{equation}
Given a triplet of real numbers $(\alpha_1, \alpha_2,\alpha_3)$ let $\mathrm{min}\,(\alpha_1, \alpha_2,\alpha_3)$, $\mathrm{max}\,(\alpha_1, \alpha_2,\alpha_3)$, and $\mathrm{med}\,(\alpha_1, \alpha_2,\alpha_3)$ denote the minimum,
the maximum, and the median (i.e. $\mathrm{med}\,(\alpha_1, \alpha_2,\alpha_3)=\alpha_1+\alpha_2+\alpha_3-\mathrm{max}\,(\alpha_1, \alpha_2,\alpha_3)-\mathrm{min}\,(\alpha_1, \alpha_2,\alpha_3)$) of the numbers $\alpha_1$, $\alpha_2$,
and $\alpha_3$.

\newtheorem{Lemmad1}{Lemma}[section]
\begin{Lemmad1}\label{Lemmad1}
Assume $k_1,k_2,k_3\in\mathbb{Z}$,
$j_1,j_2,j_3\in\mathbb{Z}_+$, and $f_{k_i,j_i}\in
L^2(\mathbb{R}\times\mathbb{R})$ are functions supported in
$I_{k_i}\times\widetilde{I}_{j_i}$, $i=1,2,3$.

(a) For any $k_1,k_2,k_3\in\mathbb{Z}$ and $j_1,j_2,j_3\in\mathbb{Z}_+$,
\begin{equation}\label{om31}
|J(f_{k_1,j_1},f_{k_2,j_2},f_{k_3,j_3})|\leq C
2^{\mathrm{min}\,(k_1,k_2,k_3)/2}2^{\mathrm{min}\,(j_1,j_2,j_3)/2}
\prod_{i=1}^3||f_{k_i,j_i}||_{L^2}.
\end{equation}

(b) If $\mathrm{max}\,(k_1,k_2,k_3)\geq \mathrm{min}\,(k_1,k_2,k_3)+5$ and $i\in\{1,2,3\}$ then
\begin{equation}\label{om32}
|J(f_{k_1,j_1},f_{k_2,j_2},f_{k_3,j_3})|\leq
C2^{(j_1+j_2+j_3)/2}2^{-(j_i+k_i)/2}\prod_{i=1}^3||f_{k_i,j_i}||_{L^2}.
\end{equation}

(c) For any $k_1,k_2,k_3\in\mathbb{Z}$ and $j_1,j_2,j_3\in\mathbb{Z}_+$,
\begin{equation}\label{om3}
|J(f_{k_1,j_1},f_{k_2,j_2},f_{k_3,j_3})|\leq
C2^{\mathrm{min}\,(j_1,j_2,j_3)/2+\mathrm{med}\,(j_1,j_2,j_3)/4}\prod_{i=1}^3||f_{k_i,j_i}||_{L^2}.
\end{equation}
\end{Lemmad1}

\begin{proof}[Proof of Lemma \ref{Lemmad1}] Let
$A_{k_i}(\xi)=\Big[\int_{\mathbb{R}}|f_{k_i,j_i}(\xi,\mu)|^2\,d\mu\Big]^{1/2}$,
$i=1,2,3$. Using the Cauchy--Schwartz inequality and the support
properties of the functions $f_{k_i,j_i}$,
\begin{equation}\label{ela3}
\begin{split}
|J(f_{k_1,j_1},f_{k_2,j_2},f_{k_3,j_3})|&\leq
C2^{\mathrm{min}\,(j_1,j_2,j_3)/2}\int_{\mathbb{R}^2}A_{k_1}(\xi_1)A_{k_2}(\xi_2)
A_{k_3}(\xi_1+\xi_2)\,d\xi_1d\xi_2\\
&\leq
C2^{\mathrm{min}\,(k_1,k_2,k_3)/2}2^{\mathrm{min}\,(j_1,j_2,j_3)/2}
\prod_{i=1}^3||f_{k_i,j_i}||_{L^2},
\end{split}
\end{equation}
which is part (a).

For part (b) we observe that
\begin{equation}\label{om20}
|\Omega(\xi_1,\xi_2)|=2\,\mathrm{min}\,(|\xi_1|,|\xi_2|,|\xi_1+\xi_2|)\cdot
\mathrm{med}\,(|\xi_1|,|\xi_2|,|\xi_1+\xi_2|).
\end{equation}
Also, by examining the supports of  the functions, $J(f_{k_1,j_1},f_{k_2,j_2},f_{k_3,j_3})\equiv 0$ unless
\begin{equation}\label{new1}
\mathrm{max}\,(k_1,k_2,k_3)\leq \mathrm{med}\,(k_1,k_2,k_3)+2,
\end{equation}
and
\begin{equation}\label{new2}
\begin{cases}
&\mathrm{max}\,(j_1,j_2,j_3)\in[\widetilde{k}-5,\widetilde{k}+5]\,\text{ or}\\
&\mathrm{max}\,(j_1,j_2,j_3)\geq \widetilde{k}+5\text{ and }\mathrm{max}\,(j_1,j_2,j_3)-\mathrm{med}\,(j_1,j_2,j_3)\leq 5,
\end{cases}
\end{equation}
where $\widetilde{k}=\mathrm{min}\,(k_1,k_2,k_3)+\mathrm{med}\,(k_1,k_2,k_3)$.

Simple changes of variables and the observation that the function
$\omega$ is odd show that
\begin{equation}\label{ela1}
|J(f,g,h)|=|J(g,f,h)|\text{ and }|J(f,g,h)|=|J(\widetilde{f},h,g)|,
\end{equation}
where $\widetilde{f}(\xi,\mu)=f(-\xi,-\mu)$. Thus, by symmetry, in proving \eqref{om32} we
may assume $i=3$. Let
\begin{equation*}
B_{k_3}(\xi,\mu)=\Big[\frac{1}{2^{j_1}2^{j_2}}\int_{\mathbb{R}^2}
|f_{k_3,j_3}(\xi,\mu+\alpha+\beta)|^2(1+\alpha/2^{j_1})^{-2}(1+\beta/2^{j_2})^{-2}\,d\alpha
d\beta\Big]^{1/2}.
\end{equation*}
Clearly,
\begin{equation}\label{om6}
||B_{k_3}||_{L^2}=C||f_{k_3,j_3}||_{L^2}\text{ and }B_{k_3}\text{ is supported in }I_{k_3}\times\mathbb{R}.
\end{equation}
Also, by the Cauchy--Schwartz inequality,
\begin{equation}\label{om7}
\begin{split}
|J&(f_{k_1,j_1},f_{k_2,j_2},f_{k_3,j_3})|\\
&\leq
C2^{(j_1+j_2)/2}\int_{\mathbb{R}^2}A_{k_1}(\xi_1)A_{k_2}(\xi_2)B_{k_3}(\xi_1+\xi_2,\Omega(\xi_1,\xi_2))\,d\xi_1d\xi_2.
\end{split}
\end{equation}
We have three cases depending on the relative sizes of $|\xi_1|$,
$|\xi_2|$ and $|\xi_1+\xi_2|$. Let
\begin{equation*}
\begin{cases}
&R_1=\{(\xi_1,\xi_2):|\xi_1+\xi_2|\leq |\xi_1| \text{ and
}|\xi_2|\leq |\xi_1|\},\\
&R_2=\{(\xi_1,\xi_2):|\xi_1+\xi_2|\leq |\xi_2| \text{ and
}|\xi_1|\leq |\xi_2|\},\\
&R_3=\{(\xi_1,\xi_2):|\xi_1|\leq|\xi_1+\xi_2|\text{ and }|\xi_2|\leq|\xi_1+\xi_2|\}.
\end{cases}
\end{equation*}

For $(\xi_1,\xi_2)\in R_1$, using \eqref{om20},
$\Omega(\xi_1,\xi_2)=\pm 2\xi_2(\xi_1+\xi_2)$. We define
$B'_{k_3}(\xi,\mu)=B_{k_3}(\xi,2\xi\cdot\mu)$,
$||B'_{k_3}||_{L^2}\approx
2^{-k_3/2}||B_{k_3}||_{L^2}$. The integral over
$R_1$ in the right-hand side of \eqref{om7} is dominated
by
\begin{equation}\label{ela6}
\begin{split}
C\int_{\mathbb{R}^2}A_{k_1}(\xi_1)&A_{k_2}(\xi_2)[B'_{k_3}
(\xi_1+\xi_2,\xi_2)+B'_{k_3}
(\xi_1+\xi_2,-\xi_2)]\,d\xi_1d\xi_2\\
&\leq
C2^{-k_3/2}||A_{k_1}||_{L^2}||A_{k_2}||_{L^2}||B_{k_3}||_{L^2},
\end{split}
\end{equation}
which gives \eqref{om32} in this case (see \eqref{om6}).

The bound for the integral over $(\xi_1,\xi_2)\in R_2$ is
identical. We consider now the integral over
$(\xi_1,\xi_2)\in R_3$, in which case
$\Omega(\xi_1,\xi_2)=\pm 2\xi_1\xi_2$. By symmetry, to bound
the right-hand side of \eqref{om7} is suffices to bound
\begin{equation}\label{om10}
\int_{R_3}A_{k_1}(\xi_1)A_{k_2}(\xi_2)B_{k_3}(\xi_1+\xi_2,2\xi_1\xi_2)\,d\xi_1d\xi_2.
\end{equation}
We define $B''_{k_3}(\xi,\mu)=B_{k_3}(\xi,\mu+\xi^2/2)$,
so $||B''_{k_3}||_{L^2}=||B_{k_3}||_{L^2}$. Using \eqref{new1} and the assumption $\mathrm{max}\,(k_1,k_2,k_3)\geq \mathrm{min}\,(k_1,k_2,k_3)+5$, if $\xi_1\in I_{k_1}$, $\xi_2\in I_{k_2}$, $(\xi_1,\xi_2)\in R_3$, and $\xi_1+\xi_2\in I_{k_3}$, then $|\xi_1-\xi_2|\geq 2^{k_3-100}$. The integral in \eqref{om10} is dominated by
\begin{equation}\label{om11}
\begin{split}
\int_{\{|\xi_1-\xi_2|\geq
2^{k_3-100}\}}A_{k_1}(\xi_1)A_{k_2}(\xi_2)B''_{k_3}(\xi_1+\xi_2,-(\xi_1-\xi_2)^2/2)\,d\xi_1d\xi_2.
\end{split}
\end{equation}
Using the Cauchy--Schwartz inequality and a simple change of
variables, the integral in \eqref{om11} is dominated by
$C2^{-k_3/2}||A_{k_1}||_{L^2}||A_{k_2}||_{L^2}||B''_{k_3}||_{L^2}$, which completes the proof of \eqref{om32}.

For part (c), using part (a), we may assume
\begin{equation}\label{ela2}
\mathrm{med}\,(j_1,j_2,j_3)\leq 2\,\mathrm{min}\,(k_1,k_2,k_3).
\end{equation}
Using \eqref{ela1}, we may also assume $j_1=\mathrm{min}\,(j_1,j_2,j_3)$ and $j_2=\mathrm{med}\,(j_1,j_2,j_3)$. Let
\begin{equation*}
\widetilde{R}_{j_2}=\{(\xi_1,\xi_2):|\xi_1-\xi_2|\geq 2^{j_2/2}\}.
\end{equation*}
For the integral over $(\xi_1,\xi_2)\in{}^c\negmedspace\,\widetilde{R}_{j_2}=\mathbb{R}^2\setminus \widetilde{R}_{j_2}$ we use a bound similar to \eqref{ela3}:
\begin{equation*}
\begin{split}
\Big|\int_{{}^c\negmedspace\widetilde{R}_{j_2}\times\mathbb{R}^2}&f_{k_1,j_1}(\xi_1,\mu_1)f_{k_2,j_2}(\xi_2,\mu_2)f_{k_3,j_3}(\xi_1+\xi_2,\mu_1+\mu_2+\Omega(\xi_1,\xi_2))\,d\xi_1d\xi_2d\mu_1d\mu_2\Big|\\
&\leq C2^{j_1/2}\int_{{}^c\negmedspace\widetilde{R}_{j_2}}A_{k_1}(\xi_1)A_{k_2}(\xi_2)A_{k_3}(\xi_1+\xi_2)\,d\xi_1d\xi_2 \\
&\leq C2^{j_1/2}\iint_{|\mu|\leq 2^{j_2/2}}A_{k_1}(\xi_2+\mu)A_{k_2}(\xi_2)A_{k_3}(2\xi_2+\mu)\,d\xi_2d\mu\\
&\leq C2^{j_1/2}\int_{|\mu|\leq 2^{j_2/2}}\Big(\int_\mathbb{R}|A_{k_1}(\xi_2+\mu)|^2|A_{k_2}(\xi_2)|^2\,d\xi_2\Big)^{1/2}||A_{k_3}||_{L^2}d\mu\\
&\leq C2^{j_1/2}2^{j_2/4}||A_{k_1}||_{L^2}||A_{k_2}||_{L^2}||A_{k_3}||_{L^2},
\end{split}
\end{equation*}
which suffices for \eqref{om3}. For the integral over $(\xi_1,\xi_2)\in\widetilde{R}_{j_2}$ we use a bound similar to \eqref{om7}
\begin{equation}\label{ela4}
\begin{split}
&\Big|\int_{\widetilde{R}_{j_2}\times\mathbb{R}^2}f_{k_1,j_1}(\xi_1,\mu_1)f_{k_2,j_2}(\xi_2,\mu_2)f_{k_3,j_3}(\xi_1+\xi_2,\mu_1+\mu_2+\Omega(\xi_1,\xi_2))\,d\xi_1d\xi_2d\mu_1d\mu_2\Big|\\
&\leq
C2^{(j_1+j_2)/2}\int_{\widetilde{R}_{j_2}}A_{k_1}(\xi_1)A_{k_2}(\xi_2)B_{k_3}(\xi_1+\xi_2,\Omega(\xi_1,\xi_2))\,d\xi_1d\xi_2.
\end{split}
\end{equation}
We further decompose the integral in the right-hand side of \eqref{ela4} into three parts, corresponding to the regions $R_1$, $R_2$, and $R_3$. Using \eqref{ela6}, the integrals over the regions $\widetilde{R}_{j_2}\cap R_1$ and $\widetilde{R}_{j_2}\cap R_2$  are dominated by $C2^{-k_3/2}||A_{k_1}||_{L^2}||A_{k_2}||_{L^2}||B_{k_3}||_{L^2}$, which suffices in view of the assumption \eqref{ela2}. For the integral over the region $\widetilde{R}_{j_2}\cap R_3$, by symmetry it suffices to control
\begin{equation}\label{ela8}
\int_{\widetilde{R}_{j_2}\cap R_3}A_{k_1}(\xi_1)A_{k_2}(\xi_2)B_{k_3}(\xi_1+\xi_2,2\xi_1\xi_2)\,d\xi_1d\xi_2.
\end{equation}
As in the estimate of the integral in \eqref{om10}, the integral in \eqref{ela8} is dominated by
\begin{equation*}
\begin{split}
\int_{\{|\xi_1-\xi_2|\geq
2^{j_2}\}}A_{k_1}(\xi_1)A_{k_2}(\xi_2)B''_{k_3}(\xi_1+\xi_2,-(\xi_1-\xi_2)^2/2)\,d\xi_1d\xi_2.
\end{split}
\end{equation*}
The bound \eqref{om3} follows using the Cauchy--Schwartz inequality and a simple change of
variables.
\end{proof}

We restate now Lemma \ref{Lemmad1} in a form that is suitable for the bilinear estimates in the next sections.

\newtheorem{Lemmad2}[Lemmad1]{Corollary}
\begin{Lemmad2}\label{Lemmad2}
Assume $k_1,k_2,k_3\in\mathbb{Z}$,
$j_1,j_2,j_3\in\mathbb{Z}_+$, and $f_{k_i,j_i}\in
L^2(\mathbb{R}\times\mathbb{R})$ are functions supported in
$D_{k_i,j_i}$, $i=1,2$.

(a) For any  $k_1,k_2,k_3\in\mathbb{Z}$ and $j_1,j_2,j_3\in\mathbb{Z}_+$,
\begin{equation}\label{on31}
||\mathbf{1}_{D_{k_3,j_3}}(\xi,\tau)(f_{k_1,j_1}\ast f_{k_2,j_2})(\xi,\tau)||_{L^2}\leq C
2^{\mathrm{min}\,(k_1,k_2,k_3)/2}2^{\mathrm{min}\,(j_1,j_2,j_3)/2}
\prod_{i=1}^2||f_{k_i,j_i}||_{L^2}.
\end{equation}

(b) If $\mathrm{max}\,(k_1,k_2,k_3)\geq \mathrm{min}\,(k_1,k_2,k_3)+5$ and $i\in\{1,2,3\}$ then
\begin{equation}\label{on32}
||\mathbf{1}_{D_{k_3,j_3}}(\xi,\tau)(f_{k_1,j_1}\ast f_{k_2,j_2})(\xi,\tau)||_{L^2}\leq
C2^{(j_1+j_2+j_3)/2}2^{-(j_i+k_i)/2}\prod_{i=1}^2||f_{k_i,j_i}||_{L^2}.
\end{equation}

(c) For any  $k_1,k_2,k_3\in\mathbb{Z}$ and $j_1,j_2,j_3\in\mathbb{Z}_+$,
\begin{equation}\label{on3}
||\mathbf{1}_{D_{k_3,j_3}}(\xi,\tau)(f_{k_1,j_1}\ast f_{k_2,j_2})(\xi,\tau)||_{L^2}\leq
C2^{\mathrm{min}\,(j_1,j_2,j_3)/2+\mathrm{med}\,(j_1,j_2,j_3)/4}\prod_{i=1}^2||f_{k_i,j_i}||_{L^2}.
\end{equation}

(d) In addition, $\mathbf{1}_{D_{k_3,j_3}}(\xi,\tau)(f_{k_1,j_1}\ast f_{k_2,j_2})(\xi,\tau)\equiv 0$ unless
\begin{equation}\label{nev1}
\mathrm{max}\,(k_1,k_2,k_3)\leq \mathrm{med}\,(k_1,k_2,k_3)+2,
\end{equation}
and
\begin{equation}\label{nev2}
\begin{cases}
&\mathrm{max}\,(j_1,j_2,j_3)\in[\widetilde{k}-8,\widetilde{k}+8]\,\text{ or}\\
&\mathrm{max}\,(j_1,j_2,j_3)\geq \widetilde{k}+8\text{ and }\mathrm{max}\,(j_1,j_2,j_3)-\mathrm{med}\,(j_1,j_2,j_3)\leq 10,
\end{cases}
\end{equation}
where $\widetilde{k}=\mathrm{min}\,(k_1,k_2,k_3)+\mathrm{med}\,(k_1,k_2,k_3)$.
\end{Lemmad2}

\begin{proof}[Proof of Corollary \ref{Lemmad2}] Clearly,
\begin{equation*}
||\mathbf{1}_{D_{k_3,j_3}}(\xi,\tau)(f_{k_1,j_1}\ast f_{k_2,j_2})(\xi,\tau)||_{L^2}=\sup_{||f||_{L^2}=1}\Big|\int_{D_{k_3,j_3}}f\cdot(f_{k_1,j_1}\ast f_{k_2,j_2})\,d\xi d\tau\Big|.
\end{equation*}
Let $f_{k_3,j_3}=\mathbf{1}_{D_{k_3,j_3}}\cdot f$, and then $f^\#_{k_i,j_i}(\xi,\mu)=f_{k_i,j_i}(\xi,\mu+\omega(\xi))$, $i=1,2,3$. The functions $f^\#_{k_i,j_i}$ are supported in $I_{k_i}\times\bigcup_{|m|\leq 3}\widetilde{I}_{j_i+m}$, $||f^\#_{k_i,j_i}||_{L^2}=||f_{k_i,j_i}||_{L^2}$, and, using simple changes of variables,
\begin{equation*}
\int_{D_{k_3,j_3}}f\cdot(f_{k_1,j_1}\ast f_{k_2,j_2})\,d\xi d\tau=J(f^\#_{k_1,j_1},f^\#_{k_2,j_2},f^\#_{k_3,j_3}).
\end{equation*}
Corollary \ref{Lemmad2} follows from Lemma \ref{Lemmad1}, \eqref{new1}, and \eqref{new2}.
\end{proof}

\section{Bilinear estimates I}\label{bilinear}

In this section we prove two bilinear estimates, which correspond to $\mathrm{Low}\times\mathrm{High}\to\mathrm{High}$ interactions:

\newtheorem{Lemmac1}{Proposition}[section]
\begin{Lemmac1}\label{Lemmac1}
Assume $k\geq 20$, $k_2\in[k-2,k+2]$, $f_{k_2}\in Z_{k_2}$, and $f_{0}\in Z_{0}$. Then
\begin{equation}\label{bj1}
2^k\big|\big|\eta_k(\xi)\cdot(\tau-\omega(\xi)+i)^{-1}f_{k_2}\ast f_0\big|\big|_{Z_k}\leq C||f_{k_2}||_{Z_{k_2}}||f_{0}||_{Z_{0}}.
\end{equation}
\end{Lemmac1}

\newtheorem{Lemmac2}[Lemmac1]{Proposition}
\begin{Lemmac2}\label{Lemmac2}
Assume $k\geq 20$, $k_2\in[k-2,k+2]$, $f_{k_2}\in Z_{k_2}$, and $f_{k_1}\in Z_{k_1}$ for any $k_1\in[1,k-10]\cap\mathbb{Z}$. Then
\begin{equation}\label{bj2}
\begin{split}
2^k\big|\big|\eta_k(\xi)(\tau-&\omega(\xi)+i)^{-1}f_{k_2}\ast\sum_{k_1=1}^{k-10}f_{k_1}\big|\big|_{Z_k}\leq C||f_{k_2}||_{Z_{k_2}}\sup_{k_1\in[1,k-10]}||(I-\partial_\tau^2)f_{k_1}||_{Z_{k_1}}.
\end{split}
\end{equation}
\end{Lemmac2}

The main ingredients in the proofs of Propositions \ref{Lemmac1} and \ref{Lemmac2} are the definitions, the representations \eqref{repr1}, \eqref{repr2}, and \eqref{repr3}, Lemma \ref{Lemmaa1}, Lemma \ref{Lemmaa2} (b), (c), Corollary \ref{Lemmad2}, and the $L^2$ estimates in Lemma \ref{Lemmac3} below.

\newtheorem{Lemmac3}[Lemmac1]{Lemma}
\begin{Lemmac3}\label{Lemmac3}
Assume that $k\geq 20$, $k_1\in(-\infty,k-10]\cap\mathbb{Z}$, $k_2\in[k-2,k+2]$, $j,j_1,j_2\in\mathbb{Z}_+$, $f_{k_1,j_1}$ is an $L^2$ function supported in $D_{k_1,j_1}$, and $f_{k_2,j_2}$ is an $L^2$ function supported in $D_{k_2,j_2}$. Then, with $\gamma_{k,k_1}=(2^{k_1/2}+2^{-k/2})^{-1}$,
\begin{equation}\label{bo11}
\begin{split}
2^{k}2^{j/2}\beta_{k,j}||\eta_k(\xi)&\eta_j(\tau-\omega(\xi))(\tau-\omega(\xi)+i)^{-1}(f_{k_1,j_1}\ast f_{k_2,j_2})||_{L^2}\\
&\leq C\gamma_{k,k_1}\cdot2^{j_1/2}\beta_{k_1,j_1}||f_{k_1,j_1}||_{L^2}\cdot 2^{j_2/2}\beta_{k_2,j_2}||f_{k_2,j_2}||_{L^2},
\end{split}
\end{equation}
where, by definition, $\beta_{k_1,j_1}=2^{j_1/2}$ if $k_1\leq 0$. In addition, $\mathbf{1}_{D_{k,j}}(\xi,\tau)(f_{k_1,j_1}\ast f_{k_2,j_2})\equiv 0$ unless
\begin{equation}\label{bo5}
\begin{cases}
&\mathrm{max}\,(j,j_1,j_2)\in[k+k_1-10,k+k_1+10]\,\text{ or}\\
&\mathrm{max}\,(j,j_1,j_2)\geq k+k_1+10\text{ and }\mathrm{max}\,(j,j_1,j_2)-\mathrm{med}\,(j,j_1,j_2)\leq 10.
\end{cases}
\end{equation}
\end{Lemmac3}

{\bf{Remark:}} The bound \eqref{bo11} holds for $k_1$ both positive and negative. However, when $k_1\leq 0$, the right-hand side contains the large factor $\gamma_{k,k_1}$. This factor is the main reason why interactions between ``general'' $L^2$ functions of very low frequency and derivatives of $L^2$ functions of high frequency cannot be estimated using our bilinear estimates.

\begin{proof}[Proof of Lemma \ref{Lemmac3}] The restriction \eqref{bo5} follows directly from \eqref{nev2}. For \eqref{bo11} we use the bounds \eqref{on31}, \eqref{on32}, and \eqref{on3} in Corollary \ref{Lemmad2}. The left-hand side of \eqref{bo11} is dominated by
\begin{equation*}
2^{k}2^{-j/2}\beta_{k,j}||\mathbf{1}_{D_{k,j}}(\xi,\tau)(f_{k_1,j_1}\ast f_{k_2,j_2})||_{L^2}.
\end{equation*}
For \eqref{bo11} it suffices to prove that
\begin{equation}\label{bo4}
\begin{split}
||&\mathbf{1}_{D_{k,j}}(\xi,\tau)(f_{k_1,j_1}\ast f_{k_2,j_2})||_{L^2}\\
&\leq C2^{-k}\gamma_{k,k_1}2^{(j+j_1+j_2)/2}\beta_{k_1,j_1}\beta_{k_2,j_2}\beta_{k,j}^{-1}||f_{k_1,j_1}||_{L^2}||f_{k_2,j_2}||_{L^2}.
\end{split}
\end{equation}

Let $\Pi=||f_{k_1,j_1}||_{L^2}||f_{k_2,j_2}||_{L^2}$. We have several cases: if $j=\mathrm{max}\,(j,j_1,j_2)$ then, using \eqref{on32}, the left-hand side of \eqref{bo4} is dominated by $C2^{-k/2}2^{(j_1+j_2)/2}\Pi$; in addition $\beta_{k_1,j_1}\beta_{k_2,j_2}\beta_{k,j}^{-1}\geq C^{-1}$ and $2^{j/2}\geq C^{-1}(2^{(k+k_1)/2}+1)$, using \eqref{bo5}, so the bound \eqref{bo4} follows in this case.

If $j_2=\mathrm{max}\,(j,j_1,j_2)$ then, using \eqref{on32}, the left-hand side of \eqref{bo4} is dominated by $C2^{-k/2}2^{(j+j_1)/2}\Pi$; in addition $\beta_{k_1,j_1}\beta_{k_2,j_2}\beta_{k,j}^{-1}\geq C^{-1}$ and $2^{j_2/2}\geq C^{-1}(2^{(k+k_1)/2}+1)$, using \eqref{bo5}, so the bound \eqref{bo4} follows in this case.

If $j_1=\mathrm{max}\,(j,j_1,j_2)\geq k+k_1-20$ and $k_1\geq 0$ then, using \eqref{on32} and \eqref{on3}, the left-hand side of \eqref{bo4} is dominated by
\begin{equation*}
C2^{-j_1/2}(2^{k_1/2}+2^{\max(j,j_2)/4})^{-1}2^{(j+j_1+j_2)/2}\Pi;
\end{equation*}
in addition $2^{j_1/2}\beta_{k_1,j_1}\geq C^{-1}2^{j_1-k_1}$, $\beta_{k_2,j_2}\geq 1$, and $\beta_{k,j}\leq C\beta_{k,j_1}$. Using \eqref{bo5}, $2^{j_1}\beta_{k,j_1}^{-1}\geq C^{-1}2^{k+k_1}$, and the bound \eqref{bo4} follows. We notice also that the restriction $j_1=\mathrm{max}\,(j,j_1,j_2)$ was not important. For later use, we restate the stronger estimate we obtain in this case: if $k_1\geq 0$ and $j_1\geq k+k_1-20$ then
\begin{equation}\label{bo161}
\begin{split}
2^k&2^{j/2}\beta_{k,j}||\eta_k(\xi)\eta_j(\tau-\omega(\xi))(\tau-\omega(\xi)+i)^{-1}(f_{k_1,j_1}\ast f_{k_2,j_2})||_{L^2}\\
&\leq C(2^{k_1/2}+2^{\max(j,j_2)/4})^{-1}\cdot2^{j_1/2}\beta_{k_1,j_1}||f_{k_1,j_1}||_{L^2}\cdot 2^{j_2/2}\beta_{k_2,j_2}||f_{k_2,j_2}||_{L^2}.
\end{split}
\end{equation}

If $j_1=\mathrm{max}\,(j,j_1,j_2)\geq k+k_1-20$ and $k_1\leq 1$ then, using \eqref{on31}, the left-hand side of \eqref{bo4} is dominated by $C2^{k_1/2}2^{(j+j_2)/2}2^{-\max(j,j_2)/2}\Pi$; in addition $2^{j_1/2}\beta_{k_1,j_1}=2^{j_1}$, $\beta_{k_2,j_2}\geq 1$, and $\beta_{k,j}\leq C\beta_{k,j_1}$. Using \eqref{bo5}, $2^{j_1}\beta_{k,j_1}^{-1}\geq C^{-1}(2^{k+k_1}+1)$, and the bound \eqref{bo4} follows since $2^{k_1}+2^{-k}\geq C^{-1}\gamma_{k,k_1}^{-2}$. For later use, we restate the stronger estimate we obtain in this last case: if $k_1\leq 1$ and $j_1\geq k+k_1-20$ then
\begin{equation}\label{bo16}
\begin{split}
2^k2^{j/2}\beta_{k,j}&||\eta_k(\xi)\eta_j(\tau-\omega(\xi))(\tau-\omega(\xi)+i)^{-1}(f_{k_1,j_1}\ast f_{k_2,j_2})||_{L^2}\\
&\leq C2^{-\max(j,j_2)/2}\gamma_{k,k_1}\cdot 2^{j_1}||f_{k_1,j_1}||_{L^2}\cdot 2^{j_2/2}\beta_{k_2,j_2}||f_{k_2,j_2}||_{L^2}.
\end{split}
\end{equation}
\end{proof}

We prove now Proposition \ref{Lemmac1} and \ref{Lemmac2}.

\begin{proof}[Proof of Proposition \ref{Lemmac1}] We use the representations \eqref{repr1} and \eqref{repr2} and analyze three cases.

{\bf{Case 1:}} $f_0=f_{0,j_1}^{k_1}$ is supported in $D_{k_1,j_1}$, $f_{k_2}=f_{k_2,j_2}$ is supported in $D_{k_2,j_2}$, $j_1,j_2\geq  0$, $k_1\leq 1$, $||f_0||_{Z_0}\approx 2^{j_1-k_1}||f_{0,j_1}^{k_1}||_{L^2}$, and $||f_{k_2}||_{Z_{k_2}}\approx 2^{j_2/2}\beta_{k_2,j_2}||f_{k_2,j_2}||_{L^2}$. The bound \eqref{bj1} which we have to prove becomes
\begin{equation}\label{ht1}
2^k\big|\big|\eta_k(\xi)\cdot(\tau-\omega(\xi)+i)^{-1}f_{k_2,j_2}\ast f_{0,j_1}^{k_1}\big|\big|_{Z_k}\leq C2^{j_1-k_1}||f_{0,j_1}^{k_1}||_{L^2}\cdot2^{j_2/2}\beta_{k_2,j_2}||f_{k_2,j_2}||_{L^2}.
\end{equation}
Let $h_k(\xi,\tau)=\eta_k(\xi)(\tau-\omega(\xi)+i)^{-1}(f_{k_2,j_2}\ast f_{0,j_1}^{k_1})(\xi,\tau)$. The first observation is that for most choices of $j_1$ and $j_2$, depending on $k$ and $k_1$, the function $h_k$ is supported in a bounded number of regions $D_{k,j}$, so \eqref{bo11} suffices to control $2^k||h_k||_{X_k}$. In view of \eqref{bo5}, the function $h_k$ is supported in a bounded number of regions $D_{k,j}$, and \eqref{ht1} follows from \eqref{bo11}, unless
\begin{equation}\label{bo6}
\begin{cases}
&|j_1-(k+k_1)|\leq 10\text{ and }j_2\leq k+k_1+10\text{ or }\\
&|j_2-(k+k_1)|\leq 10\text{ and }j_1\leq k+k_1+10\text{ or }\\
&j_1,j_2\geq k+k_1-10\text{ and }|j_1-j_2|\leq 10.
\end{cases}
\end{equation}

Assume \eqref{bo6} holds. Using \eqref{bo5}, $\mathbf{1}_{D_{k,j}}(\xi,\tau)\cdot h_k\equiv 0$ unless $j\leq\max(j_1,j_2)+C$. We have two cases: if $j_1\geq k+k_1-20$, then, in view of \eqref{bo6}, $j_2\leq j_1+C$ and the function $h_k$ is supported in $\bigcup_{j\leq j_1+C}D_{k,j}$. By \eqref{bo16},
\begin{equation*}
\begin{split}
2^k&||h_k||_{X_k}\leq C2^k\sum_{j\leq j_1+C}2^{j/2}\beta_{k,j}||\eta_j(\tau-\omega(\xi))h_k(\xi,\tau)||_{L^2}\\
&\leq C\big[\sum_{j\leq j_1+C}2^{-\max(j,j_2)/2}\big]2^{-k_1/2}\cdot 2^{j_1}||f^{k_1}_{0,j_1}||_{L^2}\cdot 2^{j_2/2}\beta_{k_2,j_2}||f_{k_2,j_2}||_{L^2},
\end{split}
\end{equation*}
which suffices for \eqref{ht1}. Assume now that $j_1\leq k+k_1-20$, so, in view of \eqref{bo6}, $|j_2-(k+k_1)|\leq 10$ and the function $h_k$ is supported in $\bigcup_{j\leq k+k_1+C}D_{k,j}$. Then, using Lemma \ref{Lemmaa1} (b) and (c) (in fact the proof of part (b)),
\begin{equation*}
\begin{split}
2^{k}||h_k||_{Z_k}&\leq C2^{k/2}||\mathcal{F}^{-1}[(\tau-\omega(\xi)+i)h_k(\xi,\tau)]||_{L^1_xL^2_t}\\
&\leq C2^{k/2}||\mathcal{F}^{-1}(f_{0,j_1}^{k_1})||_{L^2_xL^\infty_t}||\mathcal{F}^{-1}(f_{k_2,j_2})||_{L^2_xL^2_t}\\
&\leq C2^{(j_1-k_1)/2}||f^{k_1}_{0,j_1}||_{L^2}\cdot 2^{(k+k_1)/2}||f_{k_2,j_2}||_{L^2},
\end{split}
\end{equation*}
which suffices for  \eqref{ht1} since $|j_2-(k+k_1)|\leq 10$. For later use we notice that we proved the slightly stronger estimate, with the factor $2^{-k_1}$ in the right-hand side of \eqref{ht1} replaced by $2^{-k_1/2}$,
\begin{equation}\label{ht1'}
2^k\big|\big|\eta_k(\xi)(\tau-\omega(\xi)+i)^{-1}f_{k_2,j_2}\ast f_{0,j_1}^{k_1}\big|\big|_{Z_k}\leq C2^{j_1-k_1/2}||f_{0,j_1}^{k_1}||_{L^2}\cdot2^{j_2/2}\beta_{k_2,j_2}||f_{k_2,j_2}||_{L^2}.
\end{equation}

{\bf{Case 2:}} $f_0=f_{0,j_1}^{k_1}$ is supported in $D_{k_1,j_1}$, $j_1\geq  0$, $k_1\leq 1$, $f_{k_2}=g_{k_2}$ is supported in $\bigcup_{j_2\leq k_2-1}D_{k_2,j_2}$, $||f_0||_{Z_0}\approx 2^{j_1-k_1}||f_{0,j_1}^{k_1}||_{L^2}$, and $||f_{k_2}||_{Z_{k_2}}\approx ||g_{k_2}||_{Y_{k_2}}$. The bound \eqref{bj1} which we have to prove becomes
\begin{equation}\label{ht2}
2^k\big|\big|\eta_k(\xi)\cdot(\tau-\omega(\xi)+i)^{-1}g_{k_2}\ast f_{0,j_1}^{k_1}\big|\big|_{Z_k}\leq C2^{j_1-k_1}||f_{0,j_1}^{k_1}||_{L^2}\cdot||g_{k_2}||_{Y_{k_2}}.
\end{equation}
As before, let $h_k(\xi,\tau)=\eta_k(\xi)(\tau-\omega(\xi)+i)^{-1}(g_{k_2}\ast f_{0,j_1}^{k_1})(\xi,\tau)$. In view of Lemma \ref{Lemmaa1} (b), (c), and the bound \eqref{ht1'}, we may assume that $g_{k_2}$ is supported in the set $\{(\xi_2,\tau_2):\xi_2\in I_{k_2},\,|\tau_2-\omega(\xi_2)|\leq 2^{k+k_1-20}\}$. We have two cases: if $j_1\geq k+k_1-20$ then let $g_{k_2,j_2}(\xi_2,\tau_2)=g_{k_2}(\xi_2,\tau_2)\eta_{j_2}(\tau_2-\omega(\xi_2))$. Using $X_k$ norms, Lemma \ref{Lemmaa1} (b), and \eqref{bo16}, the left-hand side of \eqref{ht2} is dominated by
\begin{equation*}
\begin{split}
&C\sum_{j,j_2\leq j_1+C}2^{k}2^{j/2}\beta_{k,j}||\eta_k(\xi)\eta_j(\tau-\omega(\xi))(\tau-\omega(\xi)+i)^{-1}(f_{0,j_1}^{k_1}\ast g_{k_2,j_2})||_{L^2}\\
&\leq C\gamma_{k,k_1}\cdot2^{j_1}||f^{k_1}_{0,j_1}||_{L^2}\sum_{j,j_2\leq j_1+C}2^{-\max(j,j_2)/2}\cdot 2^{j_2/2}\beta_{k_2,j_2}||g_{k_2,j_2}||_{L^2}\\
&\leq C\gamma_{k,k_1}\cdot2^{j_1}||f_{0,j_1}^{k_1}||_{L^2}\cdot ||g_{k_2}||_{Y_{k_2}},
\end{split}
\end{equation*}
which suffices to prove \eqref{ht2} in this case. Assume now that $j_1\leq k+k_1-20$. In view of \eqref{bo5}, the function in the left-hand side of \eqref{ht2} is supported in the union of a bounded number of dyadic regions $D_{k,j}$, $|j-(k+k_1)|\leq C$. Then, using $X_k$ norms in the left-hand side of \eqref{ht2} and Lemma \ref{Lemmaa2} (c), the left-hand side of \eqref{ht2} is dominated by
\begin{equation*}
\begin{split}
C2^k2^{-(k+k_1)/2}||f^{k_1}_{0,j_1}\ast g_{k_2}||_{L^2}&\leq C2^{(k-k_1)/2}||\mathcal{F}^{-1}(f^{k_1}_{0,j_1})||_{L^2_xL^\infty_t}||\mathcal{F}^{-1}(g_{k_2})||_{L^\infty_xL^2_t}\\
&\leq C2^{(k-k_1)/2}\cdot 2^{j_1/2}||f^{k_1}_{0,j_1}||_{L^2}\cdot 2^{-k/2}||g_{k_2}||_{Y_{k_2}}\\
&\leq C2^{(j_1-k_1)/2}||f^{k_1}_{0,j_1}||_{L^2}\cdot||g_{k_2}||_{Y_{k_2}},
\end{split}
\end{equation*}
which completes the proof of \eqref{ht2}.

{\bf{Case 3:}} $f_0=g_{0,j}$ is supported in $\widetilde{I}_0\times\widetilde{I}_{j_1}$, $j_1\geq 0$, $||f_0||_{Z_0}\approx 2^{j_1}||\mathcal{F}^{-1}(g_{0,j_1})||_{L^1_xL^2_t}$. The bound \eqref{bj1} which we have to prove becomes
\begin{equation}\label{mu1}
2^k\big|\big|\eta_k(\xi)\cdot(\tau-\omega(\xi)+i)^{-1}f_{k_2}\ast g_{0,j_1}\big|\big|_{Z_k}\leq C2^{j_1}||\mathcal{F}^{-1}(g_{0,j_1})||_{L^1_xL^2_t}\cdot||f_{k_2}||_{Z_{k_2}}.
\end{equation}
Using the representation \eqref{ar209}, we see easily that
\begin{equation}\label{mu2}
||\mathcal{F}^{-1}(g_{0,j_1})||_{L^1_xL^\infty_t}+||\mathcal{F}^{-1}(g_{0,j_1})||_{L^2_xL^\infty_t}\leq C2^{j_1/2}||\mathcal{F}^{-1}(g_{0,j_1})||_{L^1_xL^2_t}.
\end{equation}
Thus, using the definitions, Lemma \ref{Lemmaa1} (b), (c), and Lemma \ref{Lemmaa2} (c),
\begin{equation*}
\begin{split}
2^k\big|\big|\eta_k(\xi)&\eta_{\leq k+C}(\tau-\omega(\xi))(\tau-\omega(\xi)+i)^{-1}f_{k_2}\ast g_{0,j_1}\big|\big|_{Z_k}\\
&\leq C2^{k/2}||\mathcal{F}^{-1}(f_{k_2}\ast g_{0,j_1})||_{L^1_xL^2_t}\\
&\leq C2^{k/2}||\mathcal{F}^{-1}(f_{k_2})||_{L^\infty_xL^2_t}||\mathcal{F}^{-1}(g_{0,j_1})||_{L^1_xL^\infty_t}\\
&\leq C2^{j_1/2}||\mathcal{F}^{-1}(g_{0,j_1})||_{L^1_xL^2_t}\cdot||f_{k_2}||_{Z_{k_2}}.
\end{split}
\end{equation*}
Thus, for \eqref{mu1}, it suffices to prove that
\begin{equation}\label{mu3}
2^k\sum_{j\geq k+C}2^{-j/2}\beta_{k,j}\big|\big|\eta_k(\xi)\eta_{j}(\tau-\omega(\xi))f_{k_2}\ast g_{0,j_1}\big|\big|_{L^2}\leq C2^{j_1}||\mathcal{F}^{-1}(g_{0,j_1})||_{L^1_xL^2_t}\cdot||f_{k_2}||_{Z_{k_2}}.
\end{equation}
Using Lemma \ref{Lemmaa2} (c) and \eqref{mu2} again,
\begin{equation*}
\begin{split}
\big|\big|\eta_k(\xi)\eta_{j}(\tau-\omega(\xi))f_{k_2}\ast g_{0,j_1}\big|\big|_{L^2}&\leq C||\mathcal{F}^{-1}(f_{k_2})||_{L^\infty_xL^2_t}||\mathcal{F}^{-1}(g_{0,j_1})||_{L^2_xL^\infty_t}\\
&\leq C2^{j_1/2}||\mathcal{F}^{-1}(g_{0,j_1})||_{L^1_xL^2_t}\cdot 2^{-k/2}||f_{k_2}||_{Z_{k_2}}.
\end{split}
\end{equation*}
We use this bound to control the sum over $j\leq 2k+j_1+C$ in \eqref{mu3}. For $j\geq 2k+j_1+C$, $2^{-j/2}\beta_{k,j}\approx 2^{-k}$, and, for \eqref{mu1}, it suffices to prove that
\begin{equation}\label{mu4}
\sum_{j\geq 2k+j_1+C}\big|\big|\eta_k(\xi)\eta_{j}(\tau-\omega(\xi))f_{k_2}\ast g_{0,j_1}\big|\big|_{L^2}\leq C2^{j_1}||\mathcal{F}^{-1}(g_{0,j_1})||_{L^1_xL^2_t}\cdot||f_{k_2}||_{Z_{k_2}}.
\end{equation}
By examining the supports of the functions, $\eta_k(\xi)\eta_{j}(\tau-\omega(\xi))f_{k_2}\ast g_{0,j_1}\equiv 0$ if $f_{k_2}\in Y_{k_2}$ and $j\geq 2k+j_1+C$. So, in \eqref{mu4}, we may assume $f_{k_2}=f_{k_2,j_2}$ is supported in $D_{k_2,j_2}$, $j_2\geq 2k+j_1+C$. The sum in $j$ in \eqref{mu4} is taken over $|j-j_2|\leq C$. Using Lemma \ref{Lemmaa2} (c) and \eqref{mu2}, the left-hand side of \eqref{mu4} is dominated by
\begin{equation*}
\begin{split}
C||\mathcal{F}^{-1}(f_{k_2,j_2})||_{L^\infty_xL^2_t}||\mathcal{F}^{-1}(g_{0,j_1})||_{L^2_xL^\infty_t}\leq C2^{j_1/2}||\mathcal{F}^{-1}(g_{0,j_1})||_{L^1_xL^2_t}\cdot ||f_{k_2,j_2}||_{Z_{k_2}}.
\end{split}
\end{equation*}
This completes the proof of \eqref{mu4} and \eqref{mu1}.
\end{proof}

For later use, we  notice that a simplified version of our argument can be used to prove the following: if $k\geq 20$, $k_2\in[k-2,k+2]$, $f_{k_2}\in Z_{k_2}$, and $f_{0}\in \overline{Z}_{0}$, then
\begin{equation}\label{hh91}
\big|\big|\eta_k(\xi)\cdot(\tau-\omega(\xi)+i)^{-1}f_{k_2}\ast f_0\big|\big|_{Z_k}\leq C||f_{k_2}||_{Z_{k_2}}||f_{0}||_{\overline{Z}_{0}}.
\end{equation}
To prove \eqref{hh91}, we use Lemma \ref{Lemmaa1} (b) to bound $||f_{k_2}||_{Z_{k_2}}\geq C^{-1}k^{-1}||f_{k_2}||_{X_{k_2}}$. Then, we write $f_0=\sum_{j_1\geq 0}\sum_{k_1\leq 1} f_{k_1,j_1}$, $f_{k_1,j_1}$ supported in $D_{k_1,j_1}$ and $||f_0||_{\overline{Z}_0}\geq \sum_{j_1\geq 0}\sum_{k_1\leq 1}2^{j_1}2^{k_1/4}||f_{k_1,j_1}||_{L^2}$. In view of the definitions, for \eqref{hh91} it suffices to prove that if $f_{k_2,j_2}$ is supported in $D_{k_2,j_2}$ then
\begin{equation*}
\sum_{j}2^{-j/2}\beta_{k,j}||\mathbf{1}_{D_{k,j}}f_{k_2,j_2}\ast f_{k_1,j_1}||_{L^2}\leq Ck^{-1}2^{j_2/2}\beta_{k_2,j_2}||f_{k_2,j_2}||_{L^2}\cdot 2^{j_1}2^{k_1/4}||f_{k_1,j_1}||_{L^2}.
\end{equation*}
Using \eqref{on31} we bound $||\mathbf{1}_{D_{k,j}}f_{k_2,j_2}\ast f_{k_1,j_1}||_{L^2}\leq C2^{k_1/2}2^{j_1/2}||f_{k_2,j_2}||_{L^2}\cdot ||f_{k_1,j_1}||_{L^2}$. So, it suffices to prove that
\begin{equation*}
2^{k_1/4}k\sum_{j}2^{-j/2}\beta_{k,j}\leq C2^{(j_1+j_2)/2},
\end{equation*}
where the sum is taken over $j$ satisfying \eqref{bo5}. This follows easily by examining the cases $\max(j_1,j_2)\leq k+k_1-20$ and $\max(j_1,j_2)\geq k+k_1-20$ (in the second case we estimate $2^{-j/2}\beta_{k,j}\leq C$).

\begin{proof}[Proof of Proposition \ref{Lemmac2}] The proof is similar to the proof of Proposition \ref{Lemmac1}, with an additional technical difficulty related to the sum in $k_1$ in the left-hand side of \eqref{bj2}. Our main tools are the bounds \eqref{bo11} (with $\gamma_{k,k_1}\approx 2^{-k_1/2}$  if $k_1\geq 1$) and \eqref{bo161}. For any $k_1\in[1,k-10]$ we decompose
\begin{equation}\label{vb1}
f_{k_1}=f_{k_1}^h+f_{k_1}^l=f_{k_1}\cdot[1-\eta_{\leq k+k_1-20}(\tau-\omega(\xi))]+f_{k_1}\cdot\eta_{\leq k+k_1-20}(\tau-\omega(\xi)).
\end{equation}
We show first that
\begin{equation}\label{vb2}
\begin{split}
2^k\big|\big|\eta_k(\xi)(\tau-&\omega(\xi)+i)^{-1}f_{k_2}\ast f_{k_1}^h\big|\big|_{X_k}\leq C2^{-k_1/4}||f_{k_2}||_{Z_{k_2}}||f_{k_1}^h||_{Z_{k_1}}.
\end{split}
\end{equation}
Assuming \eqref{vb2}, we can use the factor $2^{-k_1/4}$ to sum in $k_1$ and obtain
\begin{equation}\label{vb3}
\begin{split}
2^k\big|\big|\eta_k(\xi)(\tau-&\omega(\xi)+i)^{-1}f_{k_2}\ast\sum_{k_1=1}^{k-10}f_{k_1}^h\big|\big|_{X_k}\leq C||f_{k_2}||_{Z_{k_2}}\sup_{k_1\in[1,k-10]}||f_{k_1}||_{Z_{k_1}}.
\end{split}
\end{equation}
To prove \eqref{vb2} we use the representation \eqref{repr1} and \eqref{bo161}. We may assume $f_{k_1}^h=f_{k_1,j_1}$ is supported in $D_{k_1,j_1}$, $j_1\geq k+k_1-20$, $||f_{k_1}^h||_{Z_{k_1}}\approx 2^{j_1/2}\beta_{k_1,j_1}||f_{k_1,j_1}||_{L^2}$. We have two cases: if $f_{k_2}=f_{k_2,j_2}$ is supported in $D_{k_2,j_2}$, $j_2\geq 0$, $||f_{k_2}||_{Z_{k_2}}\approx 2^{j_2/2}\beta_{k_2,j_2}||f_{k_2,j_2}||_{L^2}$, then, using \eqref{bo161} and the definitions, the left-hand side of \eqref{vb2} is dominated by
\begin{equation*}
\begin{split}
C\big[\sum_j(2^{k_1/2}+2^{j/4})^{-1}\big]&\cdot2^{j_1/2}\beta_{k_1,j_1}||f_{k_1,j_1}||_{L^2}\cdot 2^{j_2/2}\beta_{k_2,j_2}||f_{k_2,j_2}||_{L^2}\\
&\leq C2^{-k_1/4}\cdot2^{j_1/2}\beta_{k_1,j_1}||f_{k_1,j_1}||_{L^2}\cdot 2^{j_2/2}\beta_{k_2,j_2}||f_{k_2,j_2}||_{L^2},
\end{split}
\end{equation*}
which gives \eqref{vb2} in this case. If $f_{k_2}=g_{k_2}$ is supported in $\bigcup_{j_2\leq k_2-1}D_{k_2,j_2}$, $||f_{k_2}||_{Z_{k_2}}\approx ||g_{k_2}||_{Y_{k_2}}$, then let $g_{k_2,j_2}(\xi_2,\tau_2)=g_{k_2}(\xi_2,\tau_2)\eta_{j_2}(\tau_2-\omega(\xi_2))$. In view of Lemma \ref{Lemmaa1} (b), \eqref{bo5}, and \eqref{bo161}, the left-hand side of \eqref{vb2} is dominated by
\begin{equation*}
\begin{split}
&C\sum_{j,j_2\leq j_1+C}2^{k}2^{j/2}\beta_{k,j}||\eta_k(\xi)\eta_j(\tau-\omega(\xi))(\tau-\omega(\xi)+i)^{-1}(g_{k_2,j_2}\ast f_{k_1,j_1})||_{L^2}\\
&\leq C2^{j_1/2}\beta_{k_1,j_1}||f_{k_1,j_1}||_{L^2}\sum_{j,j_2\leq j_1+C}(2^{k_1/2}+2^{\max(j,j_2)/4})^{-1}2^{j_2/2}||g_{k_2,j_2}||_{L^2}\\
&\leq C2^{-k_1/4}\cdot2^{j_1/2}\beta_{k_1,j_1}||f_{k_1,j_1}||_{L^2}\cdot ||g_{k_2}||_{Y_{k_2}},
\end{split}
\end{equation*}
which completes the proof of \eqref{vb2}.

In view of \eqref{vb3}, for \eqref{bj2} it suffices to prove that
\begin{equation}\label{vt5}
2^k\big|\big|\eta_k(\xi)(\tau-\omega(\xi)+i)^{-1}f_{k_2}\ast\sum_{k_1=1}^{k-10}f_{k_1}^l\big|\big|_{Z_k}\leq C||f_{k_2}||_{Z_{k_2}}\sup_{k_1\in[1,k-10]}||(I-\partial_\tau^2)f_{k_1}^l||_{Z_{k_1}}
\end{equation}
for any functions $f_{k_1}^l$ supported in $\bigcup_{j_1\leq k+k_1-19}{D_{k_1,j_1}}$. Using the representation \eqref{repr1}, we analyze two cases.

{\bf{Case 1:}} $f_{k_2}=f_{k_2,j_2}$ is supported in $D_{k_2,j_2}$, $||f_{k_2}||_{Z_{k_2}}\approx 2^{j_2/2}\beta_{k_2,j_2}||f_{k_2,j_2}||_{L^2}$, $j_2\geq  0$. The bound \eqref{vt5} which we have to prove becomes
\begin{equation}\label{vt6}
\begin{split}
2^k\big|\big|&\eta_k(\xi)(\tau-\omega(\xi)+i)^{-1}f_{k_2,j_2}\ast\sum_{k_1=1}^{k-10}f_{k_1}^l\big|\big|_{Z_k}\\
&\leq C2^{j_2/2}\beta_{k_2,j_2}||f_{k_2,j_2}||_{L^2}\sup_{k_1\in[1,k-10]}||(I-\partial_\tau^2)f_{k_1}^l||_{Z_{k_1}},
\end{split}
\end{equation}
for any functions $f_{k_1}^l$ supported in $\bigcup_{j_1\leq k+k_1-19}{D_{k_1,j_1}}$. Notice that $j_2$ is fixed in \eqref{vt6}. We divide the set of indices $k_1$ into two sets:
\begin{equation*}
\begin{cases}
&A_{k,j_2}=\{k_1\in[1,k-10]:|k+k_1-j_2|\leq 15\};\\
&B_{k,j_2}=\{k_1\in[1,k-10]:|k+k_1-j_2|\geq 16\}.
\end{cases}
\end{equation*}
The set $A_{k,j_2}$ has at most $31$ elements, so, for \eqref{vt6} it suffices to prove that
\begin{equation}\label{vt7}
\begin{split}
2^k\big|\big|&\eta_k(\xi)(\tau-\omega(\xi)+i)^{-1}f_{k_2,j_2}\ast f_{k_1}^l\big|\big|_{Z_k}\leq C2^{j_2/2}\beta_{k_2,j_2}||f_{k_2,j_2}||_{L^2}||(I-\partial_\tau^2)f_{k_1}^l||_{Z_{k_1}}
\end{split}
\end{equation}
for $k_1\in A_{k,j_2}$, and
\begin{equation}\label{vt8}
\begin{split}
2^k\big|\big|&\eta_k(\xi)(\tau-\omega(\xi)+i)^{-1}f_{k_2,j_2}\ast f_{k_1}^l\big|\big|_{X_k}\leq C2^{-k_1/4}\cdot2^{j_2/2}\beta_{k_2,j_2}||f_{k_2,j_2}||_{L^2}||f_{k_1}^l||_{Z_{k_1}}
\end{split}
\end{equation}
for $k_1\in B_{k,j_2}$.

We prove first \eqref{vt7}. In view of the restriction on the support of $f_{k_1}^l$, the condition $k_1\in A_{k,j_2}$, and \eqref{bo5}, the function $\eta_k(\xi)(\tau-\omega(\xi)+i)^{-1}f_{k_2,j_2}\ast f_{k_1}^l$  is supported in $\bigcup_{j\leq k+k_1+C}D_{k,j}$. In view  of the definition of the space $Z_k$, for \eqref{vt7} it suffices to prove that
\begin{equation}\label{vt9}
\begin{split}
2^k\big|\big|&\eta_{\leq k-1}(\tau-\omega(\xi))\eta_k(\xi)(\tau-\omega(\xi)+i)^{-1}f_{k_2,j_2}\ast f_{k_1}^l\big|\big|_{Y_k}\\
&\leq C2^{j_2/2}\beta_{k_2,j_2}||f_{k_2,j_2}||_{L^2}||(I-\partial_\tau^2)f_{k_1}^l||_{Z_{k_1}}
\end{split}
\end{equation}
and
\begin{equation}\label{vt10}
\begin{split}
2^k\sum_{j=k}^{k+k_1+C}2^{j/2}\beta_{k,j}\big|\big|&\eta_{j}(\tau-\omega(\xi))\eta_k(\xi)(\tau-\omega(\xi)+i)^{-1}f_{k_2,j_2}\ast f_{k_1}^l\big|\big|_{L^2}\\
&\leq C2^{j_2/2}\beta_{k_2,j_2}||f_{k_2,j_2}||_{L^2}||f_{k_1}^l||_{Z_{k_1}}.
\end{split}
\end{equation}
For \eqref{vt9} we use Lemma \ref{Lemmaa1} (a), (c), and Lemma \ref{Lemmaa2} (b). Since $|k+k_1-j_2|\leq 10$, the left-hand side of \eqref{vt9} is dominated by
\begin{equation*}
\begin{split}
C2^{k/2}||\mathcal{F}^{-1}(f_{k_2,j_2}\ast f_{k_1}^l)||_{L^1_xL^2_t}&\leq C2^{k/2}||\mathcal{F}^{-1}(f_{k_2,j_2})||_{L^2}||\mathcal{F}^{-1}(f_{k_1}^l)||_{L^2_xL^\infty_t}\\
&\leq C2^{k/2}||f_{k_2,j_2}||_{L^2}\cdot 2^{k_1/2}||(I-\partial_\tau^2)f_{k_1}^l||_{Z_{k_1}},
\end{split}
\end{equation*}
which completes the proof of \eqref{vt9}. For \eqref{vt10}, we  notice that the sum in the left-hand side contains at most $k_1+C$ terms. In addition, using Lemma \ref{Lemmaa1} (b), $||f_{k_1}^l||_{Z_{k_1}}\geq Ck_1^{-1}||f_{k_1}^l||_{X_{k_1}}$, and, using \eqref{bo11}, for any $j\in[k,k+k_1+C]$
\begin{equation*}
\begin{split}
2^k2^{j/2}\beta_{k,j}\big|\big|&\eta_{j}(\tau-\omega(\xi))\eta_k(\xi)(\tau-\omega(\xi)+i)^{-1}f_{k_2,j_2}\ast f_{k_1}^l\big|\big|_{L^2}\\
&\leq C2^{-k_1/2}2^{j_2/2}\beta_{k_2,j_2}||f_{k_2,j_2}||_{L^2}\cdot||f_{k_1}^l||_{X_{k_1}}.
\end{split}
\end{equation*}
This completes the proof of \eqref{vt10} and \eqref{vt7}.

We prove now the bound \eqref{vt8}. The main observation is that the function $\eta_k(\xi)(\tau-\omega(\xi)+i)^{-1}f_{k_2,j_2}\ast f_{k_1}^l$ is  supported in a bounded number of regions $D_{k,j}$ (assuming $j_2$ and $k_1$ fixed). This  is due to the support property of the function $f_{k_1}^l$, the assumption $|k+k_1-j_2|\geq 16$, and \eqref{bo5}. Thus, using \eqref{bo11}, the  left hand side of  \eqref{vt8} is dominated by
\begin{equation*}
\begin{split}
C\sup_{j}2^k2^{j/2}\beta_{k,j}\big|\big|&\eta_{j}(\tau-\omega(\xi))\eta_k(\xi)(\tau-\omega(\xi)+i)^{-1}f_{k_2,j_2}\ast f_{k_1}^l\big|\big|_{L^2}\\
&\leq C2^{-k_1/2}2^{j_2/2}\beta_{k_2,j_2}||f_{k_2,j_2}||_{L^2}\cdot||f_{k_1}^l||_{X_{k_1}},
\end{split}
\end{equation*}
which suffices  for \eqref{vt8} since $||f_{k_1}^l||_{Z_{k_1}}\geq Ck_1^{-1}||f_{k_1}^l||_{X_{k_1}}$ (see Lemma \ref{Lemmaa1} (b)).

{\bf{Case 2:}} $f_{k_2}=g_{k_2}$ is supported in $\bigcup_{j_2\leq k_2-20}D_{k_2,j_2}$, $||f_{k_2}||_{Z_{k_2}}\approx ||g_{k_2}||_{Y_{k_2}}$. The bound \eqref{vt5} which we have to prove becomes
\begin{equation*}
\begin{split}
2^k\big|\big|&\eta_k(\xi)(\tau-\omega(\xi)+i)^{-1}g_{k_2}\ast\sum_{k_1=1}^{k-10}f_{k_1}^l\big|\big|_{Z_k}\leq C||g_{k_2}||_{Y_{k_2}}\sup_{k_1\in[1,k-10]}||(I-\partial_\tau^2)f_{k_1}^l||_{Z_{k_1}}
\end{split}
\end{equation*}
for any functions $f_{k_1}^l$ supported in $\bigcup_{j_1\leq k+k_1-19}{D_{k_1,j_1}}$. Using Lemma \ref{Lemmaa1} (b)  again, it suffices to  prove that
\begin{equation}\label{vt12}
\begin{split}
2^k\big|\big|&\eta_k(\xi)(\tau-\omega(\xi)+i)^{-1}g_{k_2}\ast f_{k_1}^l\big|\big|_{X_k}\leq C2^{-k_1/2}||g_{k_2}||_{Y_{k_2}}||f_{k_1}^l||_{X_{k_1}}.
\end{split}
\end{equation}
Using \eqref{bo5} and the support properties of $g_{k_2}$ and $f_{k_1}^l$, the function $\eta_k(\xi)(\tau-\omega(\xi)+i)^{-1}g_{k_2}\ast f_{k_1}^l$  is supported in a bounded number of regions $D_{k,j}$, $|k+k_1-j|\leq  C$. Thus, for \eqref{vt12} it suffices to prove that if $f_{k_1,j_1}$ is supported in $D_{k_1,j_1}$, $j_1\leq k+k_1-19$ and $|j-k-k_1|\leq C$, then
\begin{equation}\label{vt13}
\begin{split}
2^{k/2}\big|\big|\mathbf{1}_{D_{k,j}}\cdot (g_{k_2}\ast f_{k_1,j_1})\big|\big|_{L^2}\leq C||g_{k_2}||_{Y_{k_2}}\cdot 2^{j_1/2}||f_{k_1,j_1}||_{L^2}.
\end{split}
\end{equation}

To prove \eqref{vt13} we may assume $k_2\geq 100$. For $j_2\leq k_2$ let $g_{k_2,j_2}(\xi,\tau)=\eta_{j_2}(\tau-\omega(\xi))g_{k_2}(\xi,\tau)$. Notice that in view of \eqref{on32} and Lemma \ref{Lemmaa1} (b)
\begin{equation}\label{vt14}
\begin{split}
2^{k/2}\big|\big|\mathbf{1}_{D_{k,j}}\cdot (g_{k_2,j_2}\ast f_{k_1,j_1})\big|\big|_{L^2}&\leq C2^{j_2/2}||g_{k_2,j_2}||_{L^2}\cdot 2^{j_1/2}||f_{k_1,j_1}||_{L^2}\\
&\leq C||g_{k_2}||_{Y_{k_2}}\cdot 2^{j_1/2}||f_{k_1,j_1}||_{L^2},
\end{split}
\end{equation}
for any $j_2\leq k_2$. To prove \eqref{vt13} we have to avoid the logarithmic divergence that appears when summing the bound above over $j_2\leq k_2$. In view of \eqref{pr10}, we may assume
\begin{equation}\label{vt15}
\begin{cases}
&g_{k_2}(\xi,\tau)=2^{k_2/2}\chi_{[k_2-1,k_2+1]}(\xi)(\tau-\omega(\xi)+i)^{-1}
\eta_{\leq k_2}(\tau-\omega(\xi))h(\tau);\\
&||g_{k_2}||_{Y_{k_2}}= C||h||_{L^2_\tau}.
\end{cases}
\end{equation}
We argue as  in the proof of Lemma \ref{Lemmaa2} (b). Let $h_+=h\cdot\mathbf{1}_{[0,\infty)}$, $h_-=h\cdot\mathbf{1}_{(-\infty,0]}$, and define the corresponding functions $g_{k_2,+}$ and $g_{k_2,-}$ as in \eqref{vt15}. By symmetry, it suffices to prove the bound \eqref{vt13} for the function $g_{k_2,+}$, which is supported in the set $\{(\xi,\tau):\xi\in[-2^{k_2+2},-2^{k_2-2}],\,\tau\in[2^{2k_2-10},2^{2k_2+10}]\}$. In view of \eqref{omega}, $\tau-\omega(\xi)=\tau-\xi^2$ on the support of $g_{k_2,+}$, and $g_{k_2,+}(\xi,\tau)=0$ unless $|\sqrt\tau+\xi|\leq C$. Let
\begin{equation}\label{vt16}
\begin{split}
g'_{k_2,+}(\xi,\tau)=2^{k_2/2}\chi_{[k_2-1,k_2+1]}(-\sqrt\tau)(\tau-&\xi^2+(\sqrt\tau+\xi)^2+i\sqrt\tau2^{-k_2})^{-1}\\
&\eta_{0}(\sqrt\tau+\xi)\cdot h_+(\tau).
\end{split}
\end{equation}
Using Lemma \ref{Lemmaa1} (b), it is easy to see that $||g_{k_2,+}-g'_{k_2,+}||_{X_{k_2}}\leq C||h_+||_{L^2}$. In view of \eqref{vt14}, for \eqref{vt13} it suffices to prove that
\begin{equation}\label{vt20}
\begin{split}
2^{k/2}\big|\big|g'_{k_2,+}\ast f_{k_1,j_1}\big|\big|_{L^2}\leq C||h_+||_{L^2}\cdot 2^{j_1/2}||f_{k_1,j_1}||_{L^2}.
\end{split}
\end{equation}
We substitute the formula \eqref{vt16} and make the change of variables $\xi_2=-\sqrt{\tau_2}+\mu_2$. The left-hand side of \eqref{vt20} is dominated by
\begin{equation*}
\begin{split}
\Big|\Big|\int_{\mathbb{R}^2}f_{k_1,j_1}(\xi+&\sqrt{\tau_2}-\mu_2, \tau-\tau_2)\cdot\eta_0(\mu_2)\frac{1}{\mu_2+i/2^{k_2+1}}\cdot h'_+(\tau_2)\,d\mu_2d\tau_2\Big|\Big|_{L^2_{\xi,\tau}},
\end{split}
\end{equation*}
where $h'_+(\tau_2)=h_+(\tau_2)\chi_{[k_2-1,k_2+1]}(-\sqrt{\tau_2})(2^k/\sqrt{\tau_2})$ is supported in $[2^{2k_2-4},2^{2k_2+4}]$, $||h'_+||_{L^2}\approx ||h_+||_{L^2}$. By duality, for \eqref{vt20} it suffices to prove that for any $m\in L^2$
\begin{equation*}
\begin{split}
\int_{\mathbb{R}^4}&f_{k_1,j_1}(\xi_1,\tau_1)h'_+(\tau_2)\cdot\eta_0(\mu_2)\frac{1}{\mu_2+i/2^{k_2+1}}\\
&\times m(\xi_1-\sqrt{\tau_2}+\mu_2, \tau_1+\tau_2)\,d\mu_2d\tau_2d\xi_1d\tau_1\leq C||m||_{L^2}||h'_+||_{L^2}2^{j_1/2}||f_{k_1,j_1}||_{L^2}.
\end{split}
\end{equation*}
Let $\widetilde{m}(\xi,\tau)=\int_{\mathbb{R}}m(\xi+\mu_2,\tau)\eta_0(\mu_2)(\mu_2+i/2^{k_2+1})^{-1}\,d\mu_2$, $||\widetilde{m}||_{L^2}\leq C||m||_{L^2}$. In the left-hand side of the expression above we make the change of variable $\tau_1=\mu_1+\omega(\xi_1)$, $f_{k_1,j_1}^\#(\xi_1,\mu_1)=f_{k_1,j_1}(\xi_1,\mu_1+\omega(\xi_1))$. It suffices to prove that
\begin{equation}\label{vt30}
\begin{split}
\int_{\mathbb{R}^3}f_{k_1,j_1}^\#(\xi_1,\mu_1)&h'_+(\tau_2)\cdot \widetilde{m}(\xi_1-\sqrt{\tau_2}, \mu_1+\omega(\xi_1)+\tau_2)\,d\tau_2d\xi_1d\mu_1\\
&\leq C||\widetilde{m}||_{L^2}||h'_+||_{L^2}\cdot2^{j_1/2}||f_{k_1,j_1}^\#||_{L^2}.
\end{split}
\end{equation}
The integral in the left-hand side of \eqref{vt30} is over the set $$(\xi_1,\mu_1,\tau_2)\in\widetilde{I}_{k_1}\times\widetilde{I}_{j_1}\times[2^{2k_2-4},2^{2k_2+4}].$$
Using the Cauchy--Schwartz inequality, for \eqref{vt30} it suffices to prove that
\begin{equation*}
\sup_{\mu_1\in\mathbb{R}}\int_{\widetilde{I}_{k_1}\times[2^{2k_2-4},2^{2k_2+4}]}|\widetilde{m}(\xi_1-\sqrt{\tau_2}, \mu_1+\omega(\xi_1)+\tau_2)|^2\,d\tau_2d\xi_1\leq C||\widetilde{m}||_{L^2}^2,
\end{equation*}
which is easy to see by changing variables and recalling that $k_1\leq k_2-8$. This completes the proof of \eqref{vt13}.
\end{proof}

\section{Bilinear estimates II}\label{bilinear2}

In this section we prove two bilinear estimates, which correspond to $\mathrm{High}\times\mathrm{High}\to\mathrm{Low}$ interactions.

\newtheorem{Lemmae1}{Proposition}[section]
\begin{Lemmae1}\label{Lemmae1}
Assume $k,k_1,k_2\in\mathbb{Z}_+$ have the property that
$\mathrm{max}\,(k,k_1,k_2)\leq\mathrm{min}\,(k,k_1,k_2)+30$,
$f_{k_1}\in Z_{k_1}$, and $f_{k_2}\in Z_{k_2}$. Then
\begin{equation}\label{bk1}
2^k\big|\big|\eta_k(\xi)\cdot A_k(\xi,\tau)^{-1}f_{k_1}\ast
f_{k_2}\big|\big|_{Z_k}\leq
C||f_{k_1}||_{Z_{k_1}}||f_{k_2}||_{Z_{k_2}}.
\end{equation}
Moreover, any spaces $Z_0$ in the right-hand side of \eqref{bk1}
can be replaced with $\overline{Z}_0$.
\end{Lemmae1}

\newtheorem{Lemmae2}[Lemmae1]{Proposition}
\begin{Lemmae2}\label{Lemmae2}
Assume $k,k_1,k_2\in\mathbb{Z}_+$, $k_1,k_2\geq k+10$, $|k_1-k_2|\leq 2$, $f_{k_1}\in Z_{k_1}$, and $f_{k_2}\in Z_{k_2}$. Then
\begin{equation}\label{bk2}
\big|\big|\xi\cdot\eta_k(\xi)\cdot A_k(\xi,\tau)^{-1}f_{k_1}\ast f_{k_2}\big|\big|_{X_k}\leq C2^{-k/4}||f_{k_1}||_{Z_{k_1}}||f_{k_2}||_{Z_{k_2}}.
\end{equation}
\end{Lemmae2}

The main ingredients in the proofs of Propositions \ref{Lemmae1} and \ref{Lemmae2} are the definitions, the representations \eqref{repr1} and \eqref{repr2}, Lemma \ref{Lemmaa1}, and Corollary \ref{Lemmad2}.

\begin{proof}[Proof of Proposition \ref{Lemmae1}] We analyze two cases.

{\bf{Case 1: $\mathrm{min}\,(k,k_1,k_2)\geq 200$}}. In this case we prove the (stronger) bound \eqref{bk1} with the space $Z_k$ replaced by $X_k$ in the left-hand side. We show first that if $j_1,j_2\geq 0$, $f_{k_1,j_1}$ is an $L^2$ function supported in $D_{k_1,j_1}$, and $f_{k_2,j_2}$ is an $L^2$ function supported in $D_{k_2,j_2}$, then
\begin{equation}\label{bo40}
\begin{split}
2^{k}\sum_{j}2^{j/2}&\beta_{k,j}||\eta_k(\xi)\eta_j(\tau-\omega(\xi))(\tau-\omega(\xi)+i)^{-1}(f_{k_1,j_1}\ast f_{k_2,j_2})||_{L^2}\\
&\leq C\gamma(j_1,j_2,k)2^{j_1/2}\beta_{k_1,j_1}||f_{k_1,j_1}||_{L^2}\cdot 2^{j_2/2}\beta_{k_2,j_2}||f_{k_2,j_2}||_{L^2},
\end{split}
\end{equation}
where
\begin{equation}\label{bk5}
\gamma(j_1,j_2,k)=
\begin{cases}
2^{-\max(j_1,j_2)/4}&\text{ if }\max(j_1,j_2)\leq 2k-80;\\
2^{-\min(j_1,j_2)/8}&\text{ if }\max(j_1,j_2)\geq 2k-80.
\end{cases}
\end{equation}
To prove \eqref{bo40}, we notice that, in view of \eqref{nev2},
$\eta_k(\xi)\eta_j(\tau-\omega(\xi))(\tau-\omega(\xi)+i)^{-1}(f_{k_1,j_1}\ast
f_{k_2,j_2})\equiv 0$ unless
\begin{equation}\label{bo41}
\begin{cases}
&\mathrm{max}\,(j,j_1,j_2)\in[2k-70,2k+70]\text{ or}\\
&\mathrm{max}\,(j,j_1,j_2)\geq 2k+70\text{ and }
\mathrm{max}\,(j,j_1,j_2)-\mathrm{med}\,(j,j_1,j_2)\leq 10.
\end{cases}
\end{equation}
We notice that for $j,j_1,j_2$ as in \eqref{bo41},
$\beta_{k,j}\leq C\beta_{k_1,j_1}\beta_{k_2,j_2}$. Also, using
\eqref{on3},
\begin{equation*}
\begin{split}
||\eta_k(\xi)&\eta_j(\tau-\omega(\xi))(\tau-\omega(\xi)+i)^{-1}(f_{k_1,j_1}\ast f_{k_2,j_2})||_{L^2}\\
&\leq C2^{-j}2^{(j+j_1+j_2)/2}2^{-\mathrm{max}\,(j,j_1,j_2)/2}
2^{-\mathrm{med}\,(j,j_1,j_2)/4}||f_{k_1,j_1}||_{L^2}||f_{k_2,j_2}||_{L^2}.
\end{split}
\end{equation*}
Thus, for \eqref{bo40}, it suffices to prove that
\begin{equation}\label{bo42}
2^{k}\sum_{j}2^{-\mathrm{max}\,(j,j_1,j_2)/2}2^{-\mathrm{med}\,(j,j_1,j_2)/4}\leq
C\gamma(j_1,j_2,k),
\end{equation}
where the sum in \eqref{bo42} is taken over $j$ satisfying
\eqref{bo41}. If $\max(j_1,j_2)\leq 2k-80$ then $j\in
[2k-70,2k+70]$ and the bound \eqref{bo42} follows easily from the
definition \eqref{bk5}. If $j_1=\max(j_1,j_2)\geq 2k-80$ then the
sum in \eqref{bo42} is taken over $j\leq j_1+C$ and is dominated
by
\begin{equation*}
C2^k\sum_{j\leq j_1+C}2^{-j_1/2}2^{-\max(j,j_2)/4}\leq C(j_2+1)2^{-j_2/4},
\end{equation*}
which suffices. The case $j_2=\max(j_1,j_2)\geq 2k-80$ is
identical. This completes the proof of \eqref{bo40}.

We turn to the proof of \eqref{bk1}. We use the representation
\eqref{repr1}. If $f_{k_1}=f_{k_1,j_1}\in X_{k_1}$ and
$f_{k_2}=f_{k_2,j_2}\in X_{k_2}$ then \eqref{bk1} follows directly
from \eqref{bo40} and the definitions. Assume now that
$f_{k_1}=g_{k_1}\in Y_{k_1}$, $f_{k_2}=g_{k_2}\in Y_{k_2}$,
$||f_{k_1}||_{Z_{k_1}}\approx||g_{k_1}||_{Y_{k_1}}$, and
$||f_{k_2}||_{Z_{k_2}}\approx ||g_{k_2}||_{Y_{k_2}}$. For
$j_1\in[0,k_1]$ and $j_2\in[0,k_2]$ let
$g_{k_1,j_1}(\xi,\tau)=\eta_{j_1}(\tau-\omega(\xi))g_{k_1}(\xi,\tau)$
and
$g_{k_2,j_2}(\xi,\tau)=\eta_{j_2}(\tau-\omega(\xi))g_{k_2}(\xi,\tau)$.
We use \eqref{bo40}, Lemma \eqref{Lemmaa1} (b), and the definition
\eqref{bk5} in the case $\max(j_1,j_2)\leq 2k-80$ to write
\begin{equation*}
\begin{split}
2^k\big|\big|&\eta_k(\xi)\cdot(\tau-\omega(\xi)+i)^{-1}g_{k_1}\ast g_{k_2}\big|\big|_{X_k}\\
&\leq C\sum_{j_1,j_2\leq k+30}2^k\big|\big|\eta_k(\xi)\cdot(\tau-\omega(\xi)+i)^{-1}g_{k_1,j_1}\ast g_{k_2,j_2}\big|\big|_{X_k}\\
&\leq C\sum_{j_1,j_2\leq k+30}\gamma(j_1,j_2,k)2^{j_1/2}||g_{k_1,j_1}||_{L^2}2^{j_2/2}||g_{k_2,j_2}||_{L^2}\\
&\leq C||g_{k_1}||_{Y_{k_1}}||g_{k_2}||_{Y_{k_2}},
\end{split}
\end{equation*}
as desired.

Finally, assume $f_{k_1}=f_{k_1,j_1}\in X_{k_1}$,
$f_{k_2}=g_{k_2}\in Y_{k_2}$,
$||f_{k_2}||_{Z_{k_2}}\approx||g_{k_2}||_{Y_{k_2}}$, and
$||f_{k_1}||_{Z_{k_1}}\approx2^{j_1/2}\beta_{k_1,j_1}||f_{k_1,j_1}||_{L^2}$,
and write $g_{k_2}=\sum_{j_2=0}^{k_2}g_{k_2,j_2}$ as  before. If
$j_1\leq 2k-80$ then we can use the same computation as before. If
$j_1\geq 2k-80$ then we use \eqref{bo40}, Lemma \eqref{Lemmaa1}
(b), and the definition \eqref{bk5} to write
\begin{equation*}
\begin{split}
2^k\big|\big|&\eta_k(\xi)\cdot(\tau-\omega(\xi)+i)^{-1}f_{k_1,j_1}\ast g_{k_2}\big|\big|_{X_k}\\
&\leq C\sum_{j_2\leq k_2}2^k\big|\big|\eta_k(\xi)\cdot(\tau-\omega(\xi)+i)^{-1}f_{k_1,j_1}\ast g_{k_2,j_2}\big|\big|_{X_k}\\
&\leq C\sum_{j_2\leq k_2}2^{-j_2/8}2^{j_1/2}\beta_{k_1,j_1}||f_{k_1,j_1}||_{L^2}2^{j_2/2}||g_{k_2,j_2}||_{L^2}\\
&\leq C2^{j_1/2}\beta_{k_1,j_1}||f_{k_1,j_1}||_{L^2}||g_{k_2}||_{Y_{k_2}},
\end{split}
\end{equation*}
as desired. This completes the proof of \eqref{bk1} in the case $\mathrm{min}\,(k,k_1,k_2)\geq 200$.
\medskip

{\bf{Case 2: $\mathrm{min}\,(k,k_1,k_2)\leq 200$}}. In view of the
hypothesis, $\mathrm{max}\,(k,k_1,k_2)\leq 230$. If $k_1=0$ or
$k_2=0$ we may replace the spaces $Z_0$ in the right-hand of
\eqref{bk1} with the larger spaces $\overline{Z}_0$, see the
definition \eqref{def1'''}. Clearly, the proofs are identical to
the proofs in the corresponding cases $k_1=1$ or $k_2=1$.
Therefore we may assume $k_1,k_2\geq 1$. In view of  Lemma
\ref{Lemmaa1} (b) and the representation \eqref{repr1}, we may
assume $f_{k_1}=f_{k_1,j_1}$ is supported in $D_{k_1,j_1}$,
$f_{k_2}=f_{k_2,j_2}$ is supported in $D_{k_2,j_2}$,
$||f_{k_1}||_{Z_{k_1}}\approx2^{j_1/2}\beta_{k_1,j_1}||f_{k_1,j_1}||_{L^2}\approx
2^{j_1}||f_{k_1,j_1}||_{L^2}$, and
$||f_{k_2}||_{Z_{k_2}}\approx2^{j_2/2}\beta_{k_2,j_2}||f_{k_2,j_2}||_{L^2}\approx
2^{j_2}||f_{k_2,j_2}||_{L^2}$. Using the definitions and the fact
that $k\leq 230$, for \eqref{bk1} it suffices to prove that
\begin{equation}\label{bm1}
\begin{split}
\sum_{j}2^j\big|\big|\mathcal{F}^{-1}
[\eta_j(\tau)\eta_k(\xi)(\tau+i)^{-1}&f_{k_1,j_1}\ast f_{k_2,j_2}]\big|\big|_{L^1_xL^2_t}\\
&\leq C2^{j_1}||f_{k_1,j_1}||_{L^2}\cdot 2^{j_2}||f_{k_2,j_2}||_{L^2}.
\end{split}
\end{equation}
By examining the supports of the functions, we may assume that the
sum in \eqref{bm1} is taken over
\begin{equation}\label{bm2}
j\leq\max(j_1,j_2)+C.
\end{equation}
Assume $j_1=\max(j_1,j_2)$ (the case $j_2=\max(j_1,j_2)$ is
identical). The left-hand side of \eqref{bm1} is dominated by
\begin{equation*}
\begin{split}
C\sum_{j\leq j_1+C}\big|\big|\mathcal{F}^{-1}(f_{k_1,j_1}\ast f_{k_2,j_2})\big|\big|_{L^1_xL^2_t}&
\leq (j_1+C)||\mathcal{F}^{-1}(f_{k_1,j_1})||_{L^2}||\mathcal{F}^{-1}(f_{k_2,j_2})||_{L^2_xL^\infty_t}\\
&\leq C2^{j_1}||f_{k_1,j_1}||_{L^2}\cdot 2^{j_2/2}||f_{k_2,j_2}||_{L^2},
\end{split}
\end{equation*}
which completes the proof of \eqref{bm1}.
\end{proof}
\medskip

For later use we rewrite the stronger bound we proved in this last case: if $k,k_1,k_2\in\mathbb{Z}_+$ have the property that
$\mathrm{max}\,(k,k_1,k_2)\leq\mathrm{min}\,(k,k_1,k_2)+30\leq 230$,
$f_{k_1}\in \overline{Z}_{k_1}$, and $f_{k_2}\in \overline{Z}_{k_2}$, then
\begin{equation}\label{hh90}
2^k\big|\big|\eta_k(\xi)\cdot A_k(\xi,\tau)^{-1}f_{k_1}\ast
f_{k_2}\big|\big|_{Z_k}\leq
C||f_{k_1}||_{\overline{Z}_{k_1}}||f_{k_2}||_{\overline{Z}_{k_2}}.,
\end{equation}
where $\overline{Z}_k=Z_k$ if $k\geq 1$ and $\overline{Z}_k=\overline{Z}_0$ if $k=0$.

\begin{proof}[Proof of Proposition \ref{Lemmae2}]  We analyze two cases.

{\bf{Case 1: $k\geq 1$}}. We show first that if $j_1,j_2\geq 0$,
$f_{k_1,j_1}$ is an $L^2$ function supported in $D_{k_1,j_1}$, and
$f_{k_2,j_2}$ is an $L^2$ function supported in $D_{k_2,j_2}$,
then
\begin{equation}\label{bl40}
\begin{split}
2^{k}\sum_{j}2^{j/2}&\beta_{k,j}||\eta_k(\xi)\eta_j(\tau-\omega(\xi))(\tau-\omega(\xi)+i)^{-1}(f_{k_1,j_1}\ast f_{k_2,j_2})||_{L^2}\\
&\leq C\gamma'(j_1,j_2,k)2^{j_1/2}\beta_{k_1,j_1}||f_{k_1,j_1}||_{L^2}\cdot 2^{j_2/2}\beta_{k_2,j_2}||f_{k_2,j_2}||_{L^2},
\end{split}
\end{equation}
where
\begin{equation}\label{bl5}
\gamma'(j_1,j_2,k)=(2^{k/2}+2^{\max(j_1,j_2)/4})^{-1}.
\end{equation}
To prove \eqref{bl40}, we notice that, in view of \eqref{nev2},
$\eta_k(\xi)\eta_j(\tau-\omega(\xi))(\tau-\omega(\xi)+i)^{-1}(f_{k_1,j_1}\ast
f_{k_2,j_2})\equiv 0$ unless
\begin{equation}\label{bl41}
\begin{cases}
&\mathrm{max}\,(j,j_1,j_2)\in[k+k_1-10,k+k_1+10]\text{ or}\\
&\mathrm{max}\,(j,j_1,j_2)\geq k+k_1+10\text{ and }\mathrm{max}\,(j,j_1,j_2)-\mathrm{med}\,(j,j_1,j_2)\leq 10.
\end{cases}
\end{equation}
Also, combining \eqref{on32} and \eqref{on3},
\begin{equation*}
\begin{split}
||&\eta_k(\xi)\eta_j(\tau-\omega(\xi))(\tau-\omega(\xi)+i)^{-1}(f_{k_1,j_1}\ast f_{k_2,j_2})||_{L^2}\leq C2^{-j}2^{(j+j_1+j_2)/2}\\
&\times[2^{(j+k)/2}+2^{(\max(j_1,j_2)+k_1)/2}+2^{\mathrm{max}\,(j,j_1,j_2)/2}
2^{\mathrm{med}\,(j,j_1,j_2)/4}]^{-1}||f_{k_1,j_1}||_{L^2}||f_{k_2,j_2}||_{L^2}.
\end{split}
\end{equation*}
Thus, for \eqref{bl40}, it suffices to prove that
\begin{equation}\label{bl42}
\begin{split}
2^{k}\sum_{j}\beta_{k,j}[2^{(j+k)/2}+
2^{(\max(j_1,j_2)+k_1)/2}+&2^{\mathrm{max}\,(j,j_1,j_2)/2}
2^{\mathrm{med}\,(j,j_1,j_2)/4}]^{-1}\\
&\leq C\gamma'(j_1,j_2,k)\beta_{k_1,j_1}\beta_{k_2,j_2},
\end{split}
\end{equation}
where the sum in \eqref{bl42} is taken over $j$ satisfying
\eqref{bl41}. If $\max(j_1,j_2)\leq k+k_1-20$ then $j\in
[k+k_1-10,k+k_1+10]$; we  ignore the term
$2^{(\max(j_1,j_2)+k_1)/2}$ and the bound \eqref{bl42} follows
easily from the definitions. If $j_1=\max(j_1,j_2)\geq k+k_1-20$
then the sum in \eqref{bl42} is taken over $j\leq j_1+C$. The
left-hand side of \eqref{bl42} is dominated by
\begin{equation*}
C2^k\sum_{j\leq j_1+C}\beta_{k,j}2^{-(\max(j_1,j_2)+k_1)/2}\leq C2^{-k_1/2}
\leq C\gamma'(j_1,j_2,k)\beta_{k_1,j_1}.
\end{equation*}
The case $j_2=\max(j_1,j_2)\geq  k+k_1-20$ is identical, which
completes the proof of \eqref{bl40}.

We turn to the proof of \eqref{bk2}. We use the representation
\eqref{repr1}. If $f_{k_1}=f_{k_1,j_1}\in X_{k_1}$ and
$f_{k_2}=f_{k_2,j_2}\in X_{k_2}$ then \eqref{bk2} follows directly
from \eqref{bl40} and the definitions. Assume now that
$f_{k_1}=g_{k_1}\in Y_{k_1}$, $f_{k_2}=g_{k_2}\in Y_{k_2}$,
$||f_{k_1}||_{Z_{k_1}}\approx||g_{k_1}||_{Y_{k_1}}$, and
$||f_{k_2}||_{Z_{k_2}}\approx||g_{k_2}||_{Y_{k_2}}$. For
$j_1\in[0,k_1]$ and $j_2\in[0,k_2]$ let
$g_{k_1,j_1}(\xi,\tau)=\eta_{j_1}(\tau-\omega(\xi))g_{k_1}(\xi,\tau)$
and
$g_{k_2,j_2}(\xi,\tau)=\eta_{j_2}(\tau-\omega(\xi))g_{k_2}(\xi,\tau)$.
We use \eqref{bl40} and Lemma \eqref{Lemmaa1} (b) to write
\begin{equation*}
\begin{split}
2^k\big|\big|&\eta_k(\xi)\cdot(\tau-\omega(\xi)+i)^{-1}g_{k_1}\ast g_{k_2}\big|\big|_{X_k}\\
&\leq C\sum_{j_1,j_2\leq k_1+10}2^k\big|\big|\eta_k(\xi)\cdot(\tau-\omega(\xi)+i)^{-1}g_{k_1,j_1}\ast g_{k_2,j_2}\big|\big|_{X_k}\\
&\leq C\sum_{j_1,j_2\leq k_1+10}\gamma'(j_1,j_2,k)2^{j_1/2}||g_{k_1,j_1}||_{L^2}2^{j_2/2}||g_{k_2,j_2}||_{L^2}\\
&\leq C2^{-k/4}||g_{k_1}||_{Y_{k_1}}||g_{k_2}||_{Y_{k_2}},
\end{split}
\end{equation*}
as desired. Finally, if $f_{k_1}=f_{k_1,j_1}\in X_{k_1}$,
$f_{k_2}=g_{k_2}\in Y_{k_2}$,
$||f_{k_2}||_{Z_{k_2}}\approx||g_{k_2}||_{Y_{k_2}}$, and
$||f_{k_1}||_{Z_{k_1}}\approx2^{j_1/2}\beta_{k_1,j_1}||f_{k_1,j_1}||_{L^2}$,
we write $g_{k_2}=\sum_{j_2=0}^{k_2}g_{k_2,j_2}$ as  before and
repeat the same estimate, without the sum in $j_1$. This completes
the proof of \eqref{bk2} in the case $k\geq 1$.

{\bf{Case 2: $k=0$}}. We show first that if $j_1,j_2\geq 0$,
$f_{k_1,j_1}$ is an $L^2$ function supported in $D_{k_1,j_1}$, and
$f_{k_2,j_2}$ is an $L^2$ function supported in $D_{k_2,j_2}$,
then
\begin{equation}\label{bc40}
\begin{split}
\sum_{k'=-\infty}^1&\sum_{j=0}^\infty2^{j-k'}||\chi_{k'}(\xi)\eta_j(\tau)
\cdot\xi(\tau+i)^{-1}(f_{k_1,j_1}\ast f_{k_2,j_2})||_{L^2}\\
&\leq C2^{-\max(j_1,j_2)/4}\cdot 2^{j_1/2}\beta_{k_1,j_1}||f_{k_1,j_1}
||_{L^2}\cdot 2^{j_2/2}\beta_{k_2,j_2}||f_{k_2,j_2}||_{L^2}.
\end{split}
\end{equation}
To prove \eqref{bc40}, we notice that, in view of \eqref{nev2},
$\chi_{k'}(\xi)\eta_j(\tau)\cdot\xi(\tau+i)^{-1}(f_{k_1,j_1}\ast
f_{k_2,j_2})\equiv 0$ unless
\begin{equation}\label{bc41}
\begin{cases}
&\mathrm{max}\,(j,j_1,j_2)\in[k'+k_1-10,k'+k_1+10]\text{ or}\\
&\mathrm{max}\,(j,j_1,j_2)\geq k'+k_1+10\text{ and }
\mathrm{max}\,(j,j_1,j_2)-\mathrm{med}\,(j,j_1,j_2)\leq 10.
\end{cases}
\end{equation}
Also, using \eqref{on31},
\begin{equation*}
\begin{split}
||\chi_{k'}(\xi)&\eta_j(\tau)\cdot\xi(\tau+i)^{-1}
(f_{k_1,j_1}\ast f_{k_2,j_2})||_{L^2}\\
&\leq C2^{k'-j}2^{k'/2}2^{(j_1+j_2)/2}2^{-\max(j_1,j_2)/2}
||f_{k_1,j_1}||_{L^2}||f_{k_2,j_2}||_{L^2}.
\end{split}
\end{equation*}
Thus, for \eqref{bc40}, it suffices to prove that
\begin{equation}\label{bc42}
\begin{split}
\sum_{k'=-\infty}^1\sum_{j}2^{k'/2}\leq C2^{\max(j_1,j_2)/4},
\end{split}
\end{equation}
where the sum in \eqref{bc42} is taken over $j$ satisfying
\eqref{bc41}. If $\max(j_1,j_2)\leq k'+k_1-20$ then $j\in
[k'+k_1-10,k'+k_1+10]$, so \eqref{bc42} is clear. If
$\max(j_1,j_2)\geq k'+k_1-20$, then the sum in \eqref{bc42} is
taken over $j\leq\max(j_1,j_2)+C$, and  \eqref{bc42} follows
easily.

Given \eqref{bc40}, the bound \eqref{bk2} follows as  in the case
$k\geq 1$, using the definition of the space $X_0$. This completes
the proof of  Proposition \ref{Lemmae2}.
\end{proof}

For later use, we notice
that the bound \eqref{bc40} also shows that
\begin{equation}\label{bc50}
\big|\big|\eta_0(\xi)\cdot (\tau+i)^{-1}f_{k_1}\ast f_{k_2}
\big|\big|_{\overline{Z}_0}\leq
C||f_{k_1}||_{Z_{k_1}}||f_{k_2}||_{Z_{k_2}}.
\end{equation}

\section{Multiplication by smooth bounded functions}\label{mult}

In this section we consider operators on $Z_k$ given by
convolutions with Fourier transforms of certain smooth bounded
functions. For integers $N\geq 100$ we define the space of {\it{admissible factors}}
\begin{equation}\label{ar3}
\begin{split}
&S^\infty_N=\{m:\mathbb{R}^2\to\mathbb{C}: m\text{ is supported in }\mathbb{R}\times[-2,2]\text{ and }\\
&||m||_{S^\infty_N}:=\sum_{\sigma_1=0}^N||\partial_t^{\sigma_1}m||_{L^\infty_{x,t}}+\sum_{\sigma_1=0}^N\sum_{\sigma_2=1}^N||\partial_t^{\sigma_1}
\partial_x^{\sigma_2}m||_{L^2_{x,t}}<\infty\}.
\end{split}
\end{equation}
The precise value of $N$ is not important (in fact, we will always take $N=100$ or $N=110$). Notice that bounded functions such as $\psi(t)e^{iqU_0}$, $q\in\mathbb{R}$, $U_0$ as in \eqref{fg9}, are in $S^\infty_N$. We also define the space of {\it{restricted admissible factors}}
\begin{equation}\label{ar3'}
\begin{split}
S^2_N=\{m:&\mathbb{R}^2\to\mathbb{C}: m\text{ is supported in }\mathbb{R}\times[-2,2]\text{ and }\\
&||m||_{S^2_N}:=\sum_{\sigma_1=0}^N\sum_{\sigma_2=0}^N||\partial_t^{\sigma_1}
\partial_x^{\sigma_2}m||_{L^2_{x,t}}<\infty\}.
\end{split}
\end{equation}
Using Sobolev imbedding theorem, it is easy to verify the following properties:
\begin{equation}\label{ar3''}
\begin{cases}
&S^2_N\subseteq S^\infty_{N-10};\\
&S^\infty_N\cdot S^\infty_N\subseteq S^\infty_{N-10};\\
&S^2_N\cdot S^\infty_N\subseteq S^2_{N-10};\\
&\partial_x S^\infty_N\subseteq S^2_{N-10}.
\end{cases}
\end{equation}

For $k\in\mathbb{Z}_+$ we define
\begin{equation}\label{ar1}
Z_k^{\mathrm{high}}=\{f_k\in Z_k:f_k\text{ is supported in }
\{\tau-\omega(\xi)\in\bigcup_{j\geq k-20}\widetilde{I}_j\}\}.
\end{equation}
Clearly, $Z_k^{\mathrm{high}}=Z_k$ if $k\leq 20$. 
For $k\in\mathbb{Z}_+$ and $\epsilon\in\{-1,0\}$ let
$A_k^\epsilon(\xi,\tau)=[A_k(\xi,\tau)]^\epsilon$.

\newtheorem{Lemmab4}{Lemma}[section]
\begin{Lemmab4}\label{Lemmab4}
Assume $k_1,k_2\in\mathbb{Z}_+$, $|k_1-k_2|\leq 10$, and
$f_{k_1}^{\mathrm{high}}\in Z_{k_1}^{\mathrm{high}}$. Then, for
$m\in S^\infty_{100}$ and $\epsilon\in\{-1,0\}$,
\begin{equation}\label{ar4}
\Big|\Big|\eta_{k_2}(\xi_2)A_{k_2}^\epsilon(\xi_2,\tau_2)\cdot
\mathcal{F}[m\cdot\mathcal{F}^{-1}(f_{k_1}^{\mathrm{high}})](\xi_2,\tau_2)
\Big|\Big|_{Z_{k_2}}\leq
C ||m||_{S^\infty_{100}}\cdot||A_{k_1}^{\epsilon} f_{k_1}^{\mathrm{high}}||_{Z_{k_1}}.
\end{equation}
\end{Lemmab4}

{\bf{Remark:}} It is easy to see that a sharp bound like \eqref{ar4} cannot hold for functions $f_{k_1}$ of low modulation. Fortunately, we do not need to consider convolutions of low-modulation functions and Fourier transforms of admissible factors, in view of the identity \eqref{mv10}. 

\begin{proof}[Proof of Lemma \ref{Lemmab4}] We may assume $||m||_{S^\infty_{100}}=1$. For any $j''\in\mathbb{Z}_+$
and $k''\in\mathbb{Z}$ let
\begin{equation}\label{ar12}
m_{k'',j''}=\mathcal{F}^{-1}\big[\eta_{j''}(\tau)\chi_{k''}(\xi)\mathcal{F}(m)\big],
\end{equation}
and $m_{\leq k'',j''}=\sum_{k'''\leq k''}m_{k''',j''}$. Using
\eqref{ar3} and the Sobolev imbedding theorem
\begin{equation*}
\big|\big|\partial_t^{\sigma_1}\partial_x^{\sigma_2}m\big|\big|_{L^2_xL^\infty_t}
\leq C\text { for any
}\sigma_1\in\mathbb{Z}\cap[0,90],\,\sigma_2\in\mathbb{Z}\cap[1,90].
\end{equation*}
Thus, for any $j''\in\mathbb{Z}_+$ and $k''\in\mathbb{Z}$,
\begin{equation}\label{ar10}
\begin{cases}
&||m_{\leq k'',j''}||_{L^\infty_{x,t}}\leq C2^{-80j''};\\
&2^{k''}||m_{k'',j''}||_{L^2_xL^\infty_t}+||m_{k'',j''}||_{L^\infty_{x,t}}\leq
C(1+2^{-k''})^{-80}2^{-80j''}.
\end{cases}
\end{equation}

We turn now to the proof of \eqref{ar4}. Assume first that
$k_1,k_2\geq 1$. In view of the definition of
$Z_k^{\mathrm{high}}$ and  Lemma \ref{Lemmaa1} (b), we may assume that $f_{k_1}^{\mathrm{high}}=f_{k_1,j_1}$ is an $L^2$ function
supported in $D_{k_1,j_1}$, $j_1\geq k_1-20$, $||A_{k_1}^\epsilon
f_{k_1}^\mathrm{high}||_{Z_{k_1}}\approx 2^{\epsilon
j_1}2^{j_1/2}\beta_{k_1,j_1}||f_{k_1,j_1}||_{L^2}$. We write
\begin{equation}\label{ar20}
m=\sum_{j''=0}^\infty m_{\leq -100,j''}+
\sum_{k''=-99}^\infty\sum_{j''=0}^\infty m_{k'',j''}.
\end{equation}
For \eqref{ar4} it suffices to prove that for $\epsilon\in\{-1,0\}$
\begin{equation}\label{ar21}
\begin{split}
&\sum_{j''\geq 0}\big|\big|\eta_{k_2}(\xi_2)A^\epsilon(\xi_2,\tau_2)
\cdot[f_{k_1,j_1}\ast\mathcal{F}(m_{\leq -100,j''})](\xi_2,\tau_2)
\big|\big|_{Z_{k_2}}\\
&+\sum_{k''\geq -99}\sum_{j''\geq 0}\big|\big|\eta_{k_2}(\xi_2)
A^\epsilon(\xi_2,\tau_2)\cdot
[f_{k_1,j_1}\ast\mathcal{F}(m_{k'',j''})](\xi_2,\tau_2)\big|\big|_{Z_{k_2}}\\
&\leq C2^{\epsilon
j_1}\cdot2^{j_1/2}\beta_{k_1,j_1}||f_{k_1,j_1}||_{L^2}.
\end{split}
\end{equation}

To bound the first sum in \eqref{ar21} we make the changes of
variables $\tau_2=\mu_2+\omega(\xi_2)$,
$\tau_1=\mu_1+\omega(\xi_1)$, and write
\begin{equation*}
\begin{split}
&f_{k_1,j_1}\ast\mathcal{F}(m_{\leq -100,j''})(\xi_2,\mu_2+\omega(\xi_2))\\
&=\int_{\mathbb{R}^2}f_{k_1,j_1}(\xi_1,\mu_1+\omega(\xi_1))\mathcal{F}(m_{\leq
-100,j''})(\xi_2-\xi_1,\mu_2-\mu_1+\omega(\xi_2)-\omega(\xi_1))\,d\xi_1
d\tau_1.
\end{split}
\end{equation*}
By examining the supports of the functions and using the fact that
$|\omega(\xi_2)-\omega(\xi_1)|\leq 2^{k_1-50}$ if
$|\xi_2-\xi_1|\leq 2^{-99}$, together with $j_1\geq k_1-20$, we
see that
$\eta_{j_2}(\tau_2-\omega(\xi_2))\cdot[f_{k_1,j_1}\ast\mathcal{F}(m_{\leq
-100,j''})](\xi_2,\tau_2)\equiv 0$ unless
\begin{equation}\label{ar25}
|j_1-j_2|\leq C\,\text{ or }j_1,j_2\leq j''+C.
\end{equation}
We use the $X_{k_2}$ norm to bound the first sum in \eqref{ar21}.
Using Plancherel theorem and \eqref{ar10},
\begin{equation*}
\big|\big|f_{k_1,j_1}\ast\mathcal{F}(m_{\leq -100,j''})
\big|\big|_{L^2_{\xi_2,\tau_2}}\leq
C2^{-80j''}||f_{k_1,j_1}||_{L^2}.
\end{equation*}
Thus, the $X_{k_2}$ norm of the first sum in \eqref{ar21} is dominated by
\begin{equation*}
C\sum_{j''\geq 0}\sum_{j_2\geq 0}2^{\epsilon j_2}2^{j_2/2}
\beta_{k_2,j_2}2^{-80j''}||f_{k_1,j_1}||_{L^2},
\end{equation*}
where the sum is over $j_2,j''$ satisfying \eqref{ar25}. The bound
\eqref{ar21} for the first sum follows easily (recall that
$|k_1-k_2|\leq 10$).

To bound the second sum in \eqref{ar21} assume first that
$\epsilon=0$. We notice that if
$|\xi_2-\xi_1|\in[2^{k''-1},2^{k''+1}]$ then
$|\omega(\xi_2)-\omega(\xi_1)|\leq C2^{k_1+k''}$, so
$\eta_{j_2}(\tau_2-\omega(\xi_2))\cdot[f_{k_1,j_1}\ast\mathcal{F}(m_{k'',j''})]
(\xi_2,\tau_2)\equiv 0$ unless
\begin{equation}\label{ar26}
|j_1-j_2|\leq 4\,\text{ or }\,j_1,j_2\leq k_1+k''+j''+C\text{ and }\\
\end{equation}
Using Plancherel theorem and \eqref{ar10},
\begin{equation}\label{ng2}
\big|\big|f_{k_1,j_1}\ast\mathcal{F}(m_{k'',j''})\big|\big|_{L^2_{\xi_2,\tau_2}}
\leq C2^{-80k''}2^{-80j''}||f_{k_1,j_1}||_{L^2}.
\end{equation}
The bound \eqref{ar21} for the second sum follows by using the
$X_{k_2}$ norm since
\begin{equation*}
\sum_{j_2\leq j_1+k''+j''+C}2^{j_2/2}\beta_{k_2,j_2}\leq
C2^{10k''}2^{10j''} \cdot 2^{j_1/2}\beta_{k_1,j_1}.
\end{equation*}

We bound now the second sum in \eqref{ar21} when $\epsilon=-1$.
The main difficulty is the presence of the indices $j_2\ll j_1$.
In fact, for indices $j_2\geq j_1-10$, the argument above applies
since the left-hand side is multiplied by $2^{-j_2}$ and the
right-hand side is multiplied by $2^{-j_1}$. In view of
\eqref{ar26}, it suffices to prove that
\begin{equation}\label{ng1}
\begin{split}
&\sum_{k''+j''\geq j_1-k_1-C}\big|\big|\eta_{k_2}(\xi_2)
\eta_{\leq k_2-1}(\tau_2-\omega(\xi_2))A_{k_2}^{-1}(\xi_2,\tau_2)
[f_{k_1,j_1}\ast\mathcal{F}(m_{k'',j''})](\xi_2,\tau_2)\big|\big|_{Y_{k_2}}\\
&+\sum_{k''+j''\geq j_1-k_1-C}\sum_{j_2\geq k_2}^{j_2\leq
j_1-10}2^{-j_2/2}\beta_{k_2,j_2}\big|\big|\eta_{k_2}(\xi_2)
\eta_{j_2}(\tau_2-\omega(\xi_2))f_{k_1,j_1}\ast\mathcal{F}(m_{k'',j''})
\big|\big|_{L^2}\\
&\leq C2^{-j_1/2}\beta_{k_1,j_1}||f_{k_1,j_1}||_{L^2}.
\end{split}
\end{equation}
Using Lemma \ref{Lemmaa1} (c) and \eqref{ar10}, the first sum in
the left-hand side of \eqref{ng1} is dominated by
\begin{equation*}
C\sum_{k''+j''\geq
j_1-k_1-C}2^{-k_2/2}||f_{k_1,j_1}||_{L^2}\cdot||m_{k'',j''}||_{L^2_xL^\infty_t}\leq
C2^{-k_2/2}||f_{k_1,j_1}||_{L^2}\cdot 2^{-70(j_1-k_1)},
\end{equation*}
which clearly suffices. Using \eqref{ng2}, the second sum in the
left-hand side of \eqref{ng1} is dominated by
\begin{equation*}
C2^{-70(j_1-k_1)}||f_{k_1,j_1}||_{L^2}\cdot
\sup_{j_2\in[k_2,j_1]}2^{-j_2/2}\beta_{k_2,j_2}\leq
C2^{-j_1/2}||f_{k_1,j_1}||_{L^2},
\end{equation*}
which completes the proof of \eqref{ng1}.

We prove now the bound \eqref{ar4} in the case $k_1=k_2=0$. We use
the representation \eqref{repr2}. Assume first that
$f_0^{\mathrm{high}}=g_{0,j_1}$ is an $L^2$ function supported in
$\widetilde{I}_0\times\widetilde{I}_{j_1}$,
$||A^\epsilon_0f_0^{\mathrm{high}}||_{Z_0}\approx 2^{\epsilon
j_1}2^{j_1}||\mathcal{F}^{-1}(g_{0,j_1})||_{L^1_xL^2_t}$. We write
\begin{equation}\label{ar38}
m=\sum_{j''=0}^\infty m_{\leq 4,j''}+\sum_{k''=5}^\infty
\sum_{j''=0}^\infty m_{k'',j''}.
\end{equation}
and notice that $\eta_{0}(\xi_2)(g_{0,j_1}\ast
\mathcal{F}(m_{k'',j''}))(\xi_2,\tau_2)\equiv 0$ if $k''\geq 5$.
For \eqref{ar4}, using only the $Y_0$ norm, it suffices to prove
that for $\epsilon\in\{-1,0\}$
 \begin{equation}\label{ar40}
\begin{split}
\sum_{j''=0}^\infty\sum_{j_2=0}^\infty 2^{\epsilon j_2}&2^{j_2}
\big|\big|\mathcal{F}^{-1}[\eta_{j_2}(\tau_2)(g_{0,j_1}
\ast\mathcal{F}(m_{\leq 4,j''}))(\xi_2,\tau_2)]\big|\big|_{L^1_xL^2_t}\\
&\leq C2^{\epsilon j_1}2^{j_1}||\mathcal{F}^{-1}(g_{0,j_1})||_{L^1_xL^2_t}.
\end{split}
\end{equation}
By examining the supports of the functions,
$\eta_{j_2}(\tau_2)(g_{0,j_1}\ast\mathcal{F}(m_{\leq
4,j''}))(\xi_2,\tau_2)\equiv 0$ unless
\begin{equation}\label{ar33}
|j_2-j_1|\leq C\,\text{ or }\,j_1,j_2\leq j''+C.
\end{equation}
In addition,
\begin{equation*}
\begin{split}
\big|\big|\mathcal{F}^{-1}[\eta_{j_2}(\tau_2)(g_{0,j_1}\ast
\mathcal{F}(m_{\leq
4,j''}))(\xi_2,\tau_2)]\big|\big|_{L^1_xL^2_t}&\leq
C||\mathcal{F}^{-1}(g_{0,j_1})||_{L^1_xL^2_t}||m_{\leq
4,j''}||_{L^\infty_{x,t}}.
\end{split}
\end{equation*}
The bound \eqref{ar40} follows  from \eqref{ar10} and \eqref{ar33}.

Assume now that $f_0^{\mathrm{high}}=f_{0,j_1}^{k'}$ is a smooth
function supported in $D_{k',j_1}$, $k'\leq 1$, $||A_0^\epsilon
f_0^{\mathrm{high}}||_{Z_0}\approx 2^{\epsilon
j_1}2^{j_1-k'}||f_{0,j_1}^{k'}||_{L^2}$. We decompose
\begin{equation}\label{ar50}
m=\sum_{j''=0}^\infty m_{\leq k'-10,j''}+\sum_{k''=k'-9}^\infty
\sum_{j''=0}^\infty m_{k'',j''}.
\end{equation}
We observe that $f_{0,j_1}^{k'}\ast\mathcal{F}(m_{\leq
k'-10,j''})$ is supported in the set
$\{(\xi_2,\tau_2):|\xi_2|\in[2^{k'-2},2^{k'+2}]\}$. In addition,
$\eta_{j_2}(\tau_2)(f_{0,j_1}^{k'}\ast\mathcal{F}(m_{\leq
k'-3,j''}))(\xi_2,\tau_2)\equiv 0$ unless \eqref{ar33} holds. The
same argument as before, using Plancherel theorem and the bound
\eqref{ar10}, shows that
\begin{equation*}
\big|\big|\eta_{0}(\xi_2)A_0^\epsilon(\xi_2,\tau_2)\big[f_{0,j_1}^{k'}
\ast\sum_{j''=0}^\infty \mathcal{F}(m_{\leq
k'-10,j''})\big](\xi_2,\tau_2)\big|\big|_{X_{0}} \leq C2^{\epsilon
j_1}2^{j_1-k'}|| f_{0,j_1}^{k'} ||_{L^2}.
\end{equation*}
To  handle the part corresponding to the second sum in the
right-hand side of \eqref{ar50}, we use the space $Y_0$. It
suffices to prove that
\begin{equation}\label{ar59}
\begin{split}
\sum_{k''=k'-9}^5\sum_{j''=0}^\infty\sum_{j_2=0}^\infty
2^{\epsilon j_2}2^{j_2}
\big|\big|\mathcal{F}^{-1}[\eta_{j_2}(\tau_2)(f_{0,j_1}^{k'}
\ast\mathcal{F}(m_{k'',j''}))(\xi_2,\tau_2)]\big|\big|_{L^1_xL^2_t}&\\
\leq C2^{\epsilon j_1}2^{j_1-k'}|| f_{0,j_1}^{k'} ||_{L^2}&.
\end{split}
\end{equation}
As before, we may assume that $j_2$ satisfies the restriction
\eqref{ar33} and estimate
\begin{equation*}
\begin{split}
\big|\big|\mathcal{F}^{-1}[\eta_{j_2}(\tau_2)(f_{0,j_1}^{k'} \ast
\mathcal{F}(m_{k'',j''}))(\xi_2,\tau_2)]\big|\big|_{L^1_xL^2_t}&
\leq C||\mathcal{F}^{-1}(f_{0,j_1}^{k'} )||_{L^2_{x,t}}
||m_{k'',j''}||_{L^2_xL^\infty_t}\\
&\leq C2^{-80j''}2^{-k''}|| f_{0,j_1}^{k'} ||_{L^2_{x,t}},
\end{split}
\end{equation*}
using Plancherel theorem and \eqref{ar10}. The bound \eqref{ar59} follows.

We prove now the bound \eqref{ar4} in the case $k_2=0$ and
$k_1\in[1,10]$. As before, we may assume
$f_{k_1}^{\mathrm{high}}=f_{k_1,j_1}$ is an $L^2$ function
supported in $D_{k_1,j_1}$, $j_1\geq 0$,
$||A^\epsilon_{k_1}f_{k_1}^{\mathrm{high}}||_{Z_{k_1}}\approx
2^{\epsilon j_1}2^{j_1/2}\beta_{k_1,j_1}||f_{k_1,j}||_{L^2}\approx
2^{\epsilon j_1}2^{j_1}||f_{k_1,j}||_{L^2}$. We use the
decomposition \eqref{ar50} in the case $k'=1$.
The  proof of the bound  \eqref{ar4} is then identical to the
proof in the case considered before $k_1=0$,
$f_0^{\mathrm{high}}=f_{0,j_1}^{k'}$, $k'=1$.

Finally, in the case $k_1=0$, $k_2\in[1,10]$, we have the stronger
bound
\begin{equation}\label{ar40'}
\Big|\Big|\eta_{k_2}(\xi_2)A_{k_2}^\epsilon(\xi_2,\tau_2)\mathcal{F}
[m\cdot\mathcal{F}^{-1}(f_{0}^{\mathrm{high}})](\xi_2,\tau_2)\Big|\Big|_{Z_{k_2}}\leq
C ||A_0^\epsilon f_{0}^{\mathrm{high}}||_{\overline{Z}_{0}},
\end{equation}
where $\overline{Z}_0$ is defined in \eqref{def1'''}. The  proof
of this bound is identical to the proof of  \eqref{ar4} in the
case considered before $k_1=1$, $k_2\geq 1$.
\end{proof}

In some estimates the delicate structure of the spaces $Z_k$ is
not necessary. For $\alpha\in[-20,20]$ and $k\geq 1$ we define
\begin{equation}\label{def18}
\begin{split}
E_{k,\alpha}=&\{f\in L^2:\, f \text{ supported in }I_k\times\mathbb{R}\text{ and }\\
&||f||_{E_{k,\alpha}}:=2^{\alpha k}\sum_{j=0}^\infty
2^{j}||\eta_j(\tau)f(\xi,\tau)\,||_{L^2_{\xi,\tau}}<\infty\}.
\end{split}
\end{equation}
For $k=0$, for simplicity of notation we define
$E_{0,\alpha}=Z_0$. We notice that
\begin{equation}\label{vc50}
E_{k,4}\subseteq Z_{k}\subseteq E_{k,-4}\,\text{ for any
}k\in\mathbb{Z}_+.
\end{equation}

\newtheorem{Lemmab5}[Lemmab4]{Lemma}
\begin{Lemmab5}\label{Lemmab5}
(a) Assume $k_1\in\mathbb{Z}_+$, $k_2\in[1,\infty)\cap\mathbb{Z}$,
and $I_1\subseteq \widetilde{I}_{k_1}$, $I_2\subseteq
\widetilde{I}_{k_2}$ are intervals. Then, for $m\in S^\infty_{100}$, $\alpha\in[-20,20]$, $\epsilon\in\{-1,0\}$, and
$f_{k_1}\in E_{k_1,\alpha}$
\begin{equation}\label{vc51}
\begin{split}
\Big|\Big|\mathbf{1}_{I_2}&(\xi_2)(\tau_2+i)^\epsilon\cdot\mathcal{F}
[m\cdot\mathcal{F}^{-1}(\mathbf{1}_{I_1}(\xi_1)f_{k_1})]\Big|\Big|_{E_{k_2,\alpha}}\\
& \leq C [1+d(I_1,I_2)]^{-50}||m||_{S^\infty_{100}}\cdot||(\tau_1+i)^\epsilon
f_{k_1}||_{E_{k_1,\alpha}},
\end{split}
\end{equation}
where $d(I_1,I_2)$ denotes the distance between the sets $I_1$ and
$I_2$.

(b) Assume $k_1\in\mathbb{Z}_+$. Then, for $m\in S^\infty_{100}$, 
$\alpha\in[-20,20]$, $\epsilon\in\{-1,0\}$, and $f_{k_1}\in
E_{k_1,\alpha}$
\begin{equation}\label{vc52}
\Big|\Big|\eta_0(\xi_2)(\tau_2+i)^\epsilon\cdot\mathcal{F}
[m\cdot\mathcal{F}^{-1}(f_{k_1})]\Big|\Big|_{E_{0,\alpha}} \leq C
2^{-50k_1}||m||_{S^\infty_{100}}\cdot||(\tau_1+i)^\epsilon f_{k_1}||_{E_{k_1,\alpha}}.
\end{equation}
\end{Lemmab5}

\begin{proof}[Proof of Lemma \ref{Lemmab5}] We may assume $||m||_{S^\infty_{100}}=1$ and argue as in the
proof of Lemma \ref{Lemmab4}. We may assume $f_{k_1}=f_{k_1,j_1}$ is an $L^2$ function supported in $\widetilde{I}_{k_1}\times
\widetilde{I}_{j_1}$, $||(\tau_1+i)^\epsilon
f_{k_1}||_{E_{k_1,\alpha}}\geq C^{-1}2^{\alpha k_1}2^{\epsilon
j_1}2^{j_1}||f_{k_1,j_1}||_{L^2}$. With the notation in Lemma
\ref{Lemmab4}, we write
\begin{equation}\label{vc53}
m=\sum_{j''=0}^\infty m_{\leq 0,j''}+
\sum_{k''=1}^\infty\sum_{j''=0}^\infty m_{k'',j''}.
\end{equation}
For \eqref{vc51} it suffices to prove that
\begin{equation}\label{vc55}
\begin{split}
&2^{\alpha k_2}\sum_{j_2,j''\geq 0}2^{\epsilon
j_2}2^{j_2}\big|\big|\eta_{j_2}(\tau_2)\mathbf{1}_{I_2}(\xi_2)
\cdot[(\mathbf{1}_{I_1}(\xi_1)f_{k_1,j_1})\ast\mathcal{F}(m_{\leq
0,j''})]
\big|\big|_{L^2}\\
&+2^{\alpha k_2}\sum_{k''\geq 1}\sum_{j_2,j''\geq 0}2^{\epsilon
j_2}2^{j_2}\big|\big|\eta_{j_2}(\tau_2)\mathbf{1}_{I_2}(\xi_2)
\cdot[(\mathbf{1}_{I_1}(\xi_1)f_{k_1,j_1})\ast\mathcal{F}(m_{k'',j''})]
\big|\big|_{L^2}\\
&\leq C [1+d(I_1,I_2)]^{-50}2^{\alpha k_1}2^{\epsilon
j_1}2^{j_1}||f_{k_1,j_1}||_{L^2}.
\end{split}
\end{equation}

By examining the supports of the functions we see that
 the first sum in the left-hand side of \eqref{vc55} is nontrivial only if
 $d(I_1,I_2)\leq C$ (so $|k_1-k_2|\leq C$). In addition,
 $\eta_{j_2}(\tau_2)\mathbf{1}_{I_2}(\xi_2)
\cdot[(\mathbf{1}_{I_1}(\xi_1)f_{k_1,j_1})\ast\mathcal{F}(m_{\leq
0,j''})]\equiv 0$ unless
\begin{equation}\label{vc56}
|j_1-j_2|\leq C\,\text{ or }j_1,j_2\leq j''+C.
\end{equation}
Using Plancherel theorem and \eqref{ar10},
\begin{equation*}
\big|\big|(\mathbf{1}_{I_1}(\xi_1)f_{k_1,j_1})\ast\mathcal{F}(m_{\leq
0,j''}) \big|\big|_{L^2}\leq C2^{-80j''}||f_{k_1,j_1}||_{L^2}.
\end{equation*}
The bound \eqref{vc55} for the first sum follows easily. For the
second sum, we may assume that $2^{k''}\geq C^{-1}d(I_1,I_2)$ (so
$2^{k''}\geq C^{-1}2^{|k_1-k_2|}$) and that the restriction
\eqref{vc56} still holds. Using Plancherel theorem and
\eqref{ar10},
\begin{equation*}
\big|\big|(\mathbf{1}_{I_1}(\xi_1)f_{k_1,j_1})\ast\mathcal{F}(m_{k'',j''})
\big|\big|_{L^2}\leq C2^{-80k''}2^{-80j''}||f_{k_1,j_1}||_{L^2}.
\end{equation*}
The bound \eqref{vc55} for the second sum follows easily. This
completes the proof of part (a).

For part (b), we may assume $k_1\geq 10$ (in view of Lemma
\ref{Lemmab5}) and $f_{k_1}=f_{k_1,j_1}$ is as before. We
decompose $m$ as in \eqref{vc53}. For \eqref{vc52} it suffices to
prove that
\begin{equation}\label{vc60}
\begin{split}
\sum_{|k''-k_1|\leq 2}\sum_{j_2,j''\geq 0}2^{\epsilon
j_2}2^{j_2}\big|\big|\mathcal{F}^{-1}[\eta_{j_2}(\tau_2)\eta_0(\xi_2)
\cdot (f_{k_1,j_1}\ast\mathcal{F}(m_{k'',j''}))]&
\big|\big|_{L^1_xL^2_t}\\
\leq C2^{-50k_1}2^{\alpha k_1}2^{\epsilon
j_1}2^{j_1}||f_{k_1,j_1}||_{L^2}&.
\end{split}
\end{equation}
We may also assume that the restriction \eqref{vc56} holds. Using
Plancherel theorem and \eqref{ar10},
\begin{equation*}
\begin{split}
\big|\big|\mathcal{F}^{-1}[\eta_{j_2}(\tau_2)\eta_0(\xi_2) \cdot
(f_{k_1,j_1}\ast\mathcal{F}(m_{k'',j''}))]
\big|\big|_{L^1_xL^2_t}&\leq
C||m_{k'',j''}||_{L^2_xL^\infty_t}||\mathcal{F}^{-1}(f_{k_1,j_1})||_{L^2}\\
&\leq C2^{-80k''}2^{-80j''}||f_{k_1,j_1}||_{L^2}.
\end{split}
\end{equation*}
The bound \eqref{vc60} follows easily. This completes the proof of
part (b).
\end{proof}

We state now a slightly stronger form of Lemma \ref{Lemmab4} that will be used in the  next section.

\newtheorem{Lemmab6}[Lemmab4]{Corollary}
\begin{Lemmab6}\label{Lemmab6}
(a) If $k_1,k_2\in\mathbb{Z}_+$,
$\epsilon\in\{-1,0\}$, $f_{k_1}^{\mathrm{high}}\in
Z_{k_1}^{\mathrm{high}}$, and $m\in S^\infty_{100}$ then
\begin{equation}\label{arla}
\Big|\Big|\eta_{k_2}(\xi_2)A_{k_2}^\epsilon(\xi_2,\tau_2)\cdot
\mathcal{F}[m\cdot\mathcal{F}^{-1}(f_{k_1}^{\mathrm{high}})]\Big|\Big|_{Z_{k_2}}\leq
C 2^{-30|k_1-k_2|}||m||_{S^\infty_{100}}||A_{k_1}^{\epsilon}
f_{k_1}^{\mathrm{high}}||_{Z_{k_1}}.
\end{equation}

(b) If $k_2\in\mathbb{Z}_+$,
$\epsilon\in\{-1,0\}$, $f_{0}\in
\overline{Z}_0$, and $m'\in S^2_{100}$ then
\begin{equation}\label{arla2}
\Big|\Big|\eta_{k_2}(\xi_2)A_{k_2}^\epsilon(\xi_2,\tau_2)\cdot
\mathcal{F}[m'\cdot\mathcal{F}^{-1}(f_0)]\Big|\Big|_{Z_{k_2}}\leq
C 2^{-30k_2}||m'||_{S^2_{100}}||A_{0}^{\epsilon}
f_{0}||_{\overline{Z}_0}.
\end{equation}
\end{Lemmab6}

\begin{proof}[Proof of Corollary \ref{Lemmab6}] Part (a) follows from Lemma \ref{Lemmab4}, Lemma \ref{Lemmab5}, and
\eqref{vc50}. For part (b), we notice that $||m'_{\leq k'',j''}||_{L^2_xL^\infty_t}\leq C2^{-80j''}$ for any $k''\in\mathbb{Z}$, $j''\in\mathbb{Z}_+$. The bound \eqref{arla2} then follows from the proof of \eqref{ar40}, the bound \eqref{ar40'}, and the proof of Lemma \ref{Lemmab5} (a) with $k_1=1$.
\end{proof}

\section{Proof of Theorem \ref{Main1}}\label{proof}

In this section we complete the proof of Theorem \ref{Main1}. The
main ingredients are Lemma \ref{Lemmaw1}, Lemma \ref{Lemmab1},
Lemma \ref{Lemmab3}, Proposition \ref{Lemmac1}, Proposition
\ref{Lemmac2}, Proposition \ref{Lemmae1}, Proposition
\ref{Lemmae2}, Lemma \ref{Lemmab5}, and Corollary \ref{Lemmab6}. We
start by showing that the datum $e^{\pm
iU_0(.,0)}P_{\pm\mathrm{high}}\phi$ of the initial value problems
\eqref{fg20} and \eqref{fg21} are in $\widetilde{H}^\sigma$,
$\sigma\geq 0$.

\newtheorem{Lemmat1}{Lemma}[section]
\begin{Lemmat1}\label{Lemmat1}
Assume $U:\mathbb{R}\to\mathbb{R}$ satisfies the bounds
\begin{equation}\label{vc1}
||\partial_x^{\sigma_2}U||_{L^2_x}\leq 1\text{ for any
}\sigma_2\in[1,110]\cap\mathbb{Z}.
\end{equation}
Then, for any $\sigma\in[0,20]$ and $\phi\in H^\sigma$,
\begin{equation}\label{vc2}
||e^{\pm iU}P_{\pm\mathrm{high}}\phi||_{\widetilde{H}^\sigma}\leq C||\phi||_{H^\sigma}.
\end{equation}
\end{Lemmat1}

\begin{proof}[Proof of Lemma \ref{Lemmat1}] To fix the notation, assume that
the sign in the left-hand side of \eqref{vc2} is $+$. So we may
assume that $\widehat{\phi}$ is supported in the interval
$[2^{10},\infty)$. For any $k''\in\mathbb{Z}$ let
\begin{equation}\label{vc12}
V_{k''}=\mathcal{F}_1^{-1}\big[\chi_{k''}(\xi)\mathcal{F}_1[e^{iU(x)}]\big],
\end{equation}
and $V_{\leq k''}=\sum_{k'''\leq k''}V_{k'''}$. Using \eqref{vc1}
and the Sobolev imbedding theorem,
\begin{equation}\label{vc10}
||V_{\leq 0}||_{L^\infty_{x}}\leq C\text{ and
}||V_{k''}||_{L^2_x}+||V_{k''}||_{L^\infty_x} \leq
C2^{-80k''}\text{ for any }k''\geq 1.
\end{equation}

We turn now to the proof of \eqref{vc2}. For any $k_1\geq 10$
let $\phi_{k_1}=P_{k_1}\phi$. In view of the definitions, it
suffices to prove that
\begin{equation}\label{vc5}
\begin{cases}
&||P_{k_2}(e^{iU}\phi_{k_1})||_{L^2}\leq C2^{-40|k_1-k_2|}
||\phi_{k_1}||_{L^2}\text{ if }k_2\geq 1;\\
&||P_{0}(e^{iU}\phi_{k_1})||_{L^1}\leq
C2^{-40k_1}||\phi_{k_1}||_{L^2}.
\end{cases}
\end{equation}
For the first bound in \eqref{vc5}, if $|k_1-k_2|\leq 10$, then
$||P_{k_2}(e^{iU}\phi_{k_1})||_{L^2}\leq C||\phi_{k_1}||_{L^2}$ as
desired. If $|k_1-k_2|\geq 10$ then
\begin{equation*}
||P_{k_2}(e^{iU}\phi_{k_1})||_{L^2}
\leq\negmedspace\sum_{k''\geq |k_1-k_2|-C}||P_{k_2}(V_{k''}\phi_{k_1})||_{L^2}
\leq C\negmedspace\sum_{k''\geq |k_1-k_2|-C}||V_{k''}
||_{L^\infty}||\phi_{k_1}||_{L^2},
\end{equation*}
which suffices in view of \eqref{vc10}. For the second bound in
\eqref{vc5}, since $k_1\geq 10$,
\begin{equation*}
||P_{0}(e^{iU}\phi_{k_1})||_{L^1}\leq\sum_{|k''-k_1|\leq 2}
||P_{0}(V_{k''}\phi_{k_1})||_{L^1}\leq C\sum_{|k''-k_1|\leq 2}
||V_{k''}||_{L^2}||\phi_{k_1}||_{L^2}.
\end{equation*}
which suffices in view of \eqref{vc10}.
\end{proof}

We prove now our main bilinear estimate for functions in $F^\sigma$.

\newtheorem{Lemmat2}[Lemmat1]{Proposition}
\begin{Lemmat2}\label{Lemmat2}
If $m\in S^\infty_{110}$, $m'\in S^2_{110}$, $\sigma\in[0,20]$, and $u,v\in F^\sigma$ then
\begin{equation}\label{vc21}
\begin{split}
||\partial_x&(m\cdot uv)||_{N^\sigma}+||m'\cdot (uv)||_{N^\sigma}\\
&\leq C(||m||_{S^\infty_{110}}+||m'||_{S^2_{110}})(||u||_{F^\sigma}||v||_{F^0}+||u||_{F^0}||v||_{F^\sigma}).
\end{split}
\end{equation}
\end{Lemmat2}

\begin{proof}[Proof of Proposition \ref{Lemmat2}]
We show first that
\begin{equation}\label{vc22}
||\partial_x(uv)||_{N^\sigma}\leq C(||u||_{F^\sigma}||v||_{F^0}+
||u||_{F^0}||v||_{F^\sigma}).
\end{equation}
For $k\in\mathbb{Z}_+$  let
$F_k(\xi,\tau)=\eta_k(\xi)\mathcal{F}(u)(\xi,\tau)$  and
$G_k(\xi,\tau)=\eta_k(\xi)\mathcal{F}(v)(\xi,\tau)$. Then
\begin{equation*}
\begin{cases}
&||u||^2_{F^\sigma}=\sum_{k_1=0}^\infty2^{2\sigma k_1}
||(I-\partial_\tau^2)F_{k_1}||_{Z_{k_1}}^2;\\
&||v||^2_{F^\sigma}=\sum_{k_2=0}^\infty2^{2\sigma k_2}
||(I-\partial_\tau^2)G_{k_2}||_{Z_{k_2}}^2,
\end{cases}
\end{equation*}
and
\begin{equation*}
\eta_k(\xi)\mathcal{F}[\partial_x(u\cdot v)](\xi,\tau)=
C\xi\sum_{k_1,k_2\in\mathbb{Z}}\eta_k(\xi)[F_{k_1}\ast
G_{k_2}](\xi,\tau).
\end{equation*}
We  observe that $\eta_k(\xi)[F_{k_1}\ast G_{k_2}](\xi,\tau)\equiv
0$ unless
\begin{equation*}
\begin{cases}
&k_1\leq k-10\text{ and }k_2\in[k-2,k+2]\,\,\text {or}\\
&k_1\in[k-2,k+2]\text{ and }k_1\leq k-10\,\,\text {or}\\
&k_1,k_2\in[k-10,k+20]\,\,\text{or}\\
&k_1,k_2\geq k+10\text{ and }|k_1-k_2|\leq 2.
\end{cases}
\end{equation*}
For $k,k_1,k_2\in\mathbb{Z}$ let
\begin{equation*}
H_{k,k_1,k_2}(\xi,\tau)=\eta_k(\xi)A_k(\xi,\tau)^{-1}
\xi\cdot(F_{k_1}\ast G_{k_2})(\xi,\tau).
\end{equation*}
Using the definitions,
\begin{equation}\label{nu1}
\begin{split}
||&\partial_x(u\cdot v)||_{N^{\sigma}}^2= C\sum_{k\geq 0}2^{2\sigma k}
\Big|\Big|\sum_{k_1,k_2}H_{k,k_1,k_2}\Big|\Big|_{Z_k}^2.
\end{split}
\end{equation}
For $k\in\mathbb{Z}_+$ fixed we estimate, using Proposition
\ref{Lemmac1}, \ref{Lemmac2}, \ref{Lemmae1}, and \ref{Lemmae2},
\begin{equation*}
\begin{split}
\Big|\Big|\sum_{k_1,k_2}H_{k,k_1,k_2}\Big|\Big|_{Z_k}&\leq\sum_{|k_2-k|\leq 2}
\Big|\Big|\sum_{k_1\leq k-10}H_{k,k_1,k_2}\Big|\Big|_{Z_k}+\sum_{|k_1-k|\leq 2}
\Big|\Big|\sum_{k_2\leq k-10}H_{k,k_1,k_2}\Big|\Big|_{Z_k}\\
&+\sum_{k_1,k_2\in[k-10,k+20]}||H_{k,k_1,k_2}||_{Z_k}+
\sum_{k_1,k_2\geq k+10,\,|k_1-k_2|\leq 2}||H_{k,k_1,k_2}||_{Z_k}\\
&\leq C\big[\sum_{|k_2-k|\leq 2}||G_{k_2}||_{Z_{k_2}}\big]\cdot||u||_{F^0}
+C\big[\sum_{|k_1-k|\leq 2}||F_{k_1}||_{Z_{k_1}}\big]\cdot||v||_{F^0}\\
&+C\big[\sum_{|k_1-k|\leq 20}||F_{k_1}||_{Z_{k_1}}\big]\big[\sum_{|k_2-k|\leq 20}
||G_{k_2}||_{Z_{k_2}}\big]\\
&+C2^{-k/4}\big[\sum_{k_1\geq k}||F_{k_1}||_{Z_{k_1}}^2\big]^{1/2}
\big[\sum_{k_2\geq k}||G_{k_2}||_{Z_{k_2}}^2\big]^{1/2}.
\end{split}
\end{equation*}
The bound \eqref{vc22} follows. A similar estimate, using
Proposition \ref{Lemmae1} and \eqref{bc50}, shows that
\begin{equation}\label{mv5}
\begin{split}
||\eta_0(\xi)A_0(\xi,\tau)^{-1}\mathcal{F}(uv)||_{\overline{Z}_0}&+\big[\sum_{k\geq
1}2^{2\sigma k}||\eta_k(\xi)A_k(\xi,\tau)^{-1}\mathcal{F}(uv)||_{Z_k}^2\big]^{1/2}\\
&\leq C(||u||_{F^\sigma}||v||_{F^0}+ ||u||_{F^0}||v||_{F^\sigma}).
\end{split}
\end{equation}

We would like now to use the bound \eqref{arla} to include the
factor $m$. We may assume $||m||_{S^\infty_{110}}=1$. For any $u\in C(\mathbb{R}:H^{-2})$ we write
$u=u^{\mathrm{low}}+u^{\mathrm{high}}$,
\begin{equation*}
\begin{cases}
&u^{\mathrm{low}}=\sum\limits_{k\geq
0}\mathcal{F}^{-1}[\eta_k(\xi)\mathcal{F}(u)(\xi,\tau)\cdot\eta_{\leq
k-15}(\tau-\omega(\xi))]=\sum\limits_{k\geq 0}\mathcal{F}^{-1}(f_k^{\mathrm{low}});\\
&u^{\mathrm{high}}=\sum\limits_{k\geq
0}\mathcal{F}^{-1}[\eta_k(\xi)\mathcal{F}(u)(\xi,\tau)\cdot(1-\eta_{\leq
k-15}(\tau-\omega(\xi)))]=\sum\limits_{k\geq
0}\mathcal{F}^{-1}(f_k^{\mathrm{high}});
\end{cases}
\end{equation*}
Then, using \eqref{arla} with $\epsilon=0$,
\begin{equation}\label{mv1}
\begin{split}
||m\cdot u^{\mathrm{high}}||_{F^\sigma}^2&=\sum_{k\geq
0}2^{2\sigma
k}||\eta_k(\xi)\mathcal{F}[(t^2+1)mu^{\mathrm{high}}]||_{Z_k}^2\\
&\leq C\sum_{k\geq 0}2^{2\sigma k}\big[\sum_{k'\geq
0}||\eta_k(\xi)\mathcal{F}[(t^2+1)m\cdot\mathcal{F}^{-1}
(f_{k'}^{\mathrm{high}})]||_{Z_k}\big]^2\\
&\leq C\sum_{k\geq 0}2^{2\sigma k}\big[\sum_{k'\geq
0}2^{-30|k-k'|}||f_{k'}^{\mathrm{high}}||_{Z_{k'}}\big]^2\\
&\leq C||u||_{F^\sigma}^2
\end{split}
\end{equation}
for any $u\in F^\sigma$. A similar estimate, using \eqref{arla} with $\epsilon=1$, gives
\begin{equation}\label{mv2}
\begin{split}
||m\cdot w^{\mathrm{high}}||_{N^\sigma}\leq C||w||_{N^\sigma}
\end{split}
\end{equation}
for any $w\in N^\sigma$. We estimate now the first term in the left-hand side of
\eqref{vc21} by
\begin{equation}\label{mv4}
\begin{split}
&||\partial_x[(mu^{\mathrm{high}})v]||_{N^\sigma}
+||\partial_x[u^{\mathrm{low}}(mv^{\mathrm{high}})]||_{N^\sigma}\\
+&||m\cdot
\partial_x(u^{\mathrm{low}}v^{\mathrm{low}})||_{N^\sigma}
+||\partial_xm\cdot(u^{\mathrm{low}}v^{\mathrm{low}})||_{N^\sigma}
\end{split}
\end{equation}
In view of \eqref{vc22} and \eqref{mv1}, the first two terms in
\eqref{mv4} can be estimated by
$C(||u||_{F^\sigma}||v||_{F^0}+||u||_{F^0}||v||_{F^\sigma})$, as
desired. For the third term, we use the important observation that
the product of two low-modulation functions has high modulation:
\begin{equation}\label{mv10}
(u^{\mathrm{low}}v^{\mathrm{low}})^{\mathrm{low}}\equiv 0.
\end{equation}
Assuming \eqref{mv10}, the third term in \eqref{mv4} can be
estimated by
$C(||u||_{F^\sigma}||v||_{F^0}+||u||_{F^0}||v||_{F^\sigma})$,
using \eqref{vc22} and \eqref{mv2}. To prove \eqref{mv10}, we
write
\begin{equation*}
u^{\mathrm{low}}=\sum_{k\geq
15}\mathcal{F}^{-1}(f_k^{\mathrm{low}})\text{ and
}v^{\mathrm{low}}=\sum_{k\geq
15}\mathcal{F}^{-1}(g_k^{\mathrm{low}}),
\end{equation*}
where $f_k^{\mathrm{low}}$ and $g_k^{\mathrm{low}}$ are supported
in $\bigcup_{j\leq k-15}D_{k,j}$. For \eqref{mv10} it suffices to
prove that
\begin{equation*}
\eta_k(\xi)\eta_{\leq k-15}(\tau-\omega(\xi))(f_{k_1}\ast
g_{k_2})\equiv 0,\,k\geq 15,
\end{equation*}
which follows easily from \eqref{nev1} and \eqref{nev2}.

In view of \eqref{ar3''}, for \eqref{vc21}, it suffices to prove that if $||m'||_{S^2_{100}}=1$ then
\begin{equation}\label{mv11}
||m'\cdot (uv)||_{N^\sigma} \leq
C(||u||_{F^\sigma}||v||_{F^0}+||u||_{F^0}||v||_{F^\sigma}).
\end{equation}
We write $u=u^{\mathrm{high}}+u^{\mathrm{low}}$,
$v=v^{\mathrm{high}}+v^{\mathrm{low}}$ as before. Then, using \eqref{mv5}, the bound \eqref{arla2} with $k_2=0$, and Lemma \ref{Lemmab5} (b),
$$||P_0(m'\cdot uv)||_{N^\sigma}
\leq C(||u||_{F^\sigma}||v||_{F^0}+||u||_{F^0}||v||_{F^\sigma}).$$
Also, using \eqref{mv5} and \eqref{mv1} as before,
\begin{equation*}
\begin{split}
||(I-P_0)(m'\cdot u^{\mathrm{high}}v)||_{N^\sigma}
+&||(I-P_0)(m'\cdot u^{\mathrm{low}}v^{\mathrm{high}})||_{N^\sigma}\\
&\leq C(||u||_{F^\sigma}||v||_{F^0}+||u||_{F^0}||v||_{F^\sigma}).
\end{split}
\end{equation*}
Finally, using \eqref{mv2}, \eqref{arla2}, and the observation \eqref{mv10},
$$||(I-P_0)(m'\cdot u^{\mathrm{low}}v^{\mathrm{low}})||_{N^\sigma}
\leq C(||u||_{F^\sigma}||v||_{F^0}+||u||_{F^0}||v||_{F^\sigma}),$$
which completes the proof of \eqref{mv11}.
\end{proof}

To bound the error terms in the formulas \eqref{fg10} and
\eqref{fg11} of $E_+$ and $E_-$ we use the less demanding spaces
$E_{k,\alpha}$ defined in \eqref{def18} for $k\geq 1$,
\begin{equation*}
\begin{split}
E_{k,\alpha}=&\{f\in L^2:\, f \text{ supported in }I_k\times\mathbb{R}\text{ and }\\
&||f||_{E_{k,\alpha}}:=2^{\alpha k}\sum_{j=0}^\infty
2^{j}||\eta_j(\tau)f(\xi,\tau)\,||_{L^2_{\xi,\tau}}<\infty\}.
\end{split}
\end{equation*}. 
For $\sigma\geq 0$ and
$\alpha\in[-20,20]$ we define
\begin{equation*}
\begin{split}
F^{\sigma}_\alpha&=\Big\{u\in\mathcal{S}'(\mathbb{R}\times\mathbb{R}):
||u||_{F^{\sigma}_\alpha}^2:=\sum_{k=0}^\infty 2^{2\sigma k}
||\eta_k(\xi)(I-\partial_\tau^2)\mathcal{F}(u)||_{E_{k,\alpha}}^2
<\infty\Big\},
\end{split}
\end{equation*}
and
\begin{equation*}
\begin{split}
N^{\sigma}_\alpha&=\Big\{u\in\mathcal{S}'(\mathbb{R}\times\mathbb{R}):
||u||_{N^{\sigma}_\alpha}^2:=\sum_{k=0}^\infty 2^{2\sigma
k}||\eta_k(\xi)(\tau+i)^{-1}\mathcal{F}(u)||_{E_{k,\alpha}}^2<\infty\Big\}.
\end{split}
\end{equation*} 
In view of \eqref{vc50},
\begin{equation}\label{hh4}
F^\sigma_{6}\subseteq F^\sigma\subseteq F^{\sigma}_{-6}\text{ and }N^\sigma_{6}\subseteq N^\sigma\subseteq N^{\sigma}_{-6}.
\end{equation}

\newtheorem{Lemmat6}[Lemmat1]{Lemma}
\begin{Lemmat6}\label{Lemmat6}
If $m\in S^\infty_{110}$, $\sigma\in[0,20]$, $\alpha\in[-20,20]$, and
$u\in F^\sigma_{\alpha}$ then
\begin{equation}\label{hh100}
\begin{cases}
&||m\cdot u||_{F^\sigma_{\alpha}}\leq
C||m||_{S^\infty_{110}}||u||_{F^\sigma_{\alpha}};\\
&||m\cdot u||_{N^\sigma_{\alpha}}\leq
C||m||_{S^\infty_{110}}||u||_{N^\sigma_{\alpha}}.
\end{cases}
\end{equation}
\end{Lemmat6}

\begin{proof}[Proof of Lemma \ref{Lemmat6}] We may assume $||m||_{S^\infty_{110}}=1$. Let $f_{k'}=\eta_{k'}(\xi)\mathcal{F}(u)$, $k'\in\mathbb{Z}_+$. Using Lemma \ref{Lemmab5} with $\epsilon=0$, we have
\begin{equation*}
\begin{split}
||m\cdot u||_{F^\sigma_{\alpha}}^2&=\sum_{k\geq 0}2^{2\sigma k}||\eta_k(\xi)\mathcal{F}[(t^2+1)m\cdot u]||_{E_{k,\alpha}}^2\\
&\leq C\sum_{k\geq 0}2^{2\sigma k}\big[\sum_{k'\geq 0}||\eta_k(\xi)\mathcal{F}[(t^2+1)m\cdot
\mathcal{F}^{-1}(f_{k'})]||_{E_{k,\alpha}}\big]^2\\
&\leq C\sum_{k\geq 0}2^{2\sigma k}\big[\sum_{k'\geq
0}2^{-50|k-k'|}||f_{k'}||_{E_{k',\alpha}}\big]^2\\
&\leq C||u||_{F^\sigma_{\alpha}}^2.
\end{split}
\end{equation*}
A similar estimate using Lemma \ref{Lemmab5} with $\epsilon=1$
gives the second bound in
\eqref{hh100}.
\end{proof}

\newtheorem{Lemmat3}[Lemmat1]{Lemma}
\begin{Lemmat3}\label{Lemmat3}
(a) Assume that $I\neq
I'\in\{[-2^{10},2^{10}],[2^{10},\infty),(-\infty,-2^{10}]\}$ and
$m\in S^\infty_{110}$. Then, for any $\sigma\in[0,20]$ and
$u\in F^\sigma_{-10}$,
\begin{equation}\label{mv20}
\begin{cases}
&||P_I[m\cdot P_{I'}(u)]||_{F^\sigma_{10}}\leq
C||m||_{S^\infty_{110}}||u||_{F^\sigma_{-10}};\\
&||P_I[m\cdot P_{I'}(u)]||_{N^\sigma_{10}}\leq
C||m||_{S^\infty_{110}}||u||_{N^\sigma_{-10}},
\end{cases}
\end{equation}
where $P_I$ denotes the operator defined by the multiplier
$(\xi,\tau)\to\mathbf{1}_I(\xi)$.

(b) In addition, for any $\sigma\in[0,20]$ and $u\in
F^\sigma_{-10}$,
\begin{equation}\label{mv21}
||\partial_x^2P_-(m\cdot P_{+\mathrm{high}}(u))||_{F^\sigma_{10}}+
||\partial_x^2P_+(m\cdot
P_{-\mathrm{high}}(u))||_{F^\sigma_{10}}\leq
C||m||_{S^\infty_{110}}||u||_{F^\sigma_{-10}}.
\end{equation}
\end{Lemmat3}

\begin{proof}[Proof of Lemma \ref{Lemmat3}] We may assume $||m||_{S^\infty_{110}}=1$ and use Lemma
\ref{Lemmab5} and the definitions. The main observation is that if
$k,k'\in\mathbb{Z}_+$ then
$d(I\cap\widetilde{I}_k,I'\cap\widetilde{I}_{k'})\geq
C^{-1}(2^k+2^{k'})$. Let
$f_{k'}=\eta_{k'}(\xi)\cdot\mathcal{F}(u)(\xi,\tau)$. Using Lemma
\ref{Lemmab5} with $\epsilon=0$, we have
\begin{equation*}
\begin{split}
||P_I[m\cdot P_{I'}&(u)]||_{F^\sigma_{10}}^2=\sum_{k\geq
0}2^{2\sigma
k}||\eta_k(\xi)\mathbf{1}_I(\xi)\cdot\mathcal{F}[(t^2+1)m\cdot
P_{I'}(u)]||_{E_{k,10}}^2\\
&\leq C\sum_{k\geq 0}2^{2\sigma k}\big[\sum_{k'\geq
0}||\eta_k(\xi)\mathbf{1}_I(\xi)\cdot\mathcal{F}[(t^2+1)m\cdot
\mathcal{F}^{-1}(\mathbf{1}_{I'}f_{k'})]||_{E_{k,10}}\big]^2\\
&\leq C\sum_{k\geq 0}2^{2\sigma k}\big[\sum_{k'\geq
0}(2^k+2^{k'})^{-50}2^{20k'}2^{-\sigma k'}2^{\sigma
k'}||f_{k'}||_{E_{k',-10}}\big]^2\\
&\leq C||u||_{F^\sigma_{-10}}^2.
\end{split}
\end{equation*}
A similar estimate using Lemma \ref{Lemmab5} with $\epsilon=1$
gives the second bound in
\eqref{mv20}.

For part (b) the same argument as before works, except for the
dyadic piece corresponding to $k=0$ (in the left-hand side). To
handle this dyadic piece we need the additional observation
\begin{equation*}
||\xi^2\mathbf{1}_{\pm}(\xi)\eta_0(\xi)f||_{Z_0}\leq
||\xi^2\eta_0(\xi)f||_{X_0}\leq
C||\eta_0(\xi)f||_{\overline{Z}_0}\leq C||\eta_0(\xi)f||_{Z_0},
\end{equation*}
where $\mathbf{1}_{\pm}$ denotes the characteristic function of
the interval $\{\xi:\pm\xi\in[0,\infty)\}$.
\end{proof}
\medskip

We can now analyze the nonlinear terms $E_+$, $E_-$, and $E_0$ in
\eqref{fg10}, \eqref{fg11}, and \eqref{fg12}. We assume that $u_0,
U_0:\mathbb{R}\times[-2,2]\to\mathbb{R}$ are fixed functions that
satisfy the bounds (compare with \eqref{fg0} and \eqref{fg9})
\begin{equation}\label{po1}
\begin{cases}
&||\partial_t^{\sigma_1}\partial_x^{\sigma_2}
u_0||_{L^2_{x,t}}\leq \delta\text{ for any }\sigma_1,\sigma_2\in[0,120]\cap\mathbb{Z};\\
&||\partial_t^{\sigma_1}\partial_x^{\sigma_2}
U_0||_{L^2_{x,t}}\leq\delta\text{ for any
}\sigma_1\in[0,120]\cap\mathbb{Z},\sigma_2\in[1,120]\cap\mathbb{Z},
\end{cases}
\end{equation}
for some $\delta\ll1$, and $E_+$, $E_-$, and $E_0$ are defined as in \eqref{fg10},
\eqref{fg11}, and \eqref{fg12}. For simplicity of notation, let
$\mathbf{w}=(w_+,w_-,w_0)$ and
\begin{equation*}
\mathbf{E}(\mathbf{w})=\big(E_+(w_+,w_-,w_0),E_-(w_+,w_-,w_0),
E_0(w_+,w_-,w_0)\big).
\end{equation*}
For any Banach space $B$ let
$||\mathbf{w}||_B=||w_+||_B+||w_-||_B+||w_0||_B$ and
\begin{equation*}
||\mathbf{E}(\mathbf{w})||_B=||E_+(w_+,w_-,w_0)||_B+||E_-(w_+,w_-,w_0)||_B
+||E_0(w_+,w_-,w_0)||_B.
\end{equation*}

\newtheorem{Lemmat4}[Lemmat1]{Proposition}
\begin{Lemmat4}\label{Lemmat4}
Assume that $\sigma\in[0,20]$, $u_0,U_0$ satisfy \eqref{po1},
$\mathbf{w},\mathbf{w}'\in F^\sigma$, and
$\psi:\mathbb{R}\to[0,1]$ is the smooth function defined in
section \ref{linear}. Then
\begin{equation}\label{po2}
\begin{split}
||\psi(t)[\mathbf{E}(\mathbf{w})-\mathbf{E}(\mathbf{w}')]||_{N^\sigma}\leq
&C||\mathbf{w}-\mathbf{w}'||_{F^\sigma}
(\delta+||\mathbf{w}||_{F^0}+||\mathbf{w}'||_{F^0})\\
+&C||\mathbf{w}-\mathbf{w}'||_{F^0}
(||\mathbf{w}||_{F^\sigma}+||\mathbf{w}'||_{F^\sigma}).
\end{split}
\end{equation}
\end{Lemmat4}

\begin{proof}[Proof of Proposition \ref{Lemmat4}] Let $T_{i,+}$ and $T_{i,-}$, $i\in\{1,2,3,4,5\}$ denote the terms in line $i$ in the formulas \eqref{fg10} and \eqref{fg11}. To control $||\psi(t)[T_{1,+}(\mathbf{w})-T_{1,+}(\mathbf{w'})]||_{N^\sigma}$ it suffices to prove that
\begin{equation}\label{hh1}
||m\cdot P_{+\mathrm{high}}(\partial_x(m'uv))||_{N^{\sigma}}\leq C(||u||_{F^\sigma}||v||_{F^0}+||u||_{F^0}||v||_{F^\sigma}),
\end{equation}
for any functions $u,v\in F^\sigma$, where $||m||_{S^\infty_{110}}=||m'||_{S^\infty_{110}}=1$. We bound the left-hand of \eqref{hh1} by
\begin{equation}\label{hh2}
\begin{split}
&||(P_{-\mathrm{high}}+P_{\mathrm{low}})[m\cdot P_{+\mathrm{high}}(\partial_x(m'uv))]||_{N^{\sigma}}\\
+&||P_{+\mathrm{high}}[m\cdot (P_{-\mathrm{high}}+P_{\mathrm{low}})(\partial_x(m'uv))]||_{N^{\sigma}}\\
+&||P_{+\mathrm{high}}[m\cdot(\partial_x(m'uv))]||_{N^{\sigma}}.
\end{split}
\end{equation}
For the first two terms in \eqref{hh2} we  use Lemma \ref{Lemmat3} (a), Proposition \ref{Lemmat2}, and \eqref{hh4}. For the third term in \eqref{hh2} we use Proposition \ref{Lemmat2} and \eqref{ar3''}. The bound \eqref{hh1} follows.

To control $||\psi(t)[T_{2,+}(\mathbf{w})-T_{2,+}(\mathbf{w'})]||_{N^\sigma}$ it suffices to prove that
\begin{equation}\label{hh6}
||m\cdot P_{+\mathrm{high}}[\partial_x(u_0\cdot P_{-\mathrm{high}}(m'u))]||_{N^{\sigma}}+||m\cdot P_{+\mathrm{high}}[\partial_x(u_0\cdot P_{\mathrm{low}}(m'u))]||_{N^{\sigma}}\leq C\delta||u||_{F^\sigma}
\end{equation}
for any $u\in F^\sigma$, where $||m||_{S^\infty_{110}}=||m'||_{S^\infty_{110}}=1$. We use  Lemma \ref{Lemmat6}, Lemma \ref{Lemmat3} (a), and \eqref{hh4}. The -first term in the left-hand side of \eqref{hh6} is dominated by
\begin{equation*}
\begin{split}
||m\cdot P_{+\mathrm{high}}[\partial_x(u_0\cdot P_{-\mathrm{high}}(m'u))]||_{N^{\sigma}_6}&\leq C||P_{+\mathrm{high}}[u_0\cdot P_{-\mathrm{high}}(m'u)]||_{F^{\sigma}_7}\\
&\leq C\delta||m'u||_{F^\sigma_{-10}}\leq C\delta||u||_{F^\sigma},
\end{split}
\end{equation*}
as desired. The bound for the second term is similar. Also, the bound for $||\psi(t)[T_{3,+}(\mathbf{w})-T_{3,+}(\mathbf{w'})]||_{N^\sigma}$ follows in the same way.

To control $||\psi(t)[T_{4,+}(\mathbf{w})-T_{4,+}(\mathbf{w'})]||_{N^\sigma}$ it suffices to prove that
\begin{equation}\label{hh8}
||\partial_x^2P_{-}[((I-P_0)e^{iU_0})\cdot P_{+\mathrm{high}}(mu)]||_{N^{\sigma}}\leq C\delta||u||_{F^\sigma}
\end{equation}
for any $u\in F^\sigma$, where $||m||_{S^\infty_{110}}=1$. This follows as before, using Lemma \ref{Lemmat6}, Lemma \ref{Lemmat3} (b), \eqref{hh4}, and the observation that $||(I-P_0)e^{iU_0}||_{S^\infty_{110}}\leq C\delta$.

To control $||\psi(t)[T_{5,+}(\mathbf{w})-T_{5,+}(\mathbf{w'})]||_{N^\sigma}$ it suffices to prove that
\begin{equation}\label{hh9}
||P_+\partial_xu_0\cdot u||_{N^{\sigma}}\leq C\delta||u||_{F^\sigma}.
\end{equation}
for any $u\in F^\sigma$. We bound the left-hand side of \eqref{hh9} by
\begin{equation}\label{hh15}
||(I-P_0)[(I-P_0)(P_+\partial_xu_0)u]||_{N^{\sigma}}+||(I-P_0)[(P_0P_+\partial_xu_0)u]||_{N^{\sigma}}+||P_0[P_+\partial_xu_0\cdot u]||_{N^{\sigma}}.
\end{equation}
For the first term in \eqref{hh15} we use Proposition \ref{Lemmat2} with $m\equiv 1$, $m'\equiv 0$. For the second term in \eqref{hh15} we use the bound \eqref{hh91}. For the third term in \eqref{hh15} we use Lemma \ref{Lemmat6}:
\begin{equation*}
||P_0[P_+\partial_xu_0\cdot u]||_{N^{\sigma}}\leq C||P_+\partial_xu_0\cdot u||_{F^{\sigma}_{-10}}\leq C\delta||u||_{F^\sigma_{-10}},
\end{equation*}
as desired.

The proofs for the terms $T_{i,-}$ are identical. To control $||\psi(t)[E_0(\mathbf{w})-E_0(\mathbf{w'})]||_{N^\sigma}$ it suffices to prove that
\begin{equation}\label{hh25}
\begin{cases}
&||P_{\mathrm{low}}\partial_x(muv)||_{N^\sigma}\leq C(||u||_{F^\sigma}||v||_{F^0}+||u||_{F^0}||v||_{F^\sigma});\\
&||P_{\mathrm{low}}\partial_x(mu_0u)||_{N^\sigma}\leq C\delta||u||_{F^\sigma}
\end{cases}
\end{equation}
for any functions $u,v\in F^\sigma$, where $||m||_{S^\infty_{115}}=1$. For the first bound in \eqref{hh25} we use Proposition \ref{Lemmat2}. For the second bound we use Lemma  \ref{Lemmat6} and the observation $||mu_0||_{S^\infty_{110}}\leq C\delta$. This completes the proof of Proposition \ref{Lemmat4}.
\end{proof}

\begin{proof}[Proof of Theorem \ref{Main1}] For any interval $I\subseteq\mathbb{R}$ and $\sigma\geq 0$ we define the Banach spaces
\begin{equation*}
\begin{cases}
&F^\sigma(I)=\{u\in \mathcal{S}'(\mathbb{R}\times I):||u||_{F^\sigma(I)}:=\inf\limits_{\widetilde{u}\equiv u\text{ on }\mathbb{R}\times I}||\widetilde{u}||_{F^\sigma}<\infty\};\\
&N^\sigma(I)=\{u\in \mathcal{S}'(\mathbb{R}\times I):||u||_{N^\sigma(I)}:=\inf\limits_{\widetilde{u}\equiv u\text{ on }\mathbb{R}\times I}||\widetilde{u}||_{N^\sigma}<\infty\}.
\end{cases}
\end{equation*}
With this notation, the estimates in Lemma \ref{Lemmab1} and Lemma \ref{Lemmab3} become
\begin{equation}\label{jm1}
||W(t-t_0)\phi||_{F^\sigma([t_0-a,t_0+a])}\leq C||\phi||_{\widetilde{H}^\sigma},
\end{equation}
and
\begin{equation}\label{jm2}
\Big|\Big|\int_{t_0}^tW(t-s)(u(s))\,ds\Big|\Big|_{F^\sigma([t_0-a,t_0+a])}\leq  C||u||_{N^\sigma([t_0-a,t_0+a])},
\end{equation}
for any $\sigma\in[0,20]$, $t_0\in\mathbb{R}$ and $a\in[0,5/4]$. The estimate in Proposition \ref{Lemmat4} becomes
\begin{equation}\label{jm3}
\begin{split}
||\mathbf{E}(\mathbf{w})-\mathbf{E}(\mathbf{w}')||_{N^\sigma(I)}\leq
&C||\mathbf{w}-\mathbf{w}'||_{F^\sigma(I)}
(\delta+||\mathbf{w}||_{F^0(I)}+||\mathbf{w}'||_{F^0(I)})\\
+&C||\mathbf{w}-\mathbf{w}'||_{F^0(I)}
(||\mathbf{w}||_{F^\sigma(I)}+||\mathbf{w}'||_{F^\sigma(I)}),
\end{split}
\end{equation}
for any $\sigma\in[0,20]$ and $I\subseteq[-5/4,5/4]$, provided that \eqref{po1} holds.

Assume that $u_0$, $U_0$ are fixed and satisfy \eqref{po1}. For data $\Phi=(\phi_+,\phi_-,\phi_0)\in \widetilde{H}^{20}$ with the property
\begin{equation}\label{jm4}
||\Phi||_{\widetilde{H}^0}\leq\delta,
\end{equation}
we consider the vector-valued initial-value problem
\begin{equation}\label{jm5}
\begin{cases}
&(\partial_t+\mathcal{H}\partial_x^2)\mathbf{v}=\mathbf{E}(\mathbf{v})\text{ on }\mathbb{R}\times[-5/4,5/4];\\
&\mathbf{v}(0)=\Phi.
\end{cases}
\end{equation}
We can construct a solution of \eqref{jm5} by iteration: let $\mathbf{v}^0=(0,0,0)$ and let
\begin{equation}\label{jm9}
\mathbf{v}^{k+1}=W(t)\Phi+\int_0^tW(t-s)(\mathbf{E}(\mathbf{v}^{k})(s))\,ds,\,k=0,1,\ldots.
\end{equation}
In view of \eqref{jm1}, \eqref{jm2}, \eqref{jm3}, and \eqref{jm4}, $||\mathbf{v}^{k}||_{F^0([-5/4,5/4])}\leq C\delta$ for any $k\geq 0$. Thus, using \eqref{jm1}, \eqref{jm2}, \eqref{jm3}, and \eqref{jm4}  again,
\begin{equation}\label{jm7}
||\mathbf{v}^{k+1}-\mathbf{v}^{k}||_{F^0([-5/4,5/4])}\leq (C\delta)^{k+1}\text{ for any }k=0,1,\ldots.
\end{equation}
Using \eqref{jm1}, \eqref{jm2}, \eqref{jm3}, \eqref{jm4}, and \eqref{jm7} we obtain $||\mathbf{v}^{k}||_{F^{\sigma}([-5/4,5/4])}\leq C||\Phi||_{\widetilde{H}^{\sigma}}$, $\sigma\in[0,20]$, and then
\begin{equation*}
||\mathbf{v}^{k+1}-\mathbf{v}^{k}||_{F^{\sigma}([-5/4,5/4])}\leq (C\delta)^{k}||\Phi||_{\widetilde{H}^{\sigma}}\text{ for any }k=0,1,\ldots.
\end{equation*}
Thus the sequence $\mathbf{v}^k$ converges in the space $F^{20}([-5/4,5/4])$ to a function $\mathbf{v}=\mathbf{v}(\Phi)$. In addition, for any $\sigma\in[0,20]$,
\begin{equation}\label{jm8}
||\mathbf{v}(\Phi)||_{F^{\sigma}([-5/4,5/4])}\leq C||\Phi||_{\widetilde{H}^{\sigma}},
\end{equation}
$\mathbf{v}(\Phi)\in C([-5/4,5/4]:\widetilde{H}^{20})$ (using \eqref{hh80}), $\mathbf{v}(\Phi)$ solves the initial-value problem \eqref{jm5}, and if $||\Phi||_{\widetilde{H}^0},||\Phi'||_{\widetilde{H}^0}\leq\delta$ then
\begin{equation}\label{jm10}
\begin{split}
&||\mathbf{v}(\Phi)-\mathbf{v}(\Phi')||_{F^\sigma([-5/4,5/4])}\\
&\leq C||\Phi-\Phi'||_{\widetilde{H}^\sigma}+C(||\Phi||_{\widetilde{H}^\sigma}+||\Phi'||_{\widetilde{H}^\sigma})||\mathbf{v}(\Phi)-\mathbf{v}(\Phi')||_{F^0([-5/4,5/4])}.
\end{split}
\end{equation}
In particular, when $\sigma=0$, $||\mathbf{v}(\Phi)-\mathbf{v}(\Phi')||_{F^0([-5/4,5/4])}\leq C||\Phi-\Phi'||_{\widetilde{H}^0}$.

Assume now that we start with data $\phi\in H^\infty_r$ with the property
\begin{equation}\label{jm12}
||\phi||_{L^2}\leq\delta_0=\delta/C, \text{ where }C\text{ is sufficiently large}.
\end{equation}
We construct the functions $u_0$, $\widetilde{u}$, $U_0$, $\mathbf{w}=(w_+,w_-,w_0)$, and $$\Phi=(\phi_+,\phi_-,\phi_0)=(e^{iU_0(.,0)}P_{+\mathrm{high}}\phi,e^{-iU_0(.,0)}P_{-\mathrm{high}}\phi,0)$$ 
as in section \ref{gauge}. Clearly, \eqref{jm4} holds due to Lemma \ref{Lemmat1}, and $\Phi\in\widetilde{H}^{20}$. We show now that
\begin{equation}\label{jm15}
\mathbf{w}\equiv \mathbf{v}(\Phi)\text{ in }\mathbb{R}\times[-1,1],
\end{equation}
where $\mathbf{v}(\Phi)$ is constructed as before. This is somewhat delicate since it is not clear how to show algebraically that the function $e^{-iU_0}v_++e^{iU_0}v_-+v_0+u_0$ is a solution of the original initial-value problem.    

To prove \eqref{jm15} we show first that
\begin{equation}\label{jm16}
||\mathbf{w}(t)||_{\widetilde{H}^0}\leq C\delta_0\text{ for any }t\in[-5/4,5/4].
\end{equation}
For the functions $w_+$ and $w_-$ this follows directly using the definition \eqref{fg4} and Lemma \ref{Lemmat1}, since, in view of the conservation law \eqref{conserve},
\begin{equation}\label{jm17}
||\widetilde{u}||_{L^\infty_tL^2_x}+||u_0||_{L^\infty_tL^2_x}\leq 3\delta_0 \text{ for any }t\in[-5/4,5/4].
\end{equation}
To prove \eqref{jm16} for the function $w_0$ we use first the definition \eqref{fg4} and \eqref{jm17}, so it suffices to prove that 
\begin{equation}\label{jm25}
||\eta_0(\xi)\mathcal{F}_1(w_0(t))(\xi)||_{B_0}\leq C\delta_0,\,t\in[-5/4,5/4].
\end{equation}
For this we use the equation \eqref{fg3} (notice $w_0(0)\equiv 0$). It suffices to prove that
\begin{equation}\label{jm19}
||\eta_0(\xi)\xi^2\mathrm{sgn}(\xi)\mathcal{F}_1(\widetilde{u}(t))(\xi)||_{B_0}+||\eta_0(\xi)\xi\mathcal{F}_1(\widetilde{u}(t)(\widetilde{u}(t)/2+u_0(t)))(\xi)||_{B_0}\leq C\delta_0,
\end{equation}
for any $t\in[-5/4,5/4]$. We bound the first term in \eqref{jm19} by
\begin{equation*}
\sum_{k'\leq 1}2^{-k'}||\chi_{k'}(\xi)\xi^2\mathcal{F}_1(\widetilde{u}(t))(\xi)||_{L^2_\xi}\leq C||\widetilde{u}(t)||_{L^2_x}\leq C\delta_0,
\end{equation*}
as desired. We bound the second term in \eqref{jm19} by
\begin{equation*}
||\mathcal{F}_1^{-1}[\eta_0(\xi)\xi\mathcal{F}_1(\widetilde{u}(t)(\widetilde{u}(t)/2+u_0(t)))(\xi)]||_{L^1_x}\leq C||\widetilde{u}(t)||_{L^2_x}||\widetilde{u}(t)/2+u_0(t)||_{L^2_x},
\end{equation*}
which suffices in view of \eqref{jm17}. This completes the proof of \eqref{jm16}. 

Next, we show that there is $\varepsilon=\varepsilon(||\phi||_{H^{100}})$ with the property that
\begin{equation}\label{jm20}
||\mathbf{w}||_{F^0([t_0-\varepsilon,t_0+\varepsilon])}\leq C\delta_0\text{ for any }t_0\in[-1,1].
\end{equation}
Let $\mathbf{g}=\psi(t)(\partial_t+\mathcal{H}\partial_x^2)\mathbf{w}$. In view of \eqref{jm1}, \eqref{jm2}, \eqref{jm16}, and \eqref{hh4}, for \eqref{jm20} it suffices to prove that
\begin{equation}\label{jm40}
||\psi((t-t_0)/\varepsilon)\cdot \mathbf{g}||_{N^0_6}\leq C(||\phi||_{H^{100}})\varepsilon^{1/4}.
\end{equation}
We show first that for any $t\in[-5/4,5/4]$
\begin{equation}\label{jm21}
||(I-\partial_t^2)\mathbf{g}(t)||_{\widetilde{H}^{20}}\leq C(||\phi||_{H^{100}}).
\end{equation}
For \eqref{jm21} we notice first that $\mathcal{H}\partial_x^2:\widetilde{H}^{\sigma}\to \widetilde{H}^{\sigma-2}$ is a bounded operator. Thus it suffices to prove that $||\partial_t^{\sigma}\mathbf{w}||_{\widetilde{H}^{50}}\leq C(||\phi||_{H^{100}})$, $\sigma=0,1,2,3$. For $w_+$ and $w_-$ this is clear using the definitions $w_{\pm}=e^{\pm iU_0}P_{\pm\mathrm{high}}\widetilde{u}$ and Lemma \ref{Lemmat1}. For $w_0$ this follows using the identity \eqref{fg3}
\begin{equation*}
\partial_tw_0=-\mathcal{H}\partial_x^2w_0-P_{\mathrm{low}}\partial_x((u_0+\widetilde{u}/2)\cdot \widetilde{u}),
\end{equation*}
the bound \eqref{jm16}, and the same argument as in the proof of \eqref{jm25}. This completes the proof of \eqref{jm21}. To pass from \eqref{jm21} to \eqref{jm40}, we may assume $t_0=0$ and $\mathbf{g}=g$ is scalar valued. It suffices to prove that
\begin{equation}\label{jm60}
||\psi(t/\varepsilon)\cdot g||_{N^0_6}\leq C\varepsilon^{1/4}||(I-\partial_t^2)g||_{L^1_t\widetilde{H}^{20}}.
\end{equation}
In view of the $L^1_t$ norm in the right-hand side of \eqref{jm60}, we may assume that $g(x,t)=h(x)K(t-t_0)$, where $K(t)=\int_\mathbb{R}(\tau^2+1)^{-1}e^{it\tau}\,d\tau$ and $||(I-\partial_t^2)g||_{L^1_t\widetilde{H}^{20}}\approx ||h|||_{\widetilde{H}^{20}}$. The bound \eqref{jm60} then follows easily from the definitions.  

We can now complete the proof of \eqref{jm15}. Assume $\mathbf{w}(t_0)=\mathbf{v}(t_0)=\Psi$ for some $t_0\in[-1,1]$. Then, for $t\in[t_0-\varepsilon,t_0+\varepsilon]$ we write
\begin{equation*}
\begin{cases}
&\mathbf{w}(t)=W(t-t_0)\Psi+\int_{t_0}^tW(t-s)(\mathbf{E}(\mathbf{w})(s))\,ds;\\
&\mathbf{v}(t)=W(t-t_0)\Psi+\int_{t_0}^tW(t-s)(\mathbf{E}(\mathbf{v})(s))\,ds.
\end{cases}
\end{equation*}
We subtract the two identities and use \eqref{jm2}, \eqref{jm3}, \eqref{jm8} (all with $\sigma=0$), and \eqref{jm20}. The result is
\begin{equation*}
||\mathbf{v}-\mathbf{w}||_{F^0([t_0-\varepsilon,t_0+\varepsilon])}\leq C||\mathbf{E}(\mathbf{v})-\mathbf{E}(\mathbf{w})||_{N^0([t_0-\varepsilon,t_0+\varepsilon])}\leq C\delta||\mathbf{v}-\mathbf{w}||_{F^0([t_0-\varepsilon,t_0+\varepsilon])}.
\end{equation*}
So $\mathbf{v}\equiv\mathbf{w}$ in $\mathbf{R}\times[t_0-\varepsilon,t_0+\varepsilon]$. Since $\mathbf{w}(0)=\mathbf{v}(0)=\Phi$, this suffices to prove \eqref{jm15}.

We prove now part (a) of the theorem. Assume that 
\begin{equation*}
\phi_n\in H^\infty_r\text{ and }\lim_{n\to\infty}\phi_n=\phi\text{ in }L^2.
\end{equation*}
By rescaling\footnote{The smooth flow has invariance property $S^\infty(\phi_\lambda)=[S^\infty(\phi)]_{\lambda}$, where $\phi_\lambda(x)=\lambda\phi(\lambda x)$ and $u_\lambda(x,t)=\lambda u(\lambda x,\lambda^2t)$.}, we may assume $||\phi||_{L^2}\leq\delta_0/2$, as in \eqref{jm12}. By using the conservation law \eqref{conserve}, we may assume $T=1$. It suffices to prove that for any $\epsilon>0$
\begin{equation}\label{jm30}
||S^\infty_1(\phi_n)-S^\infty_1(\phi_m)||_{L^\infty_tL^2_x}\leq\varepsilon\text{ for }m,n\text{ sufficiently large}.
\end{equation}
We fix $M=M(\phi,\varepsilon)$ sufficiently large and define $\widehat{\phi^M}(\xi)=\mathbf{1}_{[-M,M]}(\xi)\widehat{\phi}(\xi)$ and $\widehat{\phi_n^M}(\xi)=\mathbf{1}_{[-M,M]}(\xi)\widehat{\phi_n}(\xi)$. It is known that the flow map $S^\infty_1$ extends continuously on, say, $H^2_r$ (see, for example, \cite{Po}). Since $\lim_{n\to\infty}\phi_n^M=\phi^M$ in $H^2_r$,
\begin{equation*}
\lim_{n,m\to\infty}||S^\infty_1(\phi_n^M)-S^\infty_1(\phi_m^M)||_{L^\infty_tH^2_x}=0.
\end{equation*}

We estimate now $||S^\infty_1(\phi_n)-S^\infty_1(\phi_n^M)||_{L^\infty_tL^2_x}$. As in section \ref{gauge}, we construct $u_{0,n}$, $U_{0,n}$ (which are identical for both functions $\phi_n$ and $\phi_n^M$), $$\Phi_n=(e^{iU_{0,n}}P_{+\mathrm{high}}\phi_n,e^{-iU_{0,n}}P_{-\mathrm{high}}\phi_n,0),$$ and $$\Phi_n^M=(e^{iU_{0,n}}P_{+\mathrm{high}}\phi_n^M,e^{-iU_{0,n}}P_{-\mathrm{high}}\phi_n^M,0).$$ Using Lemma \ref{Lemmaw1}, the identity \eqref{jm15}, \eqref{jm10} with $\sigma=0$, \eqref{hh80}, and Lemma \ref{Lemmat1}
\begin{equation*}
\begin{split}
||S^\infty_1(\phi_n)-S^\infty_1(\phi_n^M)||_{L^\infty_tL^2_x}&\leq C||\mathbf{v}(\Phi_n)-\mathbf{v}(\Phi_n^M)||_{L^\infty_{t\in[-1,1]}L^2_x}\leq C||\Phi_n-\Phi_n^M||_{\widetilde{H}^0}\\
&\leq C||\phi_n-\phi_n^M||_{L^2}\leq C(||\phi-\phi^M||_{L^2}+||\phi-\phi_n||_{L^2}).
\end{split}
\end{equation*}
The bound \eqref{jm30} follows if $M=M(\phi,\varepsilon)$ and $n$ are sufficiently large.

For part (b) of the theorem, we may assume that $\sigma\leq 2$. The same argument as before works, once we observe that, using \eqref{jm10},
\begin{equation*}
||\mathbf{v}(\Phi_n)-\mathbf{v}(\Phi_n^M)||_{F^\sigma([-5/4,5/4])}\leq C||\Phi_n-\Phi_n^M||_{\widetilde{H}^\sigma}(1+||\Phi_n||_{\widetilde{H}^\sigma}+||\Phi_n^M||_{\widetilde{H}^\sigma}).
\end{equation*}
\end{proof}

\end{document}